\documentclass[reqno]{amsart}
\usepackage{epsfig,graphicx}

\usepackage{float}
\usepackage{upgreek}
\usepackage{enumitem}
\usepackage{appendix}
\usepackage{amsmath}
\newtheorem{theorem}{Theorem}[section]

\newtheorem{lemma}[theorem]{Lemma}

\usepackage{tikz-cd}
\usepackage{tikz}
\usetikzlibrary{matrix,calc}
\usepackage{bbm}

\DeclareMathOperator{\dist}{dist}

\DeclareMathOperator{\ra}{\mathrm{A}}
\DeclareMathOperator{\rb}{\mathrm{B}}
\DeclareMathOperator{\rc}{\mathrm{C}}
\DeclareMathOperator{\hn}{\hat{\nu}}
\DeclareMathOperator{\ho}{\hat{\omega}}

\DeclareMathOperator{\ff}{\mathfrak{f}}
\DeclareMathOperator{\fa}{\mathfrak{a}}
\DeclareMathOperator{\fc}{\mathfrak{c}}
\DeclareMathOperator{\fe}{\mathfrak{e}}

\DeclareMathOperator{\fp}{\mathfrak{p}}
\DeclareMathOperator{\fm}{\mathfrak{m}}
\DeclareMathOperator{\fq}{\mathfrak{q}}
\DeclareMathOperator{\fg}{\mathfrak{g}}
\DeclareMathOperator{\fk}{\mathfrak{k}}
\DeclareMathOperator{\bk}{\mathbf{k}}

\DeclareMathOperator{\real}{\mathbb{R}}
\DeclareMathOperator{\complex}{\mathbb{C}}

\DeclareMathOperator{\integer}{\mathbb{Z}}

\DeclareMathOperator{\uk}{\mathfrak{U}_k}
\DeclareMathOperator{\ukp}{\mathfrak{U}_{k^\prime}}
\DeclareMathOperator{\ukpp}{\mathfrak{U}_{k^{\prime\prime}}}
\DeclareMathOperator{\knuk}{\bigcup_{k\in\mathbb{N}}\mathfrak{U}_k}

\DeclareMathOperator{\yb}{\overline{\mathcal{Y}}}

\DeclareMathOperator{\hi}{{\mathbb{H}}}

\DeclareMathOperator{\ce}{{\mathcal{E}}}

\DeclareMathOperator{\ub}{\overline{U}}

\DeclareMathOperator{\cf}{\mathcal{F}}

\DeclareMathOperator{\crr}{\mathcal{R}}
\DeclareMathOperator{\ipp}{\mathcal{I}_{\mathfrak{p}^\prime}}

\DeclareMathOperator{\ssq}{\mathbb{S}_{\mathfrak{q}}}

\DeclareMathOperator{\iq}{\mathcal{I}_{\mathfrak{q}}}
\DeclareMathOperator{\ip}{\mathcal{I}_{\mathfrak{p}}}
\DeclareMathOperator{\hp}{\mathbb{H}_{\mathfrak{p}}}
\DeclareMathOperator{\hpp}{\mathbb{H}_{\fp^{\prime}}}

\DeclareMathOperator{\hpe}{\mathbb{H}_{\mathfrak{p},E}}

\DeclareMathOperator{\cp}{\mathbb{CP}}
\DeclareMathOperator{\rp}{\mathbb{RP}}
\DeclareMathOperator{\rpb}{\overline{\mathbb{RP}}}
\DeclareMathOperator{\rpcb}{\overline{\mathbb{RP}}^{\mathfrak{C}}}
\DeclareMathOperator{\ssp}{\mathbb{S}_{\mathfrak{p}}}
\DeclareMathOperator{\ep}{\mathrm{E}_{\mathfrak{p}}}
\DeclareMathOperator{\eq}{\mathrm{E}_{\mathfrak{q}}}

\DeclareMathOperator{\fcc}{\mathfrak{F}}

\DeclareMathOperator{\ctwon}{\mathcal{C}_{2n}}
\DeclareMathOperator{\ctwo}{\mathcal{C}_{2}}

\DeclareMathOperator{\ctwob}{\overline{\mathcal{C}}_{2}}
\DeclareMathOperator{\ctwonb}{\overline{\mathcal{C}}_{2n}}
\DeclareMathOperator{\ctwonmt}{\mathcal{C}_{2n-2}}
\DeclareMathOperator{\ctwonmtb}{\overline{\mathcal{C}}_{2n-2}}

\DeclareMathOperator{\re}{\mathfrak{Re}}
\DeclareMathOperator{\sym}{\mathrm{Sym}}

\begin{document}


\title[Topological aspects of $\integer/2\integer$ eigenfunctions for the Laplacian on $S^2$]{Topological aspects of $\integer/2\integer$ eigenfunctions for the Laplacian on $S^2$}

\author[C. H.~Taubes]{C. H.~Taubes$^\dag$}
\author[Y. Wu]{Y. Wu$^\lozenge$}


\address{Department of Mathematics, Harvard University, Cambridge, MA}
\email{chtaubes@math.harvard.edu}

\address{Center of Mathematical Sciences and Applications, Harvard University, Cambridge, MA}
\email{ywu@cmsa.fas.harvard.edu}

\begin{abstract}
This paper concerns the behavior of the eigenfunctions and eigenvalues of the round sphere’s Laplacian acting on the space of sections of a real line bundle which is defined on the complement of an even numbers of points in $S^2$.  Of particular interest is how these eigenvalues and eigenvectors change when viewed as functions on the configuration spaces of points.           
\end{abstract}
\keywords{}

\maketitle

\section{Introduction}
	This paper describes an eigenvalue problem for the Laplacian on the round 2-sphere which is associated to a configuration of an even number of distinct points on that sphere.  Each such point configuration specifies a Hilbert space domain for the Laplacian which makes it a self-adjoint operator much like the standard Sobolev space domain (which is the domain in the case of the zero point configuration).  The eigenvalues and eigenfunctions change when the points in the configuration move about; and our goal is to understand how and why the eigenvalues and eigenfunctions change the way they do.  Our goal in this regard is realized only in part; and it is fair to say that our study produces more questions than answers.

\subsection{The background}
	To set the stage for what is to come fix an even, positive integer to be denoted by 2n and let $\fp$ denote a chosen collection of $2n$ distinct points in the round 2-sphere.  With $\fp$ in hand, let $\ip$ denote the real line bundle over $S^2-\fp$ with monodromy $-1$ on any embedded circle in $S^2-\fp$ linking any given point from $\fp$ and no others.   This is viewed as an associated vector bundle to the principle $\integer/2\integer$ bundle that is defined by the 2-1 branched cover of $S^2$ with $\fp$ being the branch locus.  (All real line bundles in this paper have implicit fiber metrics that associate them to principle $\integer/2\integer$ bundles).  Viewed in this way, the bundle $\ip$ has a canonical fiber metric which can and will be used to define the norm of sections and, with the round metric’s inner product, the norm of derivatives of sections.  With regards to derivatives:  There is a canonical, metric compatible `covariant derivative' for sections of $\ip$ and $\ip$ valued tensors on $S^2$.  This is because these are defined locally up to multiplication by the constant $-1$.  If $f$ denotes a section, then $df$ is used to denote the corresponding covariant derivative which is a section of $\ip\otimes T^*S^2$ over $S^2-\fp$.  (The covariant derivative on $\ip$-valued 1-forms and tensors is defined with the help of the round metric’s Levi-Civita connection.)

 	Let $\hp$ denote the Hilbert space completion of the space of smooth, compactly supported sections of $\ip$ over $S^2-\fp$ using the norm whose square is defined by the rule
\begin{equation}
f \longrightarrow \|f\|_{\mathbb{H}}^2 = \int_{S^2-\fp} |df|^2 + \int_{S^2} |f|^2\tag{1.1}
\end{equation}
As explained in \cite{TW}, there is a largest positive number to be denoted by $\ep$ such that the inequality
\begin{equation}
\int_{S^2-\fp}|df|^2 \geq E_{\fp} \int_{S^2} |f|^2\tag{1.2}
\end{equation}
holds for every $f \in\hp$.  (As a consequence, the right most integral on the right hand side of (1.1) is redundant except in the well understood case when $n = 0$.)

Of interest in this article are those elements in $\hp$ that are eigensections of the Laplacian:  A section $f$ is an eigensection if there is a number $\lambda\in[\ep, \infty)$ such that 
\begin{equation}
\int_{S^2- \fp} \langle df, dh\rangle = \lambda \int_{S^2-\fp} \langle f, h\rangle
\tag{1.3}
\end{equation}
for all $h \in \hp$.  The number $\lambda$ is called the eigenvalue.  (The notation here and subsequently uses $\langle , \rangle$ to denote the round metric's inner product on both $T^*S^2$ and on $T^*S^2 \otimes \ip$.)  Any given eigensection for the Laplacian is a smooth section of $\ip$.  This is because the eigensection (denoted by $f$) can be canonically viewed (up to multiplication by $-1$) as an ordinary function on any disk that is disjoint from $\fp$ where it obeys the spherical metric's Laplace equation
\begin{equation}
d^\dag df = \lambda f
\tag{1.4}
\end{equation}
with $\lambda$ denoting the corresponding eigenvalue.  (What is denoted by $d^\dag$ signifies the formal $L^2$ adjoint of $d$.)  

An eigensection from $\hp$ for the Laplacian is said in what follows to be an $\ip$-eigensection and its corresponding eigenvalue is said to be an $\ip$-eigenvalue.

Proposition 2.1 in \cite{TW} asserts that $\ep$ is an $\ip$-eigenvalue.  The technology used to prove that proposition can be employed in a straightforward way to prove the proposition that follows:\\

\noindent\bf Proposition 1.1:  \it Fix a set $\fp$ of distinct points in $S^2$. The corresponding Hilbert space $\hp$ has an orthonormal basis of $\ip$-eigensections.  Moreover, the corresponding set of $\ip$-eigenvalues have finite multiplicities and they form a discrete set in $[\ep, \infty)$ with no accumulation points.\\

\rm	Let $\ctwon$ denote the space of unordered $2n$-tuples of points in $S^2$.  Each such $2n$-tuple of points has a corresponding set of $\ip$-eigensections and the associated set of $\ip$-eigenvalues:  Our goal in this paper is explore how these $\ip$-eigensections and their $\ip$-eigenvalues change as $\fp$ moves in $\ctwon$.  

\subsection{Eigenvalue functions on $\ctwon$}

	We start by investigating the behavior of the eigenvalues as functions of the configurations.  Our principle example is the lowest eigenvalue function, the function $\fp \longrightarrow \ep$.  This function on $\ctwon$ has the following properties:\it
	\begin{itemize}
	\item The function $\mathrm{E}_{(\cdot)}$ is continuous; it is differentiable only where the corresponding eigenspace is 1-dimensional.
	\item The infimum of $\mathrm{E}_{(\cdot)}$ on $\ctwon$ is zero, it is only achieved on $\ctwon$ if $n = 0$.
	\item The supremum of $E_{\fp}$ on $\ctwon$ is greater than an $n$-independent multiple of $n$ and it is always achieved.  But the function $\mathrm{E}_{(\cdot)}$ is not differentiable there (nor is it at any other putative critical point).
\end{itemize}
\rm \hfill (1.5)\\

\rm By way of an example for the first bullet:  The eigenspace for the minimal eigenvalue $\ep$ on $\ctwo$ has dimension 2 when the two points in the configuration are antipodal (see Section 7).  And, as we explain later, there are configurations in any $n> 1$ version of $\ctwon$ where the minimal eigenvalue eigenspace has multiplicity at least four.  (This multiplicity business is surprising because mathematical physics/differential equation lore says that the lowest eigenvalue of a Laplacian must have multiplicity one.)

To elaborate on the second bullet:  The function $\mathrm{E}_{(\cdot)}$ limits to zero along paths in $\ctwon$ that bring all of the constituent points together.  We study the behavior of the lowest eigenvalue function $\mathrm{E}_{(\cdot)}$ (and also higher $\mathcal{I}_{(\cdot)}$-eigenvalues) along sequences where subsets of configuration points coalesce by introducing a natural compactification of $\ctwon$ as a stratified space with strata $\{\ctwon, \mathcal{C}_{2n-2}, \mathcal{C}_{2n-4}, \ldots, \mathcal{C}_{0}\}$.  We then show that $\mathrm{E}_{(\cdot)}$ (in fact, all $\mathcal{I}_{(\cdot)}$-eigenvalues) extends continuously to this compactification.  The compactification is denoted by $\ctwonb$.

To elaborate on the third bullet:  Although $\mathrm{E}_{(\cdot)}$ is differentiable where the corresponding eigenspace has multiplicity 1, it is not necessarily differentiable where the multiplicity is greater than 1.  In this regard, we have an expression for the directional derivative of an eigenvalue function $\fp\longrightarrow \lambda_{\fp}$ where the eigenspace has multiplicity 1  which implies this:\\

\it The directional derivatives of $\lambda_{(\cdot)}$ at a given configuration p, are all zero
 if and only if the corresponding $\ip$-eigensection is $\mathcal{O}\left(\mathrm{dist}(\fp,\cdot)^{3\over 2}\right)$ near $\fp$.
 
\rm \hfill (1.6)\\

\rm By way of comparison, the arguments in Part 3 of Section 4 in \cite{TW} prove that the most that can be said about an $\ip$-eigensection with no other a priori information is that its norm is $\mathcal{O}\left(\mathrm{dist}(p,\cdot)^{1\over 2}\right)$ near any given point $p\in\fp$.  As we explain below, the behavior in (1.6) implies that the eigenspace multiplicity of $\mathrm{E}_{(\cdot)}$ is greater than 1 at any configuration where all of the directional derivatives are zero.  That implies, in turn, that $\mathrm{E}_{(\cdot)}$ does not have continuous derivatives where it is maximal.  

What is said in (1.5) for the minimal eigenvalue function $\mathrm{E}_{(\cdot)}$ has an analog for the function that assigns to a configuration its $k$'th lowest eigenvalue (for any $k \geq 1$); see Section 3.4 for a precise definition of this function. The basic conclusion is that these eigenvalue functions are differentiable where they have multiplicity 1, and need only be continuous where their multiplicity is greater than 1.  Moreover, any putative critical point must occur where the multiplicity is in fact greater than 1 and thus where the function might not be differentiable.  (Ljusternick-Schnirrelman constructions -- which is to say min-max -- can be used to define the notion of a critical value and critical point for a function that is only continuous.)  
The preceding remarks about the lack of differentiability of the k’th eigenvalue functions are proved in Section 3.5.

Looking ahead:  Sections 2 and 3 of this paper investigate the properties of these eigenvalue functions on $\ctwon$ and $\ctwonb$.

\subsection{The projective space bundle}
	To understand the phenomena depicted in (1.5), we replace the eigenvalue functions on $\ctwon$ with a `universal' eigenvalue function on the space of pairs of the form $(\fp, f)$ where $\fp$ is from $\ctwon$ and $f$ is a section of $\ip$ with the $S^2$-integral of $f^2$ equal to 1.  This function is denoted by $\mathcal{E}$ and it is defined by the rule
\begin{equation}
(\fp, f) \longrightarrow \int_{S^2-\fp} |df|^2.
\tag{1.7}
\end{equation}
This function $\mathcal{E}$ is (formally) a smooth function of these sorts of pairs.  Moreover, a perturbation theoretic argument says this:  If $(\fp, f)$ is a critical point of $\mathcal{E}$, then $f$ is an $\ip$-eigensection (this is Raleigh-Ritz) and $f$ as it nears any point in $\fp$ vanishes as $\mathcal{O}\left(\mathrm{dist}(\fp,\cdot)^{3\over 2}\right)$ which is precisely the condition in (1.6).  

	Given $\fp \in \ctwon$, let $\ssp$ denote the subset in the Hilbert space $\hp$ (the norm is depicted in (1.1)) that consists of the elements whose square has $S^2$-integral equal to 1.  The function $\mathcal{E}$ is formally a function on a `fiber bundle' over $\ctwon$ whose fiber of $\fp$ is the subspace $\ssp$.  But now something interesting appears:  There is no such fiber bundle.  This is a very pretty example of what the physicists would call an anomaly, this one being a $\integer/2\integer$ version of a integer characteristic class anomaly that was discussed at some length in an unpublished paper by G. Segal \cite{Se} in the context of certain observations by L. Fadeev \cite{Fa}.  (There are many papers on anomalies in physics; see the recent paper of L. Müller for other versions of anomalies and references \cite{Mu}.)  In the context at hand, there exists an obstruction class in $H^2(\ctwon;\integer/2\integer$) that is zero if and only if the assignments $\{\fp\longrightarrow\ssp\}_{\fp\in\ctwon}$   define a fiber bundle.  This obstruction class is non-zero; it is denoted by $\omega$.

	We can none-the-less define an $\rp^\infty$ bundle over $\ctwon$ whose fiber over any given configuration $\fp$ is $\ssp/\{\pm 1\}$ which is sufficient for the purposes at hand because $\mathcal{E}$ has the same value on $(\fp, f)$ as it has on $(\fp, -f)$.  This bundle is denoted by $\rp$ and $\mathcal{E}$ descends there as a smooth function.  In particular, one can look for critical points of $\mathcal{E}$ on $\rp$.  This leads us to study the (co)homology of $\rp$ because that (in principle) dictates critical points of $\mathcal{E}$.  (This $\rp$ bundle and its function $\mathcal{E}$ are described in Section 4 of this paper.)
	
By way of a parenthetical remark:  In an alternate universe where the obstruction class is zero, one should still consider $\mathcal{E}$ on the quotient $\rp^\infty$ to study its critical points. 

\subsection{Critical values of $\mathcal{E}$}
	The observation in the third bullet of (1.5) that the minimal eigenvalue function $\mathrm{E}_{(\cdot)}$ is not everywhere differentiable is an automatic consequence of the following observation about the function $\mathcal{E}$ on $\rp$:\\

\it Suppose that $n$ > 0 and that $(\fp,[f])$ is a critical point of the $\ctwon$ version of $\mathcal{E}$.  Then the
 corresponding critical value is not the lowest $\ip$-eigenvalue (which is $\ep$);  
in fact, it is greater than the n'th $\ip$-eigenvalue.

\rm \hfill (1.8)\\

\rm As noted in the first bullet of (1.5), the lack of differentiability of the minimal eigenvalue function $\mathrm{E}_{(\cdot)}$ is due to the existence of configurations where the minimal eigenvalue has an eigenspace with dimension greater than 1.  The condition in (1.8) tells us slightly more; it tells us that there can't be a differentiable section of $\rp$ (as a bundle over $\ctwon$) that assigns any given configuration to some $\mathcal{I}_{(\cdot)}$-eigensection (modulo $\pm 1$) with minimal eigenvalue.  As explained in the next paragraph, the lack of a section of $\rp$ is dictated by the existence of a non-zero anomaly class.  

With regards to the anomaly and the existence of sections:  The key fact about the anomaly class is that it's pull-back to $\rp$ via the projection map is zero; and this can't be the case if there is a section because the composition of first the section and then the projection is the identity map on $\ctwon$.  (If the anomaly were zero, then there would be a section of the corresponding infinite dimensional sphere bundle; and thus a section of the $\rp^\infty$ bundle.  This is because a vector bundle over a finite dimensional manifold whose fiber dimension is greater than the base manifold dimension has a non-zero section.)

Morse theoretic arguments (min-max, for example) using the (co)homology of $\rp$ would yield critical points of $\mathcal{E}$ were $\mathcal{E}$ a proper function on $\rp$.  Although $\mathcal{E}$ is proper along the fibers of the projection to $\ctwon$, it is not globally proper because $\ctwon$ isn't compact. 

\subsection{Homology classes of $\rpb$}
	As noted previously, the space $\ctwon$ has a compactification as a stratified space (denoted by $\ctwonb$) that accounts for the coalescence of subsets of configuration points.  As it turns out, there is an extension of $\rp$ over this compactification (denoted by $\rpb$) that maps to the compactification $\ctwonb$ and a corresponding, continuous extension of the function $\mathcal{E}$ as a function over $\rpb$.   This $\rpb$ space has non-zero $\integer/2\integer$ homology classes in degrees $4n+2m+1$ with $m$ any non-negative integer.  It is possible that these classes yield critical points of $\mathcal{E}$ on large $n$ versions of $\ctwonb$ because each such $(n,m)$ class is the image via a fiber preserving map to $\rpb$ of the $\integer/2\integer$ fundamental class of an $\rp^{2m+1}$ bundle over $\ctwonb$.  (The class sees the whole of $\ctwonb$ because of this; in particular, the class is not carried by any single strata in $\ctwonb$.)  We show that these classes do give critical points of $\mathcal{E}$ on $\rpb$ and we show in Section 6 of this paper that the corresponding critical values (as functions of either $n$ or $m$) are unbounded.  Even so, it remains to be seen whether these critical points all lie in some fixed $\mathcal{C}_{2j}$ strata with $j$ having an $n$ and $m$ independent upper bound.  This is work for the future.	

By way of a remark:  The classes used here can be viewed as a homological work-around for the non-differentiability of the $k$'th eigenvalue function on $\ctwon$ in the case when $k$ is odd.  (There are also classes in the $\integer/2\integer$ homology of $\ctwonb$ with even degree --  any even number -- but their corresponding critical points are almost surely in the $\mathcal{C}_0$ or $\ctwo$ strata.  These strata and the corresponding critical points are described in Section 7.)

\subsection{A potential application}
 	Critical points of $\mathcal{E}$ on $\ctwon$ (should they exist) can be used to construct homogeneous singularity models for $\integer/2$ harmonic 1-forms and spinors on $\mathbb{R}^3$.  To explain (see \cite{TW} and/or \cite{T} for more), a homogeneous $\integer/2$ harmonic 1-form on $\mathbb{R}^3$ consists of a data set $(Z, \mathcal{I}, \nu)$ with $Z$ being a finite union of rays from the origin, with $\mathcal{I}$ being real line bundle over the complement of $Z$ and with $\nu$ being an $\mathcal{I}$-valued 1-form defined on $\mathbb{R}^3-Z$ that has the following properties:

\it

\begin{itemize}
\item $\nu$ is both closed and coclosed.
\item $|\nu|$ extends over $Z$ as a H\"older continuous function on $\mathbb{R}^3$.
\item $\nu$ is homogeneous with respect to pull-back by the rescaling diffeomorphism that multiplies the Euclidean coordinates by a non-zero real number.
\end{itemize}
\rm \hfill (1.9)\\

\rm As explained in \cite{TW}, the data sets of this sort with $Z$ being the union of $2n$ rays are in 1-1 correspondence with pairs $(\fp, f)$ with $\fp \in \ctwon$ and with $f\in\hp$ being a non-zero eigensection of the Laplacian with $|f| \leq \mathcal{O}\left(\mathrm{dist}(\fp,\cdot)^{3\over 2}\right)$  near $\fp$ (which is exactly the critical point criteria for the function $\mathcal{E}$).  This same \cite{TW} explains how such a pair can also be used to define a homogeneous, $\integer/2$ harmonic spinor on $\mathbb{R}^3$.  To summarize the correspondence:  The set $Z$ is the union of the rays from the origin through the points in the configuration $\fp$, the line bundle $\mathcal{I}$ is the pull-back of $\ip$ via the map $\pi: x\longrightarrow {x\over |x|}$ from $\mathbb{R}^3\setminus\{0\}$ to $S^2$, and $\nu$ is given by the rule
\begin{equation}
\nu = d\left(
|x|^\mu \pi^* f
\right),
\tag{1.10}
\end{equation}
where $\mu$  is  ${1\over 2}\left(1+(1+4\lambda)^{1\over 2} \right)$ with $\lambda$ denoting the corresponding $\ip$-eigenvalue.

Before \cite{TW}, the only example of a $\integer/2$ harmonic 1-form on $\mathbb{R}^3$ had $Z$ being a straight line through the origin (which corresponds to a critical point of $\mathcal{E}$ on $\mathcal{C}_2$).  The paper \cite{TW} used the symmetry groups of the Platonic solids to construct examples where $Z$ is the union of (respectively) four rays from the origin (pointing to the vertices of a regular tetrahedron), eight rays from the origin (pointing to the vertices of a cube) and 12 rays from the origin (pointing to the vertices of an iscosahedron.)  Going in the reverse direction, the constructions in \cite{TW} give critical points of $\mathcal{E}$ on $\mathcal{C}_4$ and $\mathcal{C}_8$ and $\mathcal{C}_{12}$.  

\subsection{Two remarks for the reader and a table of contents}

	The first remark concerns the case of $n = 1$:  This case is illuminating because much can be seen explicitly in this case.  The $n = 1$ case is discussed in some detail in Section 7 and a reader is advised to refer to that section when new notions are introduced.  
	
The second remark concerns two notational conventions that are used throughout this paper:  The first convention has $c_0$ denoting a number greater than 1 whose value is independent of relevant background such as a configuration in $\ctwon$ or the integer $n$.  What precisely it is independent of will be clear from the context.  The value of $c_0$ increases between consecutive incarnations.  The second convention has $\chi$ denoting a chosen, smooth and non-increasing function on $\mathbb{R}$ that is equal to 1 on $\left(-\infty, {1\over 4} \right]$ and equal to 0 on $\left[{3\over 4} , \infty\right)$.  This fixed version of $\chi$ is implicitly used to construct `bump functions' and `cut-off functions'.

	What follows directly is the table of contents for this paper.

\begin{enumerate}

\item[1.]   \textsc{Introduction}
\item[2.]   \textsc{Basic technology}
\item[3.]  \textsc{The lowest eigenvalue as a function on $\ctwon$}
\item[4.]  \textsc{The $\rp^\infty$ bundle}
\item[5.]  \textsc{The extension of $\rp$ and $\mathcal{E}$ to $\ctwonb$}
\item[6.]  \textsc{Min-max for $\mathcal{E}$ on $\rpb$}
\item[7.]  \textsc{The case of $\mathcal{C}_2$}
\item[a.]  \textsc{Appendix on the cohomology of a weak $\rp^\infty$ bundle}

\end{enumerate}

\subsection{Acknowledgements}
Before starting, Y-Y.~Wu wishes to acknowledge and thank the generous support of the Center of Mathematical Sciences and Applications (CMSA) at Harvard.  Meanwhile, C. ~H.~Taubes wishes to acknowledge and thank the generous support of the National Science Foundation (DMS grant number 2002771).

\section{Basic technology}

This section introduces the basic technology that is used in later sections.  To set the stage for what is to come, suppose henceforth that $\fp$ is a chosen $2n$-tuple from $\ctwon$.  
The convention in what follows has $c_0$ denoting a number that is greater than 1 and that is independent of any chosen sections of $\ip$ in any given appearance.  It can depend on $\fp$ unless stated otherwise.  As noted in Section 1, the value of $c_0$ can be assumed to increase between successive appearances. 

\subsection{$\fp$-independent a priori bounds}
	The lemma that follows summarizes some basic a priori bounds on eigensections of the Laplacian for a given $\fp \in \ctwon$.  These bounds are independent of $\fp$.\\

\noindent \bf Lemma 2.1:  \it There exists $\kappa > 1$ with the following significance:  Fix $\fp \in\ctwon$ and let $f$ denote an $\ip$-eigensection.  Denote the corresponding $\ip$-eigenvalue by $\lambda$.  Then the function $|f|$  obeys:
\begin{itemize}
\item	$\displaystyle \int_{S^2} |d|f||^2 = \lambda \int_{S^2} |f|^2.$
\item $\displaystyle |f| \leq \kappa (\lambda+1) \left( \int_{S^2} |f|^2 \right)^{1\over 2}$ on the complement of $\fp$.
\end{itemize}

\noindent \bf Proof of Lemma 2.1:  \rm The first bullet inequality follows because any $f \in \hp$ obeys the inequality $|d|f| | \leq |df|$ almost everywhere, and because the integral of $|df|^2$ for the eigensection $f$ is equal to $\lambda$ times that of $|f|^2$. To prove the second bullet, note first the function $|f|$ obeys the distributional inequality
\begin{equation}
d^\dag d|f| \leq \lambda|f|.\tag{2.1}
\end{equation}
With this understood, fix a point $q\in S^2-\fp$ and let $G_q(\cdot)$ denote the Green's function for the operator $d^\dag d + 1$ with pole at $q$.   This Green's function is smooth on $S^2-q$ and it is bounded near $q$ by an $\mathcal{O}(1)$ multiple of $-{1\over 2\pi} \ln(\mathrm{dist}(\cdot,q))$. The plan for what follows is to multiply both sides of (2.1) by $G_q$, integrate the resulting inequality over $S^2$ and then integrate by parts twice to bound $|f|(q)$ by $\lambda$ times the integral of $G_q|f|$.  The latter bound will lead to the second bullet's bound via the Cauchy-Schwarz inequality because the $S^2$ integral of $G_q^2$ is finite.  To implement this strategy, care must be taken near the points of $\fp$.  To that end, fix $\epsilon>0$ but much less than 1, and then let $\chi_\epsilon(\cdot)$ denote the function on $S^2$ that is given by the rule
\begin{equation}
\chi_\epsilon(\cdot) = \chi\left(
2{ \ln(\mathrm{dist}(\cdot, \fp)\over \ln(100\epsilon)} - 1
\right).
\tag{2.2}
\end{equation}

This function is equal to 1 where the distance to $\fp$ is greater than $(100\epsilon)^{1\over 2}$ and equal to zero where the distance to $\fp$ is less than $100\epsilon$.   If $\epsilon$ is much less than the distance from $q$ to $\fp$, then $q$ will sit where $\chi_\epsilon$ is equal to 1.  Note also that the derivative of $\chi_\epsilon$ has support only in the region where the distance to $\fp$ is between $100\epsilon$ and $(100\epsilon)^{1\over 2}$ where it is bounded by $c_0 {1\over |\ln\epsilon|} {1\over {\dist(\cdot, \fp)}}$ with $c_0$ being independent of $\fp$.  This implies in particular that the $S^2$ integral of $|d \chi_\epsilon|^2$ is bounded by $c_0{1 \over |\ln\epsilon|}$  with $c_0$ again independent of $\fp$.  With these last points understood, multiply both sides of (2.1) by $\chi_\epsilon G_q$ with $\epsilon$ very small and then integrate by parts.  The result of doing that leads to the inequality
\begin{equation}
|f| (q) \leq c_0 \lambda \left(
\int_{S^2} |f|^2
\right)^{1\over 2} + c_{\fp, q} {1\over \sqrt{|\ln \epsilon|}} \left(
\int_{S^2} |df|^2 + \int_{S^2} |f|^2
\right)^{1\over 2}
\tag{2.3}
\end{equation}
where $c_0$ is independent of $\fp$ and $q$ and $f$ and $\epsilon$ whereas $c_{\fp,q}$ is only independent of $f$ and $\epsilon$ (it does depend on the distance from $q$ to $\fp$).   Take $\epsilon\to0$ in (2.3) to obtain the second bullet's bound for $|f|$.

\subsection{Asymptotics near points in $\fp$}
	As noted in the proof of Lemma 2.1, if $f$ is from $\hp$, then its norm is in the standard $L^2_1$ Sobolev space for the whole of $S^2$; and this implies via Sobolev inequalities that $|f|^q$ is integrable on $S^2$ for any finite, non-negative value of $q$.  The following lemma concerns specifically $f$'s behavior near the points in $\fp$; it says in effect that a function in $\hp$ must vanish in a weak sense as a point in $\fp$ is approached.  \\

\noindent \bf Lemma 2.2:  \it Given $\fp \in \ctwon$, there exists $\kappa > 1$ with the following significance:  If $f \in \hp$, then
$\displaystyle \int_{S^2-\fp} {1\over \dist(\cdot, \fp)^2} |f|^2 \leq \kappa \int_{S^2-\fp} |df|^2.$\\

\noindent \bf Proof of Lemma 2.2:  \rm The asserted bound follows directly from Lemma 3.2 in \cite{TW} by invoking it for the annuli from the set $\{\mathcal{A}_n(p): p \in \fp \text{ and } n \in \mathbb{N}\}$ where any given $\mathcal{A}_n(p)$ is defined by the rule
\begin{equation}
\mathcal{A}_n(p) \equiv \{q \in S^2-\fp: 2^{-n-1}\delta < \dist(q,p) < 2^{-n}\delta\}_{n\in\mathbb{N}} 
\tag{2.4}
\end{equation}
with $\delta$ being ${1\over 100}$  times the minimum distance between the points in $\fp$.

More can be said about behavior near the points in $\fp$ when $f$ is an $\ip$-eigensection.  In this case, the necessary observations about the behavior of an $\ip$-eigensection near a point $p \in \fp$ invoke observations from Parts 2 and 3 of Section 4 in \cite{TW}, specifically Equations (4.2)-(4.5) in \cite{TW} which hold for a Laplace eigensection whether or not $\fp$ is invariant under a symmetry group action.  These observations are summarized by the upcoming Lemma 2.3.  This summary of what is said in \cite{TW} uses stereographic projection from $p$'s antipodal twin to define a complex coordinate to be denoted by $z$ near $\fp$ with both $z = 0$ and $|dz| = 1$ at $p$.  By way of notation:  The lemma use $\re (\cdot)$ to denote the real part of the indicated $\complex$-valued function or differential form. \\

\noindent \bf Lemma 2.3:  \it Let $f$ denote an $\ip$-eigensection.  Given any point $p \in\fp$, there exists the following data:
\begin{itemize}
\item A non-zero complex number to be denoted by $\mathfrak{a}$,
\item A non-negative integer to be denoted by $n$.
\item A section of $\ip$ defined near $p$ to be denoted $\mathfrak{e}$ obeying$$\lim_{|z|\to0} (|\mathfrak{e}| + |z||d\mathfrak{e}|)|z|^{-n-{1\over 2}} = 0.$$ 
\end{itemize}
These are such when $f$ is depicted using the complex coordinate $z$ and its complex conjugate, it has the form $f|_z = \re\left(\mathfrak{a}z^{n+{1\over 2}}\right) +\mathfrak{ e}$.\\

\rm

When the point $p$ is germane to the discussion, then the corresponding versions of $\fa$ and $n$ are written as $\fa_p$ and $n_p$.  When $f$ and $\fp$ are germane, they are written as $\fa_p(f)$ and $n_p(f)$.  

By way of a definition:  The subset $\fp_{f} \in \fp$ is defined to be the set of points with the integer $n$ being zero.

\subsection{The divergence identity}
	The identity that follows momentarily uses $\fg$ to denote the round metric on $S^2$ when viewed as a section of $\sym^2(T^*S^2)$.  Supposing that $f$ denotes an $\ip$-eigensection and $\lambda$ its eigenvalue, define the section $\mathbb{T}$ of $\sym^2(T^* S^2)$ on $S^2-\fp$ by the rule
\begin{equation}
\mathbb{T} = df \otimes df - {1\over 2} (|df|^2 - \lambda f^2)\fg.
\tag{2.5}
\end{equation}
By virtue of (1.4), this section is divergence free, thus
\begin{equation}
\nabla^\dag \mathbb{T} = 0.
\tag{2.6}
\end{equation}
	An integrated version of this identity will be presented momentarily in the upcoming Lemma 2.4.  An observation is needed first which is this:  Let $s$ denote a continuous section of $\sym^2(T^*S^2)$ that is defined on the whole of $S^2$.   Then the function $\langle s,\mathbb{T}\rangle$ has finite integral over $S^2$ and the assignment 
\begin{equation}
s \longrightarrow \int_{S^2} \langle s, \mathbb{T}\rangle
\tag{2.7}
\end{equation}
defines a bounded linear functional on the Banach space of continuous sections of $\sym^2(T^*S^2)$.  This follows because $\mathbb{T}$'s norm is bounded by $c_0(|df|^2 + \lambda f^2)$ whose integral over $S^2$ is bounded by $c_0$ times the $\mathbb{H}$-norm of $f$.  (It also follows from Lemma 2.3 that $\mathbb{T}$'s singularities at the points in $\fp$ are integrable.)

	The upcoming lemma also uses the convention whereby vector fields and 1-forms are identified using the metric $\fg$.  \\
 
\noindent \bf Lemma 2.4:  \it Let $f$ denote a given $\ip$-eigensection and let $\lambda$ denote the corresponding eigenvalue.  Supposing that $\nu$ is a smooth vector field on $S^2$, then
\begin{equation*}
\int_{S^2}\langle \nabla \nu, \mathbb{T} \rangle
=
-{\pi \over 4} \sum_{p\in\fp_{f}} \langle \nu|_p, \re(\fa_p^2dz)\rangle.
\end{equation*}

\noindent \bf Proof of Lemma 2.4: \rm Fix $\epsilon>0$ but much less than the distance between any two points in $\fp$.  Supposing that $p \in\fp$, use $D_\epsilon(p)$ to denote the disk centered at $p$ with radius $\epsilon$. Integrate the function $\langle \nu, \nabla^\dag \mathbb{T}\rangle$ over the complement of the union of these radius $\epsilon$ disks.  This integral is zero because the integrand is zero.  Even so, integrate by parts to write that integral (hence 0) as
\begin{equation}
\int_{S^2 - \cup_{p\in\fp}D_\epsilon(\fp)} \langle \nabla \nu, \mathbb{T}\rangle + 
\sum_{p \in \fp_{f}\partial D_\epsilon(p)} \int \langle \hat r \otimes \nu, \mathbb{T}\rangle
\tag{2.8}
\end{equation}
with $\hat r$  denoting here the outward pointing, length 1 normal vector to the relevant $\partial D_\epsilon (p)$.  The $\epsilon\to0$ limit of the left hand integral in (2.8) is the integral of $\langle\nabla \nu, \mathbb{T}\rangle$ over $S^2$.  Use Lemma 2.3 to evaluate the $\epsilon\to0$  limit of the right integral in (2.8).  The assertion that the sum of these two $\epsilon\to0$  limits is zero is Lemma 2.4. 

	The divergence identity in (2.6) for $\mathbb{T}$ as defined in (2.5) has a generalization that is also useful:  To state the more general identity, suppose now that $f$ and $f^\prime$ are two $\ip$-eigensections and let $\lambda$ and $\lambda^\prime$ denote their respective eigenvalues.  Now define the symmetric tensor $\mathbb{S}$ by the rule 
\begin{equation}
\mathbb{S} = df \otimes df^\prime + df^\prime\otimes df - \left(\langle df,df^\prime\rangle -  {1\over 2}(\lambda + \lambda^\prime )ff^\prime\right)\fg.
\tag{2.9}
\end{equation}
The tensor $\mathbb{S}$ is divergence free if $\lambda = \lambda^\prime$, but otherwise
\begin{equation}
\nabla^\dag\mathbb{S} = {1\over 2} (\lambda^\prime - \lambda) (f df^\prime - f^\prime df).
\tag{2.10}
\end{equation}
Much the same argument that led to Lemma 2.4 leads in this case to the following analog of (2.8):
\begin{equation}
\int_{S^2}\langle \nabla \nu, \mathbb{S} \rangle
- {1\over 2} (\lambda^\prime - \lambda) \int_{S^2} \langle \nu, f df^\prime - f^\prime df\rangle + 
\sum_{p \in \fp_{f}} \lim_{\epsilon \to 0} \int_{\partial D_\epsilon(p)} \langle \hat r \otimes \nu, \mathbb{S} \rangle = 0.
\tag{2.11}
\end{equation}
Now we can invoke the $f$ and $f^\prime$ versions of Lemma 2.3 to evaluate the limit in (2.11) and doing so leads to the following identity:
\begin{equation}
\int_{S^2}\langle \nabla \nu, \mathbb{S} \rangle
- {1\over 2} (\lambda^\prime - \lambda) \int_{S^2} \langle \nu, f df^\prime - f^\prime df\rangle
=
-{\pi \over 2} 
\sum_{p \in \fp_{f} \cap \fp_{f^\prime}}  \langle \nu|_p, \re (\fa_p(f)\fa_p(f^\prime)dz) \rangle.
\tag{2.12}
\end{equation}
This is the promised generalization of Lemma 2.4's identity.

\subsection{Changing the metric}

	If $\fp$ is a given $2n$-tuple of points in $S^2$ and $\fp^\prime$ is some other $2n$-tuple, then there are diffeomorphism of $S^2$ that map the set $\fp$ to $\fp^\prime$.  If $\phi$ is a diffeomorphism that does this, then $\phi^*\ipp$ is isomorphic to $\ip$.  As a consequence, the eigenvalues for the round metric's Laplacian on sections of $\ipp$ are precisely those of the Laplacian on sections of $\ip$ that is defined by the $\phi$-pull-back of the round metric (which won't be the round metric unless $\phi$ is a rotation and/or reflection).  By virtue of this pull-back correspondence, questions about the behavior of $\ip$-eigenvalues and $\ip$-eigenfunctions of the round metric's Laplacian on the $\ctwon$-dependent family of domains are questions about the metric dependence of the eigenvalues and eigenfunctions for a $\ctwon$-dependent family of Laplace operators on a fixed domain. The latter, alternate point of view is often taken in subsequent sections and it motivates the digression that follows directly about metric Laplacians acting on the sections of the fixed bundle $\ip$.

	Any given smooth Riemannian metric on $S^2$ has a corresponding version of the Hilbert space $\hp$ and a corresponding set of eigensections for its Laplacian.  To say more about this, let $\fm$ denote a given Riemannian metric.  For convenience (this is not strictly necessary), we assume that the associated area 2-form is the same as that of the standard, round metric.  This assumption is implicit in what follows.  The metric $\fm$'s inner product on $T^*S^2$ is denoted by $\langle , \rangle_{\fm}$ and the associated norm by $|\cdot|_{\fm}$.  The corresponding Hilbert space is the completion of the space of smooth, compactly supported sections of $\ip$ using the version of (1.1) with $|df|^2$ replaced by $|df|_{\fm}^2$.  An element in this Hilbert space (denoted by $f$) is an eigensection of the metric m's Laplacian with eigenvalue $\lambda$ when  
\begin{equation}
\int_{S^2 - \fp}\langle df, dh \rangle_{\fm} = \lambda \int_{S^2-\fp^2}f h
\tag{2.13}
\end{equation}
for all elements $h$ from the Hilbert space.  Eigensections are the elements in the Hilbert space that obey the Laplace equation for $\fm$'s  version of the Laplace operator. 

	With regards to this metric $\fm$ Hilbert space:  It is the original round metric's version of $\hp$ with a different but equivalent inner product.  This is to say that any section in one is in the other and vice versa.  This understood, the metric $\fm$'s  version of the eigenvalue/eigensection question will be viewed henceforth as a question concerning a bilinear form on the original $\hp$ Hilbert space.  As was the case with the round metric $\fg$, there is a complete, orthonormal basis for $\hp$ of eigensections for the bilinear form in (2.13) with the corresponding eigenvalue set being discrete with no accumulation points and with each eigenvalue having finite multiplicity.
		
	Now suppose that $\mathcal{M}$ is a real analytic, finite dimensional manifold that parameterizes in a real analytic fashion a family of metrics on $S^2$ subject to the constraint that their area forms are the same as that of $\fg$.  There is a corresponding family of bilinear forms on $\hp$ and one can ask how the sets of eigenvalues and eigenfunctions change as the domain $\mathcal{M}$ is changed.  Exactly this sort of question is addressed in Chapters 6 and 7 of Tosio Kato's classic book \it Perturbation theory for linear operators \rm \cite{Ka}.  The rest of this section summarizes the
implications 
of Kato's discussion.

	To start the story, fix a metric $\fm_0 \in\mathcal{M}$ and then fix a smooth coordinate chart for a neighborhood of $\fm_0$ that identifies the chosen neighborhood with a ball about the origin in $\real^N$.  (The integer $N$ is the dimension of $\mathcal{M}$.)  The identification is chosen so that $\fm_0$ corresponds to the origin in $\real^N$.  A point in this ball is denoted by $y$ and the corresponding metric by $\fm_y$.  If $y$ is close to the origin in $\real^N$, then $\fm_y^{-1}$ can be written as
\begin{equation}
\fm_y^{-1} = \fm_0^{-1}(1 + \fk y + \fe_y)
\tag{2.14}
\end{equation}
with $\fk$ denoting a linear map from $\real^N$ to $C^\infty(S^2;\sym^2(TS^2))$ and with $\fe_y$ being a map from the coordinate ball in $\real^N$ to $C^\infty(S^2; \sym^2(TS^2))$ with $C^k$ norm (for any given non-negative integer $k$) bounded by a $k$-dependent multiple of $|y|^2$.  \
Now suppose that $\lambda_0$ is an eigenvalue for the $\fm_0$ Laplacian on $\hp$.  As noted in \cite{Ka}, there will be an eigenvalue of the $\fm_y$ Laplacian near to $\lambda_0$
 if $y$ is sufficiently close to the origin.  Moreover, if $\lambda_0$ has multiplicity $n$, then there will be precisely $n$ linearly independent eigensections for the $\fm_y$-Laplacian with eigenvalue very close to $\lambda_0$ if $y$ is small (much closer to $\lambda_0$ than to any other $\fm_0$ eigenvalue.)  In particular, the following is true:  Given $\epsilon>0$, there exists $\delta > 0$ 
such that the eigenvalues of these $n$ eigenvectors of the $\fm_y$ Laplacian will differ from $\lambda_0$ by at most $\epsilon$ if $|y| < \delta$.  In addition, if $|y| < \delta$, then the $n$-dimensional vector space spanned by these $n$ eigensections for the $\fm_y$-Laplacian will have $\mathbb{H}_p$-distance at most $\epsilon$ from the $\lambda_0$ eigenspace of the $\fm_0$ Laplacian.  (This distance is measured by the $\mathbb{H}_p$-orthogonal projection.)   The subsequent paragraphs and the upcoming Proposition 2.5 say more about the relationship between respective $\fm_y$ and $\fm_0$ Laplacian eigenvalues and eigenspaces.

The simplest case to consider is that where $\lambda_0$ has multiplicity 1 in which case there is one nearby eigenvalue when $y$ is near 0 and that eigenvalue varies smoothly with $y$ on a neighborhood of 0 with value $\lambda_0$ at $y = 0$.  This $y$-dependent eigenvalue can be written as $\lambda_0 + y \cdot\lambda^\prime  + \mathcal{O}(|y|^2)$ with $\lambda^\prime$ given by the rule
\begin{equation}
\lambda^\prime = \int_{S^2} \langle \fk, df_0\otimes df_0\rangle
\tag{2.15}
\end{equation}

In this case, there is one corresponding eigensection of the $\fm_y$ Laplacian that is very near to a normalized $\lambda_0$ eigensection when $y$ has small norm and it also varies smoothly with $y$ near the origin.  (This $\lambda_0$ eigensection is denoted by $f_0$; it is normalized so that the $S^2$ integral of its square is equal to 1.)  In particular, this $y$-dependent element in $\hp$ can written as a Taylor's expansion as $f_0 + y\cdot \ff  + r_y$ with $\|r_y\|_{\mathbb{H}}$ bounded by $c_0|y|^2$ and with the linear map $\ff$ from $\real^N$ to $\hp$ determined via the following two rules:\\

\begin{enumerate}
\item $\int_{S^2} \ff f_0 = 0,$ 
\item $\int_{S^2} \langle d\ff, dh\rangle_{\fm_0} - \lambda_0 \int_{S^2} \ff h + \int_{S^2} \langle \fk, df_0\otimes dh\rangle = 0$ for all $h\in\mathbb{H}$ that obey $\int_{S^2} h f_0=0.$
\end{enumerate}
\hfill (2.16)\\

With this simple case understood, consider now the case when $\lambda_0$ has multiplicity $n \geq 1$.  To set the stage for this case, introduce by way of notation $\mathbb{V}$ to denote the corresponding $n$-dimensional vector space of eigensections of the Laplacian in $\hp$ with eigenvalue $\lambda_0$.  Given a pair $f$ and $f^\prime$ from $V$, define a vector in $\real^N$ which is denoted by $\mathbf{k}(f, f^\prime)$ by the rule
\begin{equation}
\mathbf{k}(f, f^\prime) \equiv \int_{S^2} \langle \fk, df\otimes df^\prime\rangle.
\tag{2.17}
\end{equation}
Supposing that $y\in \real^N$, then the assignment of any given pair of elements in $\mathbb{V}$ to the value of $y\cdot\mathbf{k}$ on that pair defines a symmetric bilinear form on $\mathbb{V}$.  An element $f\in \mathbb{V}$ is said to be an eigensection for this bilinear form when there exists a number $e(y)$ such that
\begin{equation}
y \cdot \mathbf{k}(f, f^\prime) = e(y)\int_{S^2}f f^\prime   \text{ for all } f^\prime \in \mathbb{V}.
\tag{2.18}
\end{equation}

With regards to the $y$-dependence:  If $f$ is an eigensection for $y\cdot \mathbf{k}$, then it is likewise for $r y\cdot \mathbf{k}$ for any given real number $r$, the corresponding eigenvalue $e(ry)$ being $r e(y)$.\\

\noindent \bf Proposition 2.5:  \it  Let $\lambda_0$ denote an eigenvalue of the $\fm_0$-metric Laplacian on $\hp$ and let $\mathbb{V}$ denote the corresponding eigenspace.   Given $\epsilon>0$, there exists a positive number $\delta$ with the following significance:  If $y \in\real^N$ has norm less than $\delta$, then the span of the eigensections of the $\fm_y$ metric Laplacian on $\hp$ with eigenvalue within $\epsilon$ of $\lambda_0$
 has dimension $\dim(\mathbb{V})$.  Moreover, there is a 1-1 correspondence between these $\fm_y$-metric Laplacian eigensections and the eigensections of $y\cdot \mathbf{k}$ on $\mathbb{V}$ in the following sense:  If $f$ is any one of these $\dim(\mathbb{V})$ eigensections of the metric $\fm$ Laplacian, then this eigensection $f$ can be written as 
$$
f = f_0 + y\cdot \ff + \mathcal{O}(|y|^2) 
$$
where $f_0$  signifies an eigensection for $y \cdot \mathbf{k}$ on $\mathbb{V}$; and where $\ff$ signifies a homomorphism from $\real^N$ to $\hp$ whose image is $L^2$-orthogonal to $\mathbb{V}$ and which obeys 
$$ \int_{S^2} \langle d\ff, dh \rangle_{\fm_0} - \lambda_0 \int_{S^2} \ff h + \int_{S^2} \langle \fk, df_0\otimes dh\rangle  = 0  $$
whenever $h \in \hp$ is $L^2$-orthogonal to $\mathbb{V}$.  Furthermore, the corresponding $\fm$-metric Laplace eigenvalue of the eigensection $f$ has the form $\lambda_0+e(y) + \mathcal{O}(|y|^2)$ with $e(y)$ denoting the $y \cdot \mathbf{k}$ eigenvalue of $f_0$.\\

\normalfont
As noted, this proposition summarizes 
the consequences of various
observations from Chapters 6 and 7 of \cite{Ka}.
An important point to keep in mind with regards to this proposition:  The correspondence between the eigensections of the $\fm_y$ Laplacian with eigenvalue near $\lambda_0$ and those of the bilinear form $y \cdot \mathbf{k}$ is one-to-one when $|y|$ is small if the latter has distinct eigenvalues.  If it doesn't, then there might only be a one-to-one correspondence with a particular basis of its eigensections. 

\subsection{Pulling back by a diffeomorphism}
	 Fix a positive integer to be denoted by $N$ and then fix a linear map from $\real^N$ into the subspace of divergence zero elements in $C^\infty(S^2; TS^2)$; the map is denoted by $\nu$ in what follows.  Thus, if $y \in \real^N$, then $y\cdot \nu$ is a divergence zero vector field on $S^2$.  Taking the time 1 point on the integral curves of these vector fields defines a smooth map from a ball about the origin in $\real^N$ to the subspace of area preserving diffeomorphisms of $S^2$ which is equal to the identity at $y = 0$ and whose differential at $y = 0$ is given by $\nu$.  (The differential at $y = 0$ is a linear map from $\real^N$ into the tangent space to the space of area preserving diffeomorphisms of $S^2$ at the identity element which is the space of divergence zero vector fields on $S^2$.)    Let $\mathbf{B} \subset\real^N$ denote this ball and, supposing that $y\in \mathbf{B}$, let $\phi_y$ denote the diffeomorphism that is defined by $y$.  Also:  Let $\fm_y$ denote the pull-back of the round metric on $S^2$ by $\phi_y$.  Because $\phi_y$ is area preserving, the area form for the metric $\fm_y$  is the round metric's area form.  
	 
By taking the radius of $\mathbf{B}$ to be smaller if necessary, each $y \in \mathbf{B}$ version of the inverse metric to $\fm_y$ can be written as in (2.14) with $\fk$ as follows:
\begin{equation}
\fk = -(\nabla v + (\nabla v)^T)
\tag{2.19}
\end{equation}
with it understood that the round metric has been used to identify $TS^2$ with $T^*S^2$ and to define the transpose that is used in (2.19) to make a symmetric tensor from $\nabla \nu$.  

Now suppose that $\lambda_0$ is an eigenvalue of the round metric's Laplacian on $\hp$ and let $\mathbb{V}$ denote the corresponding eigenspace.  Then the $\real^N$-dependent bilinear form on $\mathbb{V}$ that is depicted in (2.17) can be written as
\begin{equation}
\bk(f, f^\prime) = -\int_{S^2} \langle \nabla \nu, \mathbb{S}\rangle
\tag{2.20}
\end{equation}
with $\mathbb{S}$ given by (2.9).  (The contributions to $\langle\nabla v, \mathbb{S}\rangle$ from the terms proportional to $\fg$ are zero because $\nu$ is divergence free.)  This being the case, $\mathbf{k}$ can be written using (2.12):	
\begin{equation}
\mathbf{k}(f, f^\prime) = {\pi \over 2} \sum_{p \in \fp_f\cap \fp_{f^\prime} } \langle \nu|_p, \mathfrak{Re} (\fa_p(f)\fa_p(f^\prime)dz)\rangle.
\tag{2.21}
\end{equation}

This implies in particular that the eigenvalues of any $y \in\mathbf{B}$ version of $y \cdot \mathbf{k}$ are determined by the asymptotics of the eigensections in $\mathbb{V}$ at the points in $\fp$.

\subsection{ A set of useful vector fields}
	The applications of the preceding formulas will use divergence zero vector fields of the sort that are described in this subsection.  Each vector field of interest is specified by a data set consisting of a point in $S^2$ to be denoted by $q$, a number from $(-1, 1)$ to be denoted by $r$, a positive number which is denoted by $a$, and a real number to be denoted by $s$.  To describe the corresponding vector field, consider first the case where $q$ is the north pole.  Introduce spherical coordinates $(\theta, \phi)$ on $S^2$ by viewing $S^2$ as the set of length 1 vectors in $\mathbb{R}^3$ and then writing the Cartesian coordinates $(x_1, x_2, x_3)$ of any given point on $S^2$ as $(x_1= \sin \theta\cos\varphi, x_2 = \sin\theta \sin \varphi, x_3 = \cos\theta$).  Granted this definition, set 
\begin{equation}
v_{(q, r, \mathfrak{a}, s)} \equiv s \chi_{r,a}(\theta){\partial \over \partial\varphi} \text{ where } \chi_{r,a} \text{ is } \chi\left({1\over a} (\cos\theta-r)\right)\chi\left( {1\over a}(r-\cos\theta)\right)
\tag{2.22}
\end{equation}

Note in particular that $v_{(q,r,\mathfrak{a},s)}$ is equal to $s {\partial \over \partial \varphi}$  where $\cos\theta = r$ and that it has compact support where $|\cos\theta - r |< a$.  The definition of $v_{(q,r,\mathfrak{a},s)}$ for any other choice of $q$ is the push-forward of the north polar version of $v_{(q,r,\mathfrak{a},s)}$ depicted above by the rotation of $S^2$ that moves the north pole to $q$ along the short great circle between them.  (In the case of the south pole, the formula is given by the analog of (2.22) with $\theta$ replaced by $\pi - \theta$ and with the $s$ changed to $-s$.)  Let $\mathcal{V}$ denote the set of all such $v_{(q,r, a, s)}$.  The next paragraph and the subsequent ones point out some of the salient properties of the vector fields in the set $\mathcal{V}$.  

	The first point of note is this:  Given a point $p \in\fp$, and a non-zero vector $\nu\in TS^2|_p$, there are vectors fields in $\mathcal{V}$ that are zero on a neighborhood in $S^2$ of $\fp-\{p\}$ and a positive multiple of $v$ at $p$.  To find such a vector field, fix a Gaussian coordinate chart centered at $p$.  Use this chart to identify a small radius disk around $p$ with the corresponding radius disk about the origin in $\real^2$. 
	 Take the radius of this disk to be much less than the distance to any point in $\fp-p$.  The Gaussian coordinate chart is uniquely specified by requiring that it identify $v$ with a vector at the origin pointing in the positive direction along the $x$ axis.  Now fix a point (which will be $q$) in the disk in $\real^2$ on the positive $y$ axis and take $r$ to be the distance from the origin to this point.  Take the number $a$ to be any positive number that is much less than $r$.  The corresponding $v_{(q,r,\mathfrak{a},+)}$ is a non-zero, positive multiple of $\nu$ at the point $p$ and it is zero at all other points in $\fp$. 

The second point of note is an immediate consequence of the first:  If the points of $\fp$ are labeled as $(p_1, \ldots, p_{2n})$ and if $(v_1, \ldots, v_{2n})$ are chosen vectors in the tangent space to $TS^2$ at the correspondingly labeled point in $\fp$, then there are choices for $2n$ vector fields from $\mathcal{V}$ (these are labled as $\nu_1, \ldots, \nu_{2n})$ that have disjoint supports, with each $\nu_k$ being zero on a neighborhood of each point in $\fp-p_k$ and being a positive multiple of $v_k$ at $p_k$.  (This follows from the construction in the preceding paragraph).  This disjoint support property implies in particular that the corresponding 1-parameter family of diffeomorphisms are pairwise commuting.
	
	A third point is that this construction can be smoothly parametrized by certain auxiliary spaces.  The first relevant instance is this:  Fix a small radius disk about a point $p \in \fp$ such that each point in the disk is much closer to $p$ than to any other point in $\fp$.  Let $D_1$ denote this disk, and let $D_0$ denote the concentric disk whose radius is  ${1\over 100}$ times that of $D_1$.   The construction described above can be used to produce a smooth map from $D_0$ to the space of area preserving diffeomorphisms of $S^2$ such that if $q \in D_0$, then the corresponding diffeomorphism maps $p$ to $q$ and fixes $S^2-D_1$.

	The second relevant instance is a generalization of the first along the following lines:  For each $k \in\{1, \ldots, 2n\}$, let $D_{k_1} \subset S^2$ denote a very small radius disk centered on $p_k\in\fp$ whose points are much closer to $p_k$ than to any other point in $\fp$.  Let $D_{k_0} \subset D_{k_1}$ denote the concentric disk whose radius is ${1\over 100}$ times that of $D_{k_1}$.  Then there is a smooth map from $\times_{k=1,\ldots,2n} D_{k_0}$ to the space of area preserving diffeomorphisms of $S^2$ with the property that if $\fq\equiv (q_1, \ldots, q_k) \in \times_{k=1,\ldots, 2n} D_{k_0}$ is any given point, then the corresponding diffeomorphism (denoted later by $\phi_{q,p}$) maps $p_1$ to $q_1$ and $p_2$ to $q_2$, and so on.  The diffeomorphism will also fix any point not in $\times_{k=1, \ldots 2n}D_{k_1}$.

\subsection{The $\mathcal{I}_{(\cdot)}$ eigenvaues as functions on $\ctwon$}

	The proposition that follows directly is the gateway to the rest of this paper as it describes the variation of the eigenvalues of the Laplacian on the $\fp \in \ctwon$ versions of $\hp$.  To set the stage for the proposition:  Keep in mind that $\ctwon$ is a smooth manifold.  In particular, its tangent space at any given configuration $\fp$ can be identified with the vector space $\oplus_{p\in\fp} TS^2|_p$.  By way of notation:  The proposition refers to notation introduced in Lemma 2.3 and in the subsequent paragraph for the asymptotics near points in a given $\ctwon$ configuration of an associated $\mathcal{I}_{(\cdot)}$-eigensection.\\

\noindent \bf Proposition 2.6:  \it  Having fixed $\fp\in\ctwon$, let $\lambda_p$ denote a corresponding $\ip$-eigenvalue.
\begin{itemize}[leftmargin=*]
\item	If $\lambda_{\fp}$ has multiplicity 1:  There is an open neighborhood of $\fp$ in $\ctwon$ and a smooth function $\lambda_{(\cdot)}$ on this neighborhood with the following properties:
\begin{enumerate}
\item	The value of $\lambda_{(\cdot)}$ at $\fp$ is $\lambda_{\fp}$.
\item	If $\fq$ is in this neighborhood, then $\lambda_{(\fq)}$ is an $\mathcal{I}_{\fq}$-eigenvalue with multiplicity one.
\item The pairing of the differential of the function $\lambda_{(\cdot)}$ at $\fp$ with any given vector $\nu\in\oplus_{p\in\fp} T S^2|_p$ is the number
$${\pi \over 2} \sum_{p\in \fp_f} \langle \nu|_p, \mathfrak{Re}(\fa_p(f)^2dz)\rangle$$
with $f$ denoting here an $\ip$-eigensection whose eigenvalue is $\lambda_{\fp}$ and whose square has $S^2$ integral equal to one. 
\end{enumerate}
\item	If $\lambda_{\fp}$ has multiplicity $N > 1$:  There exist $\epsilon>0$ and a set of $N$ continuous functions on an open neighborhood of $\fp$, these having the following properties:
\begin{enumerate}
\item	All $N$ of the functions are equal to $\lambda_{\fp}$ at the point $\fp$.
\item	If $\fq$ is in the aforementioned open neighborhood, then the values of these $N$ functions at $\fq$ are the $\mathcal{I}_{\fq}$-eigenvalues that are in the interval $(\lambda_{\fp}-\epsilon, \lambda_{\fp}+\epsilon)$. 
\item Let $\mathbb{V}$ denote the span of the $\ip$-eigenvectors with eigenvalue $\lambda_{\fp}$.  Having fixed a vector $\nu \in \oplus_{p\in\fp} T S^2|_p$, let $(\eta_1, \ldots, \eta_N)$ denote the eigenvalues of the bilinear form on $\mathbb{V}$ whose value on any two $\ip$-eigensections $f$ and $f^\prime$ from $\mathbb{V}$ is
$${\pi \over 2} \sum_{p\in \fp_f \cap \fp_{f^\prime}} \langle \nu|_p, \mathfrak{Re}(\fa_{p}(f)\fa_{p}(f^\prime)dz)\rangle$$
  	Let $\gamma: (-\epsilon, \epsilon) \to\ctwon$ parametrize a smooth path in the open set with $\gamma(0) = \fp$ and with tangent vector $\nu$ at $\fp$.   The values of the $N$ functions at the configuration $\gamma(t)$ differ from the $N$ numbers $\{\lambda_{\fp}+ t\eta_1, \ldots, \lambda_{\fp} + t\eta_N\}$ by $o(t)$. \\
\end{enumerate}
\end{itemize}

\noindent \bf Proof of Proposition 2.6:  \rm This is an instance of Proposition 2.5 given what is said in Sections 2.7 and 2.8.

\subsection{Critical points of $\mathcal{I}_{(\cdot)}$ eigenvalues as functions on $\ctwon$}

	Suppose here that $\fp \in \ctwon$ and that $\lambda$ is an $\ip$-eigenvalue with multiplicity 1.  The top bullet of Proposition 2.6 says $\lambda$ is the restriction to $\fp$ of a smooth function on a neighborhood of $\fp$ whose restriction to any given point (call it $\fq$) in this neighborhood is an $\mathcal{I}_{\fq}$-eigenvalue with multiplicity 1.  Denote this function by $\lambda_{(\cdot)}$.  The lemma that follows makes a formal assertion to the effect that this function $\lambda_{(\cdot)}$ has no critical points.\\

\noindent \bf Lemma 2.7: \it A smooth function on an open set in $\ctwon$ whose value at any given configuration is an $\mathcal{I}_{(\cdot)}$-eigenvalue with multiplicity 1 can not have critical points on that open set.\\

\noindent \bf Proof of Lemma 2.7:  \rm
The tangent space to $\ctwon$ at a given point $\fp \in \ctwon$ is $\oplus_{p\in \fp} TS^2|_p$.  With this understood, suppose that the given eigenvalue function $\fp \to \lambda_{(\fp)}$ has multiplicity one at each point where it is defined and smoothly varying.  
Item (3) of the first bullet of Proposition 2.6 depicts the gradient of $\lambda_{(\cdot)}$ at a given configuration $\fp$.  The important point here is that the gradient is zero if and only if the set $\fp_f$ is empty, which is to say that the corresponding $\ip$-eigensection $f$ near any given $p \in \fp$ is
 $\mathcal{O}(\mathrm{dist}(\cdot, p)^{3\over 2}).$ 

	Keeping the preceding in mind, introduce now the three vector fields $\{\partial_1, \partial_2, \partial_3\}$ that generate the action of $SO(3)$ on $S^2$.  These can be written using the Cartesian coordinates $(x_1, x_2, x_3)$ on $\mathbb{R}^3$ as 

\begin{equation}
\partial_1 = x_2 {\partial \over \partial x_3}  - x_3{\partial \over \partial x_2}, \quad  \partial_2 = x_3 {\partial \over \partial x_1}  - x_1 {\partial \over \partial x_3},   \quad \partial_3 = x_1{\partial \over \partial x_2}  - x_2{\partial\over \partial x_1}.
\tag{2.23}
\end{equation}

(Although these are depicted as vector fields on $\mathbb{R}^3$, each is tangent to the spheres of constant radius centered at the origin.)  Because these vector fields commute with the Laplacian on $S^2$, the sections $\{\partial_af\}_{a=1,2,3}$ of $\ip$ obey the equation $d^\dag d(\cdot) = \lambda_p(\cdot)$ if $f$ does.  
Moreover, each of these sections of $\ip$ are in $\hp$ precisely when $\fp_f = \emptyset$ because this last condition says that $f$ is $\mathcal{O}\left(\mathrm{dist}(\fp,\cdot)^{3\over 2}\right)$ near $\fp$.  
(If $f$ were only $\mathcal{O}\left(\mathrm{dist}(p,\cdot)^{1\over 2}\right)$ near some $p\in \fp$, then this conclusion would be false for at least one section from the set $\{\partial_a f\}_{a=1,2,3}$.)  It follows as a consequence that $\lambda_{\fp}$ has multiplicity greater than one if $\fp_f$ is empty. In fact, with a little more work, one sees that $\lambda_{\fp}$ has multiplicity at least four.  

Because $\lambda_{\fp}$ is assumed to have multiplicity one, it follows from the preceding observation that $\fp_f$ can't be empty and thus $\fp$ can't be a critical point of $\lambda_{(\cdot)}$.

\section{The lowest eigenvalue as a function on $\ctwon$ }

The lemma that follows directly is the gateway for this section.  It concerns the function on $\ctwon$ that assigns to any given configuration the smallest $\mathcal{I}_{(\cdot)}$-eigenvalue.  (The value of this function at a given configuration $\fp$ is the number $\ep$ from (1.2).) \\

\noindent \bf Lemma 3.1:  \it The function $\mathrm{E}_{(\cdot)}$ on any given positive integer $n$ version of $\ctwon$ is a continuous function which is smooth on a neighborhood of any configuration where $\mathrm{E}_{(\cdot)}$ has multiplicity 1 as an $\mathcal{I}_{(\cdot)}$-eigenvalue.  Conversely, it is not differentiable where it has multiplicity greater than 1.\\

\noindent \bf  Proof of Lemma 3.1:  \rm The claim that $\mathrm{E}_{(\cdot)}$ is smooth on a neighborhood of any configuration where $\mathrm{E}_{(\cdot)}$ has multiplicity 1 is an instance of Lemma 2.7.  The claim that $\mathrm{E}_{(\cdot)}$ is not differentiable where it has multiplicity greater than 1 is a consequence of Proposition 2.6.  To see why is this, fix a configuration $\fp\in \ctwon$ where $\ep$ has multiplicity at least 2.  Let $V$ denote the linear span of the $\ip$-eigenfunction with the eigenvalue $\ep$.  The first observation is that there is at least one $\ip$-eigensection from $V$ with at least one of the numbers from the set $\{n_p(\cdot)\}_{p\in\fp}$ equal to zero (see Lemma 2.3 for the definition).  The reason is this:  If $f$ is an $\ip$-eigensection with all of the $\{n_p(f)\}_{p\in\fp}$ being positive, then the action on $f$ of rotation vector fields depicted in (2.23) done some finite number of times will produce an $\ip$-eigensection with the same eigenvalue as $f$ and with at least one of the $n_{p}(\cdot)$’s equal to zero.  (Looking ahead, Lemma 3.9 asserts that at least $n$ of the values of $n_p(\cdot)$ are non-zero for any non-trivial element from $V$.)  With the preceding understood, it follows then that there exists at least one point $p \in \fp$ and a non-zero $\nu \in TS^2|_p$ for which the corresponding quadratic form from Item (3) of the second bullet in Proposition 2.6 is not identically zero.  This quadratic form has either rank 1 or rank 2.  In any case, it follows from that same Item (3) that the function $\mathrm{E}_{(\cdot)}$ near $\fp$ is at best a Lipshitz function of local coordinates for $\ctwon$ defined near $\fp$.

The remaining parts of this section study the behavior of this function $\mathrm{E}_{(\cdot)}$.

\subsection{A compactification of $\ctwon$}

A smooth function on a compact manifold has critical points which reflect the topology of the manifold via Morse theory.  By way of comparison, a function on a non-compact manifold need not have critical points.  For example:  Take the compliment in any compact manifold of the critical points of any given function.  With the preceding pathology in mind, it proves useful to introduce a compactification of $\ctwon$ that is compatible with the function $\mathrm{E}_{(\cdot)}$.  The description of this compactification and its topology has two parts.\\
	  
	\it Part 1: \rm  Introduce some terminology:  A $\integer/2\integer$ divisor on $S^2$ is an equivalence class of finite, formal sums of the form 
\begin{equation}
\sum_{p\in\fp}m_p p
\tag{3.1}
\end{equation}
with $\fp$ denoting an unordered, finite set of distinct points and with $m_p$ denoting either 0 or 1.  (The set $\fp$ is called a \it configuration set.\rm)  The \it degree \rm of the divisor is the number of terms in the sum with $m_p =1$.   Two such formal sums are deemed to be equivalent if they differ by terms of the form $m_p p$ with $m_p = 0$.  This understood, each equivalence class is represented in a unique way by a formal sum with $m_p = 1$ for all $p$.  This representative is said to be the \it minimal representative  \rm because its configuration set has the minimal possible number of elements.  By way of an example:  An element in any positive integer $n$ version of $\ctwon$ defines a minimal representative of a $\integer/2\integer$ divisor with degree $2n$.   Conversely, the configuration set of the minimal representative of any $\integer/2\integer$ divisor defines an element in $C_{2n^\prime}$ with $n^\prime$ being half the degree of the divisor.\\  

\it 
	Part 2: \rm Let $\ctwonb$ denote the union of the set of $\integer/2\integer$ divisors with even degree at most $2n$ (including degree 0). Thus, as a point set, this is the disjoint union
\begin{equation}
\ctwonb = \ctwon \cup C_{2n-2} \cup \cdots \cup C_0  
\tag{3.2}
\end{equation}
with $C_0$ denoting here the unique minimal representative of the equivalence class of $\integer/2\integer$ divisors with degree 0.  This part of the subsection defines a topology on $\ctwonb$ that renders it compact.

A basis for the topology of $\ctwonb$ consist of sets labeled by a positive number to be denoted by $\epsilon$, an integer from the set $\{0, \ldots, n\}$ to be denoted by $k$, and a degree $2k$ divisor to be denoted by $[\fq]$.  The corresponding set is denoted by $\mathcal{U}(\epsilon, k, [\fq])$; the criteria for membership is given next.  A $\integer/2\integer$ divisor $[\fp]$ is in the set 
$\mathcal{U}(\epsilon, k, [\fq])$ when $[\fp]$ and $[\fq]$ have respective representative sums $\sum_{p\in\fp} m_p p$ and $\sum_{q\in\fq} m_q q$ with $\fp$ having a partition as a disjoint union of non-empty subsets that are indexed by the points in the configuration set $\fq$ representing the divisor $[\fq]$ such that the following is true:  The partition subset $\fp(q) \subset \fp$ labeled by any given point $q \in\fq$ obeys \\

\begin{itemize}
\item dist$(p, q) < \epsilon$  for all $p \in \fp(q)$.
\item $\sum_{p\in\fp(q)} m_p = m_q$ mod(2).
\end{itemize}
\hfill (3.3)

(The definition of membership in $\mathcal{U}(\epsilon, k, [\fq])$ does not depend on the chosen representative $\sum_{p\in\fp} m_p p$ of the equivalence class $[\fp]$ because adding a term $m_{p^\prime}p^\prime$ to the formal sum $\sum_{p\in\fp} m_p p$ with $m_{p^\prime} = 0$ can be accommodated by adding the same $m_{p^\prime}p^\prime$ to the formal sum that represents the equivalence class of the $\integer/2\integer$ divisor $[\fq]$.)  

	A point to note about these $\mathcal{U}(\cdots)$ basis sets is that they are all contractible.  Indeed, any given $\mathcal{U}(\epsilon, k,[\fq])$ deformation retracts on to corresponding $[\fq]$; the retraction moves the set of points from each partition subset $\fp(q)$ that appears in (3.3) onto $q$ by moving each point from $\fp(q)$ along the short great circle arc between the point and $q$. \\

	The proposition that follows summarizes three relevant properties of this space $\ctwonb$. \\

\noindent \bf Proposition 3.2:  \it  Fix an integer $n$ so as to define the spaces $\ctwon$ and $\ctwonb$.  Then

\begin{itemize}
\item The space $\ctwonb$ is compact

\item The subset $\ctwon$ is open and dense in $\ctwonb$.

\item 
The lowest eigenvalue as a function on $\ctwon \cup \ \mathcal{C}_{2n-2} \cup \cdots \cup \mathcal{C}_0$ is continuous with respect to the topology on $\ctwonb$ .

\end{itemize}

\rm

The rest of this subsection is occupied with the lemma's proof.\\

\noindent \bf Proof of Proposition 3.2:  \rm Some background is needed in order to prove the first two bullets:  The permutation group on $k$ letters ($k$ is a positive integer) acts on $\times_k S^2$ by permuting the factors.  The quotient is denoted by $\sym_k(S^2)$ which is a compact space.  It is homeomorphic to the space of non-trivial, homogeneous, degree $2n$ polynomials in the coordinates of $\mathbb{C}^2$ modulo the action of $\mathbb{C}^*$, thus, $\mathbb{CP}^{2n}$.  The homeomorphism sends a polynomial to its $2n$ zeros on $\mathbb{CP}^1$.  Said differently, any polynomial of the sort under consideration can be factored as
\begin{equation}
\prod_{1\leq j \leq 2n}(\alpha_j z_1 + \beta_j z_2)\tag{3.4}
\end{equation}
where the notation has $(z_1, z_2)$ denoting the coordinates for $\mathbb{C}^2$; and where each pair $(\alpha_j, \beta_j)$ are not both zero and are defined modulo $\mathbb{C}^*$.  The $k = 2n$ version is relevant because the manifold $\ctwon$ 
appears as the dense open set in $\mathbb{CP}^{2n}$ that consists of the equivalence classes of polynomials with distinct factors (the $\mathbb{C}^*$ equivalence classes of the pairs $\{(\alpha_j, \beta_j)\}$ are the constituent points of a configuration in $\ctwon$).
  
With regards to the first bullet:  Let $\{\fp_j\}_{j\in\mathbb{N}}$ denote a sequence in $\ctwon$ (degree $2n$ polynomials with distinct factors modulo $\mathbb{C}^*$) with no convergent subsequence in $\ctwon$.  It will have a subsequence that converges in $\mathbb{CP}^{2n}$, thus to a homogeneous, degree $2n$ polynomial.  That polynomial will have some factors appearing multiple times.  Discard all factors that appear an even number of times and replace all factors that appear an odd number of times with that factor appearing just once.  This new polynomial corresponds to a point in a $
\mathcal{C}_{2k}$ term in (3.1) for some $k < n$.  (But for notation, the preceding proves that any sequence in any $
\mathcal{C}_{2k}$ term in (3.1) has a limit point in some $
\mathcal{C}_{2j}$ term for $j \leq k$.) 	

With regards to the second bullet and $\ctwon$ being dense:  Any point in any $k < 2n$ version of $C_{2k}$ (a homogeneous, degree $2k$ polynomial with distinct roots) can be obtained as a limit of polynomials with $2n$ distinct roots in the manner just described by appending $2(n-k)$ linear factors with nearby but distinct roots in $\mathbb{CP}^1$ and then letting those added roots come together along the sequence.  In this way, the original polynomial when viewed as a point in any $C_{2k}$ term of (3.2), appears as a limit point of some non-convergent sequence in $\ctwon$.   

With regards to the second bullet and $\ctwon$ being open:  Let $\fp$ denote a point in $\ctwon$.  Take $\epsilon>0$ but less by a factor of  ${1\over 100}$ than the minimum of the distances between $\fp$'s constituent points.  Take $k = n$ and $[\fq] = [\fp]$.  The corresponding set $\mathcal{U}(\epsilon, n, [\fp])$ contains only points from $\ctwon$.

	The proof of the proposition's third bullet has four parts.  To set the stage for this, take $\{\fp_i\}_{i\in \mathbb{N}}$ to be a sequence in $\ctwon$ that converges in $\ctwonb$ to a point in some $C_{2k}$ term in (3.2).  This configuration is denoted by $\fq$.  The third part of the proof of the third bullet explains why $\lim_{i\to\infty} E_{\fp_i} \geq E_{\fq}$; and the fourth part explains why $\lim_{i\to\infty}E_{\fp_i}  \leq E_{\fq}$.  The first two parts of the proof set the stage for the latter two parts.\\
	
		Part 1:  The definition of convergence implies first that there exists a set of distinct points $\mathfrak{k} \subset S^2-\fq$ with at most $2(n-k)$ members with the following significance: Given $\epsilon>0$, then all sufficiently large $i$ versions of $\fp_i$ have the following two properties: To state the first, suppose for the moment that $q\in\fq$ and let $D_\epsilon(q)$ denote the radius $\epsilon$ disk centered at $q$.  This disk contains some odd number of constituents of $\fp_i$.  No generality is lost by assuming that this number is independent of $\epsilon$ and the index $i$ when $i$ is sufficiently large with `large' determined by the choice of $\epsilon$.  To state the second, let $y$ denote a point from $\fk$.  Then the corresponding radius $\epsilon$ disk centered on $y$ contains some non-zero, even number of points from $\fp_i$.  The number of points in each such disk can be taken to be independent $\epsilon$ and $i$ when $i$ is large with `large' again determined by $\epsilon$.

A sequence of choices for $\epsilon$ will be made in what follows with the sequence decreasing and converging to zero.   An upper bound for the elements in this sequence is chosen so that each version of $\sqrt\epsilon$ is much less than the distance between any two distinct points from the set $\fq \cup \mathfrak{k}$. \\

\it Part 2: \rm  Now comes a key observation:  With $\epsilon$ small (with upper bound as just described), then each sufficiently large $i$ version of the real line bundle $\mathcal{I}_{\fp_i}$ is isomorphic to the line bundle $\mathcal{I}_{\fq}$ on the complement of the union of the radius $\epsilon$ disks centered at the points in $\fq \cup \mathfrak{k}$.  This is because all of the points in any large $i$ version of $\fp_i$ are contained in these disks as are the points of $\fq$; and because a line bundle's isomorphism class on $S^2- \cup_{z\in \fq \cup \mathfrak{k}} D_\epsilon(z)$ is determined by the isomorphism classes on the boundaries of the disks from the set 
$\{D_\epsilon(z): z\in\fq\cup\mathfrak{k}\}$; and because the respective isomorphism classes in the case of $\fp_i$ and $\fq$ are the same when $i$ is large (the M\"obius bundle if $z \in \fq$ and the product bundle if $z 
\in\mathfrak{k}$).  By the same token, the line bundles  $\mathcal{I}_{\fp_i}$ and $\mathcal{I}_{\fq}$ are isomorphic on the complement of the union of the radius $2\epsilon$ disks centered at the points in $\fp_i$.  \\

Part 3:  For each index $i$, let $f_i$ denote a chosen $\mathcal{I}_{(\cdot)}$-eigensection for the $\fp_i$ version of the Laplacian with eigenvalue $\lambda$.  (In the instance of the proposition, $\lambda =E_{\fp_i}$).  It is assumed in what follows that the $S^2$ integral of each $|f_i|^2$ is equal to 1.  Define the function $\chi_\epsilon$ as in (2.2) with $\fp$ denoting $\fp_i$.  
This version will be denoted by $\chi_{\fp_i, \epsilon}$.  Use the isomorphism from Part 2 to view $\chi_{\fp_i, \epsilon} f_i$ as an element in $\mathbb{H}_{\fq}$.  Viewed in this light, the $S^2$ integrals of the square of the norms of $d(\chi_{\fp_i,\epsilon} f_i)$ and  $\chi_{\fp_i,\epsilon} f_i$ obey

\begin{itemize}
\item $\displaystyle \int_{S^2} |d(\chi_{\fp_i, \epsilon} f_i)|^2 = \int_{S^2} |df_i|^2 + \mathfrak{e}_{1,i}$,
\item $\displaystyle\int_{S^2} |\chi_{\fp_i, \epsilon} f_i|^2 = \int_{S^2} |f_i|^2 + \mathfrak{e}_{2,i}$,
\end{itemize}
\hfill (3.5)\\
with $ \mathfrak{e}_{1,i}$ and $ \mathfrak{e}_{2,i}$ obeying 
\begin{equation}
 |\mathfrak{e}_{1,i}| \leq c_0 {1\over |\ln \epsilon|} n \lambda \text{ \ \ \ and \ \ \ }
  |\mathfrak{e}_{2,i}| \leq c_0\epsilon n \lambda.
\tag{3.6}
\end{equation}
To explain:  Both bounds are obtained by invoking the bound on $|f_i|$ from the second bullet of Lemma 2.1.  The $\mathfrak{e}_{1,i}$ bound follows from that sup-norm bound and the fact that the $S^2$ integral of $|d\chi_{\fp_i,\epsilon} |^2$ is bounded by $c_0 n {1\over |\ln\epsilon|}$.  The $\mathfrak{e}_{2,i}$ bound follows from the sup-norm bound and the fact that the area of the support of $\chi_{\fp_i,\epsilon}$  is bounded by $c_0 n \epsilon$.

	The $\lambda=E_{\fp_i}$  versions of (3.5) and (3.6) imply directly that 
	$\lim_{i\to\infty }  E_{\fp_i} \geq E_{\fq}$.\\

	Part 4:  Let $f$ denote a chosen, minimal eigenvalue eigensection for the $\fq$ version of the Laplacian.  (Its eigenvalue is $E_{\fq}$.)  Having fixed positive $\epsilon$ to be very small, then $\chi_\epsilon f$ can be viewed as a section of each sufficiently large $i$ version of the line bundle $\mathcal{I}_{\fp_i}$ by choosing an isomorphism between $\mathcal{I}_{\fq}$ and the large $i$ versions of  $\mathcal{I}_{\fp_i}$  on the complement of the support of $\chi_\epsilon$.   But for terminology, the argument that lead to (3.5) and (3.6) can be repeated to see that these same inequalities hold with $f_i$ replaced by $f$ and with respective error terms $\mathfrak{e}_1$ 
	and $\mathfrak{e}_2$ that obey the bounds in (3.6) with $\lambda$ replaced by $E_{\fq}$.  Granted that these inequalities hold for each sufficiently small $\epsilon$, it then follows directly that 
	$\lim_{i\to\infty} E_{\fp_i} \leq E_{\fq}$ also.

\subsection{ How large can $\mathrm{E}_{(\cdot)}$ get?}

Proposition 3.2 implies that the supremum of $\mathrm{E}_{(\cdot)}$ on $\ctwon$ is a non-decreasing function of the integer $n$.  The following two propositions say more about this supremum.\\

\noindent \bf Proposition 3.3:  \it There are configurations in $\ctwon$ with smallest Laplace eigenvalue greater than $\kappa^{-1} n$  with $\kappa > 1$ and independent of $n$.\\ \rm

This proposition is proved momentarily.   A consequence of Propositions 3.2 and 3.3 follows directly.\\

\noindent \bf Corollary 3.4:  \it There exists an unbounded set of positive integers with the following property:  Supposing that $n$ is from this set, there is a configuration in $\ctwon$ where $\mathrm{E}_{(\cdot)}$ is equal to its supremum of $\mathrm{E}_{(\cdot)}$ over the whole of $\ctwon$. \\

\normalfont The rest of this subsection is dedicated to proving Proposition 3.3.\\

\noindent\bf Proof of Proposition 3.3:  \rm The proof has three parts.\\

\it Part 1:  \rm Fix a small positive number to be denoted by $R$ and then fix a maximal (in number) configuration of points in $S^2$ such that the distance between any two distinct points is no less than $R$.  Let $\fp_0$ denote this configuration.  Each point in $S^2$ has distance at most $R$ from some point from the set $\fp_0$ because $\fp_0$ is maximal.  The number of elements in $\fp_0$ (denoted by $|\fp_0|$) is such that
\begin{equation}
c_0^{-1} {1\over R^2} \leq |\fp_0| \leq c_0 {1\over R^2}.
\tag{3.7}
\end{equation}
Note that the disks of radius $R$ centered at the points in $\fp_0$ cover $S^2$.  Also, the disks of radius ${1\over 2}R$ centered at these points are pairwise disjoint.

Assign to each point from $\fp_0$ a pair of points that lie at distance ${1\over 8}R$ from the given point on opposite sides of a chosen great circle through that point.  When $q$ denotes the chosen point from $\fp_0$, then these other two points are denoted by $q_+$ and $q_-$.   Let $\fp$ denote the configuration consisting of all such $(q_-, q_+)$ pairs coming from the points in $\fp_0$.  Since the number of points in the configuration $\fp$ is twice the number from $\fp_0$, this number also obeys the inequality in (3.7) with a possibly different $c_0$.  An important point for what follows:  The disks of radius  ${1\over 8}R$ centered at the points from the set $\fp$ are pairwise disjoint. \\
	
\it Part 2: \rm What follows is the fundamental lemma about sections of $\mathcal{I}_{\fp}$:\\

\noindent \bf Lemma 3.5:  \it  There exists $\kappa > 1$ such that the following is true when $R < {1\over \kappa}$ .  Fix $q \in \fp_0$ and let $D$ denote the radius $R$ disk centered at $q$.  If $f$ is any section of $\ip$ over $D$ with $|df|^2$ having finite integral over $D$, then 
$$\int_D |df|^2 \geq {1\over \kappa R^2} \int_D f^2.$$

\noindent\bf Proof of Lemma 3.5:  \rm Let $D^\prime$ denote the radius  ${1\over 8}R$ disk centered at either $q_+$ or $q_-$ and let $f$ denote a section of $\mathcal{I}$ over $D^\prime$.  The integral of $|df|^2$ over $D^\prime$ is no less than $c_0^{-1} {1\over R^2}$  times that of $f^2$ because $f$ must vanish at one or more points on any concentric circle about the center point of $D^\prime$.  This is because $\ip$ restricts to any of these circles as the M\"obius line bundle. 

With the preceding understood, there are now two cases to consider:  The first is when the integral of $f^2$ over the union of the radius ${1\over 8}R$ disks centered at $q_+$ and $q_-$ is greater than ${1\over 1000}$  times its integral over $D$.  In this case, the assertion in the lemma follows from what was said in the previous paragraph about the integral of $|df|^2$ over the $q_+$ and $q_-$ versions of $D^\prime$.  In the second case, the integral of $f^2$ over $D^\prime$ is no greater than ${1\over 1000}$  times its integral over the union of the $q_+$ and $q_-$ versions of $D^\prime$.  In this case, the assertion of the lemma follows from the fact that $|f|$ isn't uniformly close to the constant function on $D$ in a suitable topology.  Indeed, by rescaling $D$ to have radius 1, the bound amounts to a straightforward observation about the Dirichelet energy of honest functions on the radius 1 disk in $\mathbb{R}^2$.\\

	\it Part 3:  \rm Proposition 3.3 follows directly from the next lemma and (3.7).  \\

\noindent \bf Lemma 3.6: \it There exists $\kappa > 1$ such that the following is true when $R < {1\over \kappa}$.  Construct $\fp$ as described above and let $f$ denote any given section of $\ip$ from the Hilbert space $\mathbb{H}_{\fp}$.  Then
$$\int_{S^2} |df|^2 \geq {1\over \kappa R^2} \int_{S^2} f^2.$$

\noindent\bf Proof of Lemma 3.6:  \rm The integral of $|df|^2$ over $S^2$ is no less than $c_0^{-1}$ times the sum of its integrals over the various $q \in \fp$ versions of $D$.  This is because at most $c_0$ of these disks intersect at any one point.  Thus, the sum of the integrals of $|df|^2$ over these disks overcounts the integral of $|df|^2$ by a factor of at most $c_0^{-1}$.
\begin{equation}
\int_{S^2} |df|^2 \geq c_0^{-1} \sum_{q\in\fp_0} \int_{D(q)}|df|^2,
\tag{3.8}
\end{equation}
with $D(q)$ denoting here the radius $R$ disk centered at $q$.  With (3.8) in hand, appeal to Lemma 3.5 to see that 
\begin{equation}
\int_{S^2} |df|^2 \geq c_0^{-1} {1\over R^2} \sum_{q\in\fp_0} \int_{D(q)} f^2,
\tag{3.9}
\end{equation}
Meanwhile, the sum of the various $D(q)$ integrals of $f^2$ for $q \in \fp_0$ is no less than the integral of $f^2$ over $S^2$ because the union of these $D(\cdot)$'s covers $S^2$.  This last observation with (3.9) implies directly what is asserted by the lemma. 

\subsection{No critical values of $\mathrm{E}_{(\cdot)}$}

	Corollary 3.4 asserts in effect that there are infinitely many values of $n$ where $\mathrm{E}_{(\cdot)}$ on $\ctwon$ has a global maximum.  However, according to Lemma 2.7, these are not configurations where the lowest eigenvalue has multiplicity 1; and so there is no guarantee that $\mathrm{E}_{(\cdot)}$ has a differential at these points (let alone, one that vanishes).  Even so, with the second bullet of Proposition 2.6 in mind, one can ask whether $\mathrm{E}_{(\cdot)}$ has a critical value in the sense of the next definition.  This definition refers to the collection of numbers $\{n_p(\cdot)\}_{p\in\fp}$ from Lemma 2.3.  To set the notation:  Supposing that $k$ is a positive integer, the definition refers to the $k$'th \it eigenvalue function, \rm $\lambda_k$, whose value at any given configuration $\fp$ is the $k$'th lowest $\ip$-eigenvalue of the Laplacian.  Thus $\lambda_1(\fp) = E_{\fp}$. \\

\noindent\bf Definition 3.7:  \it A configuration $\fp \in \ctwon$ is \underline{k-critical} if there exists an $\ip$-eigenfunction $f$ with eigenvalue $\lambda_k(\fp)$ with all numbers from the set $\{n_p(f)\}_{p\in \fp}$ being positive.  \\

\rm 
	The following lemma makes a formal statement to the effect that there are no 1-critical points in $\ctwon$.\\

\noindent \bf Lemma 3.8:  \it If $\fp \in\ctwon$, and if $f$ is an $\ip$-eigensection with eigenvalue $E_{\fp}$, then at least $n+1$ integers from the set $\{n_p(f)\}_{p\in\fp}$ are equal to zero.\\

\normalfont With a look towards finding homogeneous $\integer/2\integer$ harmonic 1-forms on $\mathbb{R}^3$:  They can't correspond (as described in Section 1.6) to $\ip$-eigensections with the lowest eigenvalue $\ep$ because $f$ is an $\ip$-eigensection and if $n_p(f) = 0$ for a given $p \in \fp$, then $|f|$ near $p$ is greater than a non-zero, constant multiple of $\text{dist}(\cdot, p)^{1/2}$.\\

\noindent \bf Proof of Lemma 3.8: \rm  Let $f$ denote a given $\ip$-eigensection with eigenvalue $E_{\fp}$, and given $f$, let $\mathfrak{c}$ denote the set of points in $S^2-\fp$ where both $f$ and $df$ are zero.  Also, let $\mathfrak{Z}$ denote the zero locus of $|f|$ (including the points in $\fp$).  As explained directly, the set $\mathfrak{Z}$ is an embedded graph in $S^2$ with  $\fp \cup \fc$ being its vertex set.  To see this, note first that the complement of 
$\fp \cup \fc$ in $\mathfrak{Z}$ consists of a disjoint union of embedded open arcs because it is a level set of $f$ where $df \neq 0$.  With regards to $\mathfrak{Z}$ near points in $\fp$:  Lemma 2.3 depicts $f$ near any given point in $\fp$, and the depiction implies that $\mathfrak{Z}$ intersects a very small, radius disk centered at $p$ as the union of $2n_p+1$ embedded,
 closed arcs with one endpoint at $p$ and the other on the boundary of the disk.  (The angle at $p$ between consecutive arcs impinging on $p$ is $2\pi/(2n_p+1)$.)  Thus, $\mathfrak{Z}$ near any given $p \in \fp$ is homeomorphic to a neighborhood of a vertex in a graph with that vertex having $2n_p+ 1$ incident half-edges.  With regards to $\mathfrak{Z}$ near points in $\fc$:  If $p\in \fc$, then $f$ near $p$ (which is an $f^{-1}(0)$ critical point) has a Taylor's expansion with respect to a complex coordinate centered at that point (it is denoted by $z$ below) which looks like
\begin{equation}
f = \mathfrak{Re}(\mathfrak{a} z^{m_p}) + \mathcal{O}(|z|^{m_p+1})
\tag{3.10}
\end{equation}
with $m_p$ being a positive integer greater than 1 and with $\fa$ being a non-zero complex number.  In this case, $\mathfrak{Z}$ near $p$ is homeomorphic to a vertex in a graph with that vertex having $2m_p$ incident half-edges.  

Given this graph structure for $\mathfrak{Z}$, it then follows that its Euler characteristic is

\begin{equation}
\chi_{\mathfrak{Z}}= 2n + |\fc| - {1\over 2} \sum_{p\in\fp} (2n_p+1) - {1\over 2} \sum_{p\in\mathfrak{c}}2m_p. 
\tag{3.11}
\end{equation}

If $\chi_3$ is non-positive, then $\mathfrak{Z}$ must have a closed cycle because a tree has Euler characteristic 1 and the Euler characteristic of a union of trees is the number of trees in the union.  As explained momentarily, $\mathfrak{Z}$ can't have any closed cycles.  Assuming this, then positivity of $\chi_3$ requires that $n$ be greater than $\sum_{p\in\fp} n_p$ which requires in turn that at least n+1 points in $\fp$ have $n_p = 0$.

To see why $\mathfrak{Z}$ can't have any closed cycles, assume to the contrary that it does to generate nonsense.  If it does, then its complement in $S^2$ is disconnected.  Assuming that, let $\Omega$ denote one of the components if $S^2-\mathfrak{Z}$.   Now let $f_\Omega$ denote the section of $\ip$ 
that is equal to $f$ on $S^2-\Omega$ and equal to $-f$ on $\Omega$.  Supposing that the $S^2$ integral of $f^2$ is 1, then this is also the case for the $S^2$ integral of $f_\Omega^2$.  And, it is also the case that the $S^2$ integral of $|df_\Omega|^2$ is equal to
 $E_{\fp}$. Even so, $f_\Omega$ can't be an eigensection for the Laplacian because $f + f_\Omega$ would then be an eigensection and $f + f_\Omega$ is zero on an open set.  But if $f_\Omega$ isn't an eigensection, then the integral of $|df_\Omega|^2$ must be greater than $\ep$ since $\ep$ is the minimal eigenvalue of the $\ip$-Laplacian and since $f_\Omega$ can, in any event, be written as a linear combination of its eigensections.  This last conclusion is required.

\subsection{The $k$'th lowest eigenvalue}

Supposing that $k$ denotes a positive integer, introduce the function $\lambda_k$ on $\ctwon$ whose value at any given configuration $\fp$ is the $k$’th lowest $\ip$-eigenvalue.  To be precise about the definition:  The value of $\lambda_k(\cdot)$ at $\fp$ is the supremum of the set of numbers $s \in[0, \infty)$ such that there are at most $k$-1 linearly independent $\ip$-eigensections with eigenvalue less than $s$.  It follows from this definition that $\lambda_k(\fp)$ is an $\ip$-eigenvalue.  It also follows from this definition that if the number of linearly independent $\ip$-eigensections with eigenvalue less than $\lambda_k(\fp)$ is less than $k-1$, then that number plus the multiplicity of $\lambda_k(\fp)$ is greater than or equal to $k$.  The function $\lambda_k(\cdot)$ is continous on $\ctwon$ and it is smooth on some neighborhood of any configuration where it has multiplicity 1.  (These assertions follow from Proposition 2.6.)  Lemma 3.8 says that there are no configurations where $\lambda_k$ is 1-critical (in the sense of Definition 3.7).  The next lemma says that $\lambda_k(\cdot)$ can’t have a $\fp$-critical configuration for small values of $k$.\\


\noindent \bf Lemma 3.9: \it  Fix $k \leq n$; and suppose that $\fp\in\ctwon$, and that $f$ is an $\ip$-eigensection with eigenvalue $\lambda_k$.  Then at least $n+1-k$ integers from the set $\{n_p(f)\}_{p\in\fp}$ are equal to zero.\\

\it

\noindent \bf Proof of Lemma 3.9:  \rm Assume that $f$ has eigenvalue $\lambda$ and that $n_p(f)\geq 1$ for a subset of $n$ or more points $p \in\fp$.  Write the size of this set as $n+\ell$.  It then follows from (3.11) that the graph $\mathfrak{Z}$ (which is the $f=0$ locus in $S^2$) has Euler characteristic $-\ell$  or less.  This implies in turn that the complement of $\mathfrak{Z}$ in $S^2$ has at least $\ell + 2$ components.  Letting $K$ denote the number of components of $S^2-\mathfrak{Z}$, label these components as $\{\Omega_\alpha\}_{\alpha = 1,2,\ldots,K}$. Now define $f_\alpha$ for $\alpha \in\{1, \ldots, K\}$ to be the section of $\ip$ that is given by the rule $f_\alpha \equiv f$ on $\Omega_\alpha$ and $f_\alpha \equiv 0$ on $S^2-\Omega_\alpha$.  These $K$-sections span a subspace of $\hp$ where the inequality
 
\begin{equation}
\int_{S^2}  |dh^2| \leq \lambda \int_{S^2} |h|^2  
\tag{3.12}
\end{equation}
holds.  Since $f$ is the only actual eigensection in this span (otherwise, there would be an eigensection vanishing on an open set), the complement of $f$ in this subspace must project injectively to the span of the $\ip$-eigensections with eigenvalue less than $\lambda$.  Thus, there must be at least  $\ell+ 1$ of the latter since $K \geq \ell+2$.  Thus $\lambda$ can't be the $k$'th eigenvalue for $k \leq  \ell+1$.  Turning this around, if $\lambda$ is the $k$'th eigenvalue, then $\ell  \leq k-1$ and thus at least $n-(k-1)$ of the integers from the set $\{n_p(f)\}_{p\in\fp}$ are zero. \\

The last lemma in this section concerns the differentiability of the $k^{\mathrm{th}}$-eigenvalue functions; it is a $k > 1$ analog of Lemma 3.1 which is about the $k = 1$ case.

\begin{lemma}
 For any positive integer $n$ and positive integer $k$:  The $k$’th eigenvalue function $\lambda_k$ is differentiable where $\lambda_k$ has multiplicity 1.  On the other hand, it is not differentiable at any configuration in $\ctwon$ where the eigenvalue $\lambda_k$ has multiplicity greater than 1 and where there are k-1 linearly independent $\ip$-eigensections with eigenvalue strictly less than $\lambda_k$.
\end{lemma}

\noindent \bf Proof of Lemma 3.10:  \rm The proof is almost the same as the proof of Lemma 3.1.

\section{The $\rp^\infty$ bundle}

	This section introduces an $\rp^\infty$ fiber bundle over $\ctwon$ whose topological properties explain (in part) why $\mathrm{E}_{(\cdot)}$ and the other $k$'th eigenvalue functions behave the way they do. 

\subsection{No universal bundle over $\ctwon$}

	To set the stage for what is to come:  Let $\mathcal{F} \longrightarrow \ctwon$ denote the fiber bundle whose fiber over any given configuration $\fp \in \ctwon$ is $S^2-\fp$.  To elaborate on the fiber bundle structure: Given a configuration $\fp$ from $\ctwon$, Section 2.6 describes a map (denoted there  and here by $\phi_{\fp}$) from the product of small radius disks about the points comprising $\fp$ to the space of area preserving diffeomorphisms of $S^2$.   By way of a reminder:  If the points that comprise $\fp$ are labeled as $(p_1, \ldots, p_{2n})$, then the disk around a given $p_k$ from $\fp$ was denoted by $D_{0k}$.  
(Note that the radii of these disks can depend on $k$.)
Set $U_{\fp} \equiv \times_{k=1,\ldots, 2n}D_{0k}$, this projecting (after symmetrization) diffeomorphically onto an open neighborhood of $\fp$ in $\ctwon$.  (Borrowing terminology from complex analysis:  Sets of this form will be called \it polydisk neighborhoods of $\fp$\rm).  The map $\phi_{\fp}$ from $U_{\fp}$ into the space of area preserving diffeomorphisms of $S^2$ has the following property:  Supposing that $\fq =(q_1, \ldots,q_{2n})$ denotes a given configuration in $U_{\fp}$, then the corresponding diffeomorphism (which was denoted in Section 2.6 by $\phi_{\fq,\fp}$) sends each $p_k$ to the corresponding $q_k$.  This map $\phi_{\fp}$ gives $\mathcal{F}$ its local product structure as a fiber bundle because it gives a identification between $U_{\fp} \times (S^2-\fp)$ and $\mathcal{F}|_{U_{\fp}}$ that is compatible with the projection map to $U_{\fp}$.
	
Each $\fp \in \ctwon$ has its corresponding line bundle $\ip$ sitting over the fiber of $\mathcal{F}$ at $\fp$ (which is $S^2-\fp$); and one might hope that these fit together to define a global line bundle over $\mathcal{F}$.  As explained in what follows, this is not the case.  

To see that there is no such universal line bundle, note first that $\mathcal{F}$ has a second fibration with base $S^2$:  The fibering map sends $(\fp, x)$ with $\fp$ from $\ctwon$ and $x$ from $S^2-\fp$ to $x$.  The fiber is diffeomorphic to the space of $2n$ unordered, distinct points in $\mathbb{R}^2$.  (This space is denoted by $C^0_{2n}$.)  Let $D_+$ and $D_-$ denote the (closed) northern and southern hemispheres of $S^2$.  Their intersection is the equatorial circle which will be denoted by $C$.  Restricting the fibration $\mathcal{F} \to S^2$ to $D_+$ and $D_-$ and $C$ leads to a Mayer-Vietoris sequence for the $\integer/2\integer$ homology of $\mathcal{ F}$ that starts with 
\begin{equation}
 0 \longrightarrow H^1(\cf)
 \longrightarrow H^1(\cf|_{D_+})\oplus  H^1(\cf|_{D_-})
 \longrightarrow H^1(\cf|_{C})\longrightarrow \cdots\\
 \tag{4.1}
\end{equation}
(Note that all cohomology henceforth is defined also using $\integer/2\integer$ coefficients.)

To exploit this exact sequence, note first that $D_+$ retracts onto any given point in $D_+$ (and similarly for $D_-$), and that these retractions are covered by corresponding deformation retractions of $\cf$ onto the fiber at that point.  (A deformation retraction of $\cf|_{D_+}$  can be constructed with parameter $t \in [0,1]$ that acts on any given $(\fp, x)$ by first pushing the points of $\fp$ radially away from $x$ to lie outside a disk concentric to $D_+$ with slightly larger radius; and having done that, the deformation retraction then moves $x$ via a linear deformation retract of $D_+$ onto any give point in $D_+$.)

There are two immediate consequences:  First, $H^1(\cf|_{D_+} )= H^1 (C^0_{2n})$ which is $\integer/2\integer$ (see \cite{Fu}), and likewise $H^1(\cf|_{D_-}) = \integer/2\integer$.  The second is that their respective homomorphisms to $H^1(\mathcal{F}|_C)$ are injective.  
(To prove injectivity:  The fiber preserving deformation retraction of $\mathcal{F}|_{D_+}$ onto $\mathcal{F}|_{\text{point}}$ can be made when the point in question is in the circle $C$.  Because of this, the isomorphism $H^1(\mathcal{F}|_{D_+}) \longrightarrow H^1(\mathcal{F}|_{\text{point}\in C})$ factors as the chain of pull-back homomorphisms $H^1(\mathcal{F}|_{D_+}) \longrightarrow H^1(\mathcal{F}|_C) \longrightarrow H^1(\mathcal{F}|_{\text{point}\in C})$ induced by the inclusions of first the point in $C$ and then $C$ in $D_+$.)  
These two consequences imply in turn that the rank of $H^1(\mathcal{F})$ is at most 1.   (This argument is more general:  Suppose that $n \geq 2$ and that $\cf$ is a fiber bundle over $S^n$ with path connected fiber whose restriction to the respective northern and southern hemispheres of $S^n$ is diffeomorphic as a fiber bundle to the product bundle.  There is a corresponding version of (4.1) in this case which implies that the rank of the first homology of $X$ is no greater than the rank of the first homology of the fiber.)

Keeping in mind that the rank of $H^1(\cf)$ is at most 1, now consider the fibration $\pi: \cf  \to \ctwon$.  The group $H^1(\ctwon)$ is also isomorphic to $\integer/2\integer$ (see \cite{Sc}, Corollary 3.1).  Denote its generator by $\alpha$.  The element $\pi^*\alpha$ is necessarily non-zero in $H^1(\mathcal{F})$ because the fiber of $\pi$ is path connected.  Thus, $\pi^*\alpha$ generates $H^1(\mathcal{F})$ and, because it is a pull-back, it restricts as 0 to each fiber.  Because of this, its restriction to any $\fp \in\ctwon$ version of $\cf|_{\fp} = S^2-\fp$ can not be the first Stieffel-Whitney class of the line bundle $\mathcal{I}_{\fp}$.

By way of a parenthetical remark:  The obstruction to the existence of this universal $\mathcal{I}$ bundle is carried by a class in $H^2(\ctwon)$ which can be constructed in the  \v{C}ech version as follows:  Fix a locally finite open cover of $\ctwon$ by sets of the form $U_{\fp}$.  Having done that, define $\hat I_{\fp}$ over $\pi^{-1} U_{\fp}$ to be the pull-back of $\ip$ by the diffeomorphism from $\pi^{-1}U_{\fp}$ to $U_{\fp} \times (S^2-\fp)$ whose inverse sends any given pair $(\fq, x)$ to $\phi_{\fq,\fp}(x)$.  
Now suppose that $\fp^\prime$ is a second configuration from $\ctwon$ and suppose that $U_{\fp}$ and $U_{\fp^\prime}$ intersect.  The bundles $\hat I_{\fp}$ and $\hat I_{\fp^\prime}$ are isomorphic over $U_{\fp} \cup U_{\fp^\prime}$ because they have the same first Stieffel-Whitney class.  Let $\eta_{\fp^\prime\fp}$  denote an isometry between them. 
 (There are precisely two such isometries.)  Now if $\fp^{\prime\prime}$ is a third configuration from $\ctwon$ with $U_{\fp^{\prime\prime}}$ sharing points with $U_{\fp} \cap U_{\fp^\prime}$, then there are corresponding isomorphisms $\eta_{_{\fp^{\prime\prime}\fp^{\prime}}}$ and $\eta_{_{\fp\fp^{\prime\prime}}}$.  
 Composing these three leads to a self-isometry of $\hat I_{\fp}$ which is denoted by $z_{\fp\fp^{\prime}\fp^{\prime\prime}} \in\{\pm 1\}$.  The collection of these indexed by 3-tuples of points whose sets define the locally finite open cover form a 2-dimensional  \v{C}ech cocycle on $\ctwon$ whose cohomology class is independent of the choices for various versions of $\eta_{_{\fp\fp^{\prime}}}$.  This cohomology class is the obstruction class in $H^2(\ctwon)$.

By way of another remark:  Fix a point $x \in\ctwon$ and let $C^x_{2n}$ denote the codimension 2 submanifold consisting of configurations that contain $x$.  The obstruction class for the universal line bundle vanishes over $\ctwon-C^x_{2n}$ because it is the Alexander dual to $C^x_{2n}$.  Note in this regard that $H^2(\ctwon)$ is isomorphic to $\integer/2\integer \oplus \integer/2\integer$ with $C^x_{2n}$ being dual to a generator (see Corollary 3.1 in \cite{Sc}).  As explained momentarily, the vanishing of this obstruction class on $\ctwon-C^x_{2n}$ is an instance of the upcoming Lemma 4.1.

To set the stage for the lemma, suppose that $X$ denotes a path connected manifold and that $\pi: \mathcal{Z} \longrightarrow  X$ denotes a smooth fiber bundle with path connected fiber.  There is an associated locally constant sheaf over $X$ whose stalk at any given point is the $\integer/2\integer$ first cohomology of its inverse image via $\pi$; thus $H^1(\pi^{-1}(x))$.  This sheaf is denoted by $\mathcal{H}^1$.  A section of this sheaf is a continuous assignment over $X$ of a class in the first cohomology of the fibers.  For example, the restriction of a class in $H^1(\mathcal{Z})$ to each fiber defines a section of $\mathcal{H}^1$.  Lemma 4.1 gives a partial converse, a condition to guarantee that a section of $\mathcal{H}^1$ comes from an $H^1$-class  (With regards to notation:  Keep in mind that all cohomology uses $\integer/2\integer$ coefficients.)  \\

\noindent \bf Lemma 4.1:  \it Let $\pi: \mathcal{Z} \longrightarrow X$ denote a fiber bundle over $X$.  A global section of the sheaf $\mathcal{H}^1$ comes from a class in $H^1(\mathcal{Z})$ when the fiber bundle $\mathcal{Z}$ admits a section. \\

\rm In the instance of the fiber bundle $\pi: \mathcal{F} \longrightarrow \ctwon$, the section of the sheaf $\mathcal{H}^1$ at any given $\fp \in\ctwon$ is the class in $H^1(\mathcal{F}|_{\fp})$ that has non-zero pairing with all sufficiently small radius linking circles around the points from $\fp$.  Because the restriction $\mathcal{F}$ to $\ctwon - C^x_{2n}$ admits the section $\fp \longrightarrow (\fp, x)$, the relevant section of the sheaf $\mathcal{H}^1$ comes from a class in $H^1(\mathcal{Z})$.  Since classes in $H^1$ correspond to real line bundles, there is a real line bundle over $\ctwon - \mathcal{C}^x_{2n}$ whose restriction to any given fiber of $\mathcal{F}$ is isomorphic to the corresponding real line bundle $I_{(\cdot)}$ on that fiber. \\

\noindent\bf Proof of Lemma 4.1:  \normalfont The argument has three parts.  By way of a heads up for what is to come:  Because the classes in $H^1$ with $\integer/2\integer$ coefficients are in 1-1 correspondence with real line bundles, the arguments can be given in terms of real line bundles.  This is what is done in the four parts of the proof.   (Also, the existence of a real line bundle is of ultimate interest in the special case where $X = \ctwon$ and $\mathcal{Z} = \ctwon$.)  \\

\it Part 1: \rm Let $Z$ denote a given fiber of $\mathcal{Z}$.  Now fix a locally finite, good open cover of $X$ (to be denoted by $\mathfrak{U}$) such that each set $U \in \mathfrak{U}$ comes with a diffeomorphism $\Psi_U: \mathcal{Z}|_U \longrightarrow U \times Z$ that intertwines the projection map $\pi$ with the product projection to $U$. (A good cover is one where all of its sets and all non-empty intersections of its sets is contractible.)  The product space $U \times Z$ has the class from $H^1(Z)$ that restricts to each fiber as the pull-back via the $\Psi_U^{-1}$ of the section from $\mathcal{H}^1$.  Fix a real line bundle on $U \times Z$ that represents this cohomology class (for example, a product cocycle).  The pull-back of this real line bundle by $\Psi_U$ defines a real line bundle on $\mathcal{Z}|_U$ whose first Stieffel-Whitney class restricts to each fiber as the given section of $\mathcal{H}^1$.  This real line bundle is denoted by $I_U$ in what follows.\\

\it Part 2:  \rm Suppose that $U$ and $U^\prime$ are sets from $\mathfrak{U}$ that share points.  The respective real line bundles $I_U$ and $I_{U^\prime}$ on $\mathcal{Z}|_{U \cap U^\prime}$ have the same Stieffel-Whitney classes; and so they are isomorphic, and hence isometric.  Fix an isometry $\eta_{_{UU^\prime}}: I_{U^\prime} \to I_U$.  Where three sets from $\mathfrak{U}$ overlap, say $U$, $U^\prime$, and $U^{\prime\prime}$, the composition of their corresponding versions of $\eta$ determine a self-isometry of $I_U$.  This is 
\begin{equation}
\eta_{_{UU^\prime U^{\prime\prime}}}
\equiv
\eta_{_{UU^\prime }}
\eta_{_{U^\prime U^{\prime\prime}}}
\eta_{_{U^{\prime\prime}U}}.
\tag{4.2}
\end{equation}
Since there are just two self-isometries (multiplication by $+1$ and by $-1$), the collection $\{\eta_{_{UU^{\prime}U^{\prime\prime}}}\}_{_{U,U^{\prime},U^{\prime\prime}\in \mathfrak{U}}}$ defines a degree 2  \v{C}ech cocycle on $X$ (it is a priori closed by virtue of its definition via (4.2)).  
The corresponding class in the second  \v{C}ech cohomology of $X$ is the obstruction class in the context of the fiber bundle $\mathcal{F}$.  (Changing the sign of any number of $\eta_{_{UU^{\prime}}}$'s changes this cocycle by a coboundary which has no effect on the cohomology class). 

If this degree 2 obstruction class is zero, then the collection $\{\eta_{_{UU^{\prime}}}\}_{_{U,U^{\prime}\in \mathfrak{U}}}$ are a collection of transition functions that glue the various $U \in \mathfrak{U}$ versions of $I_U$ together so as to define a line bundle over the whole $\mathcal{Z}$.  \\
	
\it Part 3: \normalfont The task now is to show that the class defined by the cocycle $\{\eta_{_{UU^{\prime}U^{\prime\prime}}}\}_{_{U,U^{\prime},U^{\prime\prime}\in \mathfrak{U}}}$ is the zero class if $\mathcal{Z}$ admits a section.  To this end, it is important to keep in mind that an isometry between two real line bundles is determined globally by its identification between the corresponding lines over just a single point.  (The reason is this:  If $I$ and $I^\prime$ are isomorphic line bundles, then an isometry between them is a unit length section of $I \otimes I^\prime$ which is isomorphic to the product bundle.  And, there are only two unit length sections of the product bundle, respective multiplication by $+1$ and by $-1$.)  

With the preceding understood, fix a point $x_{_U}$ in each set $U$ from the cover, and then fix a unit length point $\hat{\iota}_U$ in the fiber of $I_U$ over that point. Use the section $s: X \to\mathcal{Z}$ and the fact that $U$ is contractible to define a product structure for $I_U$ over $s(U)$ with $\hat{\iota}_U $ defining a preferred component.   Granted these product structures, choose the isomorphism $\eta_{_{UU^{\prime}}}$ to identify the $\hat{\iota}_{_{U^{\prime}}}$ component of $I_{U^{\prime}}$ with the $\hat{\iota}_{_U}$ component of $I_U$.
   
Now consider what happens in a triple intersections, $U \cap U^{\prime} \cap  U^{\prime\prime}$:  The right most $\eta_{_{U^{\prime\prime}U}}$ isomorphism in (4.2) identifies the $\hat{\iota}_U$ component of $I_U$ with the $\hat{\iota}_{_{U^{\prime\prime}}}$ component of $I_{U^{\prime\prime}}$; and the middle isomorphism in (4.2), which is $\eta_{_{U^{\prime}U^{\prime\prime}}}$, identifies that component of $I_{U^{\prime\prime}}$ with the $\hat{\iota}_{_{U^{\prime}}}$ component of $I_{U^{\prime}}$; and then the left most isomorphism in (4.2) identifies the $\hat{\iota}_{_{U^{\prime}}}$ component of $I_{U^{\prime}}$ with the $\hat{\iota}_{_U}$ component of $I_U$.  Hence, the full composition in (4.2) identifies the $\hat{\iota}_{_U}$ component in $I_U$ with itself; it is multiplication by $+1$.  Since this is the case for all of the various $\eta_{_{UU^{\prime}U^{\prime\prime}}}$, the degree 2 obstruction class vanishes.

\subsection{No universal Hilbert space bundle}

	The observation in Section 4.1 regarding the non-existence of a universal line bundle over $\ctwon$ implies in turn that there is no Hilbert space bundle over $\ctwon$ whose fiber over $\fp$ is the corresponding $\hp$.  To see this, suppose for the sake of argument that such a bundle exists (call it $\mathbb{H}$).  Given $\mathbb{H}$, one could construct (using a partition of unity) a finite set of sections of $\mathbb{H}$ to be denoted by $\{s_1, \ldots,s_N\}$ with the following properties:  \\
\begin{itemize}
\item \it	For any given $\fp \in \ctwon$ each section from $\{s_1, \ldots,s_N\}$ at $\fp$ is a smooth element in the corresponding $\hp$.
\item Given $\fp \in \ctwon$ and $x \in S^2-\fp$, there exists $s \in\{s_1, \ldots,s_N\}$ with $s|_{\fp}$ being non-zero at $x$. 
\end{itemize}
\hfill (4.3)\\
These conditions imply that the assignment of $(\fp, x)$ to the $N$-tuple $(s_1|_{\fp}(x), \ldots, s_N|_{\fp}(x))$ defines a smooth map from $\mathcal{F}$ to $\mathbb{RP}^N$.  Let $\mathcal{I}$ denote the pull-back by this map of the tautological line bundle over $\mathbb{RP}^N$.  The restriction of $\mathcal{I}$ to each $\fp\in\ctwon$ version of $\cf|_{\fp}$ is necessarily isomorphic to $\ip$.  The existence of $\mathcal{I}$ and, by implication, the existence of $\mathbb{H}$, runs afoul of the non-existence result from Section 4.1.  

By way of a summary:  The obstruction to the existence of a universal $\mathbb{H}$ bundle is the obstruction class in $H^2(\ctwon)$ that was described in the preceding subsection.

\subsection{A universal $\rp^\infty$ bundle}

	Although there is no universal Hilbert space bundle, there is a universal $\rp^\infty$ bundle over $\ctwon$ whose fiber over any given $\fp \in \ctwon$ is the quotient by $\pm 1$ of the subspace $\mathbb{S}_p \subset\hp$ that is defined by the rule whereby $f\in\mathbb{S}_p$ if and only if 
\begin{equation}
\int_{S^2} f^2 = 1.
\tag{4.4}
\end{equation}
The $\rp^\infty$ bundle will be denoted by $\rp$.  

The topology on this bundle is generated by a certain basis of open neighborhoods which is described momentarily.  To set notation for this:  Suppose for the moment that $\fp$ is a given point in $\ctwon$.  Section 4.1 introduced the notion of a polydisk neighborhood of $\fp$, denoted there by $U_{\fp}$.  By way of a reminder, the polydisk neighborhood $U_{\fp}$ is the projection to $\ctwon$ of a product of disks in $\times_{2n}S^2$ with these disks centered at the constituents of $\fp$ and chosen so that the set of concentric disks with 100 times the radius have pairwise disjoint closure.  
(The disks in this product can have different radii.)
A convention for what follows:  The radii of these disks are constrained by an upper bound that guarantees the existence of a disk of radius ${1\over 100n}$  that is disjoint from each of them.

For the purpose of what is to come, the collection of disks that define $U_{\fp}$ is denoted by $\{D_{0k}\}_{k=1,\ldots,2n}$ and the set of respective 100 times larger radius concentric disks is denoted by $\{D_{1k}\}_{k=1,\ldots,2n}$.   An important point to keep in mind:  If $\fq$ is a configuration from $U_{\fp}$, then the bundles $\mathcal{I}_{\fq}$ and $\ip$ are isomorphic on $S^2 - \cup_{k=1,\ldots ,2n} D_{1k}$.  (This is because they have the same first Stieffel-Whitney classes on this set.)  Also keep in mind that there are precisely two isometries between these bundles on $S^2-\cup_{k=1,\ldots,2n}D_{1k}$, one is obtained from the other by fiber-wise multiplication by $-1$.

By way of more notation:  Supposing that $\fp$ again denotes a configuration in $\ctwon$, and supposing that $f$ is from $\mathbb{S}_{\fp}$, then the \{$\pm 1$\} equivalence class of $f$ is denoted in what follows by $[f]$. 

With preceding notation understood, suppose that $(\fp, [f])$ is as just described.  A basis for the open neighborhoods in $\rp$ of $(\fp,[f])$ is given by sets labeled by data consisting a polydisk neighborhood in $\ctwon$ of $\fp$ which is denoted by $U_{\fp}$, an open set in $S^2-\cup_{k=1,\ldots, 2n}D_{1k}$ which is denoted by $\Omega$ containing a disk of radius  ${1\over 100n}$, and positive numbers $\delta$ and $E$ with $E$ greater than 
the integral that appears on the right hand side of (1.7), 
the $S^2$ integral of $|df|^2$.  The corresponding neighborhood in $\rp$ consists of the set of pairs $(\fq, [h])$ with $\fq \in U_{\fp}$, and with $h \in \mathbb{S}_{\fq}$ satisfying two constraints:

\begin{itemize}
\item \it $\displaystyle\int_{S^2-\fp} |dh|^2 < E$.
\item There is an isometry between $\ip$ and $\mathcal{I}_{\fq}$ on $S^2-\cup_{k=1,\ldots,2n} D_{1k}$ that identifies $h$ as a section of $\ip$ obeying $\displaystyle\int_{S^2-\cup_{k=1,\ldots,2n} D_{1k}}( | d (f-h)|^2 + |f-h|^2 )   < \delta$.
\end{itemize}
\hfill (4.5)

\normalfont
The next lemma makes the formal assertion to the effect that $\rp$ fibers over $\ctwon$.\\

\noindent \bf Lemma 4.2:  \it The space $\rp$ has the structure of a (topological) fiber bundle over $\ctwon$ with the projection map sending any given pair $(\fp,[f])$ to $\fp$.\\

\rm 
A notion of a differentiable structure is introduced in the next subsection.  \\

\noindent \bf Proof of Lemma 4.2: \rm  The proof of the lemma has three parts. The proof requires the specification of a local product structure and then the verification that the transition functions between overlapping open sets obey the required cocycle condition. \\

\it	Part 1:  \rm A product structure for $\rp$ over a given set $U_{\fp}$ in $\ctwon$ can be defined using the map $\phi_{\fp}$ from Sections 2.6 and 4.1 as follows:  If $\fq \in U_{\fp}$, then pull-back by $\phi_{\fq,\fp}$ defines an invertible, linear map from the Hilbert space $\mathbb{H}_{\fq}$ to the Hilbert space of sections of $\phi_{\fq,\fp}^*\mathcal{I}_{\fq}$ 
which preserves the integral in (4.4).  Composing this linear map with any given isometry from $\phi_{\fq,\fp}^* \mathcal{I}_{\fq}$ to $\ip$ defines an invertible, linear map from $\mathbb{S}_{\fq}$ to $\ssp$.  Choosing a different isometry from $\phi_{\fq,\fp}^* \iq$ to 
$\ip$ multiplies this linear map by $\pm 1$.  As a consequence, there is no ambiguity if pull-back by $\phi_{\fq,\fp}$ is viewed as a map from $\ssq/\{\pm 1\}$ 
to $\ssp/\{\pm 1\}$. 
 As $\fq$ varies in $U_{\fp}$, the map $\phi_{\fp}$ gives a fiberwise identification between $\rp|_{U_{\fp}}$ and $U_{\fp} \times \rp|_{\fp}$.
As explained in Parts 2 and 3 of the proof, this identification is a homeomorphism of topological spaces when $U_{\fp} \times \rp|_{\fp}$ is given the product topology using the quotient topology for $\rp|_{\fp}$ from the latter's identification with $\ssp/\{\pm 1\}$ and with $\ssp$ viewed as a submanifold in $\hp$.  Accept that this is the case for what is said in the next paragraph.

The local product structure defines a fiber bundle topology on $\rp$ if and only if the corresponding transition functions obey the required cocycle condition on triple intersections (these being 
$U_{\fp}\cap U_{{\fp}^{\prime}}\cap U_{{\fp}^{\prime\prime}})$.  To see this, suppose for the moment that $\fq\in U_{\fp}\cap U_{{\fp}^{\prime}}$.  On the one hand, $\rp|_{\fq}$ has been identified with $\rp|_{\fp}$ by the rule 
\begin{equation}
[f] \to [\phi_{\fp,\fq}^*f].
\tag{4.6}
\end{equation}
(Keep in mind that pull-back commutes with multiplication by $\pm 1$.)  There is a similar formula for the identification between $\rp|_{\fq}$ and $\rp|_{\fp^{\prime}}$.  Composing the inverse of the one depicted in (4.6) with its primed version leads to the $\fq$-dependent linear, invertible map from $\rp|_{\fp}$ to $\rp|_{\fp^{\prime}}$ that is given by the rule 
\begin{equation}
[f] \to [(\phi_{\fq,\fp^\prime} \circ \phi_{\fp,\fq}^{-1})^*f].
\tag{4.7}
\end{equation}

Supposing that $\fq$ is in a triple intersection, $U_{\fp}\cap U_{{\fp}^{\prime}}\cap U_{{\fp}^{\prime\prime}}$, then it follows from the three versions of (4.7) that the relevant cocycle for this triple intersection when evaluated at $\fq$ is the identity map (this is because $(\phi_{\fq,\fp}\circ \phi^{-1}_{{\fq,\fp}^{\prime\prime}}) \circ ( (\phi_{{\fq,\fp}^{\prime\prime}}  \circ \phi_{{\fq,\fp}^{\prime}}^{-1} ) \circ (\phi_{{\fq,\fp^{\prime}}} \circ \phi_{\fq,\fp}^{-1}))$ is the identity map.)  \\

\it	Part 2:  \rm It remains now to prove that the identification between $\rp|_{U_{\fp}}$ and $U_{\fp}\times \rp|_{\fp}$ given by (4.6) defines a homeomorphism between the two topological spaces.  The proof that such is the case is facilitated by first reinterpreting the Hilbert space topology on any given $\fp \in \ctwon$ version of $\hp$.  This part of the proof does that.

 	Consider the topology on $\hp$ whereby a basis for the open neighborhoods of a given element $f$ is labled by an open neighborhood with compact closure in $S^2-\fp$ that contains a disk of radius  ${1\over 100n}$ (to be denoted by $\Omega$); and positive numbers $\delta$ and $E$ with $E$ greater than the $S^2$ integral of $|df|^2$.  The corresponding open set consists of the subset of elements $h \in \hp$ that obey the top bullet in (4.5), and also obey the $\fq = \fp$ version of the second bullet in (4.5) with the isometry between $\ip$ and itself being multiplication by $+1$.  This basis set is denoted by $\mathcal{V}_f(\Omega,\delta,E)$.  
	
As explained directly, this new topology is equivalent to the Hilbert space topology that is defined by the $\|\cdot\|_{\mathbb{H}}$-metric balls (the norm $\|\cdot\|_{\mathbb{H}}$ is depicted in (1.1)).  Indeed, the set $\mathcal{V}_f(\Omega,\delta,E)$ is open with respect to the Hilbert Space topology because if $h$ is in $\hp$ and if $\|f-h||_{\mathbb{H}}$ is sufficiently small (given $\delta$ and $E$), then both the conditions in (4.5) will hold for $h$. 

To see that the Hilbert space balls are open in this new topology, fix $\epsilon>0$ and consider in particular the ball in $\hp$ where $\|f - (\cdot)\|_{\mathbb{H}} <\epsilon$.  This ball will be open in the new topology if each of its elements is contained in a suitable $\mathcal{V}_{(\cdot)}(\cdot)$ set from the new topology which is entirely inside the $\|f -(\cdot)\|_{\mathbb{H}}< \epsilon$ ball.  In this regard:  No generality is lost by proving that there is a $\mathcal{V}_f(\cdot)$ set inside the $\|f - (\cdot)\|_{\mathbb{H}} < \epsilon$ ball.  To do this, fix $\delta > 0$ for the moment and then an open set with compact closure in $S^2-\fp$ (denoted by $\Omega$) such that

\begin{equation}
\int_{S^2-\Omega} |df|^2 < \delta.
\tag{4.8}
\end{equation}
Let $\mathcal{E}$ denote the $S^2$ integral of $|df|^2$.  With these definitions in hand, suppose that $h$ is from the set $\mathcal{V}_f(\Omega,\delta,E=\mathcal{E}+\delta$).  
(These basis sets for the topology are defined at the start of Section 2.)
If this is the case, then 

\begin{itemize}
\item $\displaystyle\int_\Omega (|d(f-h)|^2 + |f-h|^2 ) < \delta$.
\item $\displaystyle\int_{S^2 - \Omega} |d(f-h)|^2 < 2 \int_{S^2 - \Omega} (|df|^2 + |dh|^2) <  2\delta + 2 \int_{S^2 - \Omega} |dh|^2$.
\end{itemize}
\hfill (4.9)

To exploit this:  By virtue of the top bullet in (4.9), the $\Omega$ integral of $|dh|^2$ is greater than $\mathcal{E} - 2\sqrt{ \delta}\sqrt {\mathcal{E}}$.  
Therefore, because the integral of $|dh|^2$ over $S^2 - \fp$ is at most $\mathcal{E} + \delta$ (from the definition of $\mathcal{V}_f(\Omega, \delta, \mathcal{E} = \mathrm{E} + \delta))$, its integral over $S^2 - \Omega$ must obey

\begin{equation}\int_{S^2 - \Omega} |dh|^2  \leq \delta +2\sqrt{ \delta} \sqrt{\mathcal{ E}}.
\tag{4.10}
\end{equation}

And thus, the left hand side of the inequality in the second bullet of (4.9) is no greater than $4\delta+4\sqrt \delta \sqrt \mathcal{E}$.  This last observation plus the top bullet in (4.9) implies the following:  If $\delta$ is sufficiently small given $\mathcal{E}$ (and an appropriate $\Omega$ is chosen to obey (4.8)), then the $S^2$ integral of $|d(f-h)|^2$ will be less than $\epsilon\delta^{1\over 3}$.  Meanwhile, the latter integral plus the $\Omega$ integral of $|f-h|^2$ bounds a positive multiple of the $S^2$ integral of $|f-h|^2$ (the multiple is independent of $\Omega$ because of the constraint that $\Omega$ must contain a disk of radius ${1\over 100n}$).  With the preceding understood, it follows that the $\| \cdot \|_{\hi}$ norm of $f-h$ will be less than $\epsilon$ if $\delta$ is sufficiently small and if $\Omega$ is chosen so that (4.8) holds.  This says in effect that the set $\mathcal{V}_f(\Omega, \delta, E = \ce +\delta)$ is entirely inside the $\|\cdot \|_{\hi}$ radius $\epsilon$ ball in $\hp$ centered at $f$.\\

\it Part 3:  \rm Return now to the assertion that pull-back using the map $\phi_{\fp}$ defines a homeomorphism between $\rp|_{U_{\fp}}$ and $U_{\fp} \times \rp|_{\fp}$.  To do this, it is sufficient to prove that a basis set in $\rp|_{U_{\fp}}$ corresponds to an open set in $U_{\fp}\times \rp|_{\fp}$ and vice-versa.  That this is so follows by unwinding definitions from the analysis in Part 2 and from the three observations that follow.  

The first observation is that a pair consisting of a configuration $\fq \in U_{\fp}$ and an element $[f] \in \rp|_{\fp}$ is identified by the local product map with the element $(\fq, [(\phi_{\fq,\fp}^{-1})^*f])$ in $\rp$ (keep in mind that $[(\phi_{\fq,\fp}^{-1})^*f]$ is in $\rp|_{\fq}$).  More to the point, the transformation given by the rule $f \to( \phi_{\fq,\fp}^{-1} )^*f$ defines a norm-bounded linear map from $\hp$ to $\mathbb{H}_{\fq}$ given an isometry between $(\phi_{\fq,\fp}^{-1})^*\ip$ and $\iq$.  Conversely, if $[h] \in \rp|_{\fq}$, then $\phi^*_{\fq,\fp} h$ is in $\hp$ given an identification between $\phi^*_{\fq,\fp} \iq$ and $\ip$.
  
The second observation is this:  If $\Omega$ has compact closure in $S^2-\fp$, then $\phi_{\fq,\fp}(\Omega)$ has compact closure in $S^2-\fp$.

To state the third observation, suppose that $f \in \hp$.  Then (as explained in Section 2.6), a change of coordinates identifies the $S^2$ integral of $|d(\phi_{\fq,\fp}^{-1})^*f|^2$ with an $S^2$ integral that has the form 

\begin{equation}
\int_{S^2 - \fp} \langle df, df \rangle_{\fm_{\fq}}
\tag{4.11}
\end{equation}
with $\langle , \rangle_{\fm_{\fq}}$  denoting the inner product of $df$ with itself as defined by the $\phi_{\fq,\fp}$-pull back of the round metric.  The point is that this integral defines a bounded quadratic form on $\hp$ that depends smoothly on $\fq$.  Moreover, the value of the integral in (4.11) differs from the corresponding integral of $|df|^2$ by at most $c_0 \dist(\fq,\fp)$ times the latter integral.  Thus, if the $S^2$ integral of $|df|^2$ is less than a given number (call it $E$), then the $S^2$ integral of $|d(\phi_{\fq,\fp}^{-1})^*f|^2$ will be bounded by $(1+c_0 \dist(\fq, \fp))E$.  

\subsection{On the differentiability of maps to and from $\rp$}

	The notion of a differentiable map to or from $\rp$ is problematic because the fiber-bundle structure on $\rp$ was defined in the topological category but not in the $C^1$ (let alone smooth) category.  This was done because the transition function between two product neighborhood charts, versions of $U_{\fp} \times\rp_{\fp}$ and $U_{{\fp}^{\prime}}\times \rp|_{{\fp}^{\prime}}$, are not differentiable with respect to variations of the base configuration $\fq \in U_{\fp}\cap U_{{\fp}^{\prime}}$.  To elaborate:  Differentiability fails because the transition functions (depicted in (4.7)) involve the composition of fiber elements with diffeomorphisms that depend on the configuration $\fq$.  And, the problem with this dependence is that the derivative with respect to variations in $\fq$ of the fiber-wise linear transition map depicted in (4.7) when acting on a given element $f\in\rp|_{\fp}$ costs one derivative of $f$ (this is the chain rule).  As a consequence, the derivative of the linear transition map when viewed from the context of $\hp$ does not define a bounded linear operator from $\hp$ to $\mathbb{H}_{\fq}$.  It is only bounded as a linear map from $\hp$ to the $L^2$ completion of $\mathbb{H}_{\fp^{\prime}}$.  By way of an example:  If $h \in\hp$, then the linear function

\begin{equation}
f \to \int_{S^2} hf
\tag{4.12}
\end{equation}
pulls back over $U_{\fp}\cap U_{{\fp}^{\prime}}$ via the product structure transition functions to a $\fq$-dependent linear functional on $\hpp$ which is differentiable with respect to $\fq$ (but not twice differentiable).

	One can also define the notion of a \it smooth \rm function on $\rp$:  This is a function that is smooth when written using any of the product charts $\{U_{\fp} \times \rp|_{\fp}\}_{\fp\in C_{2n}}$.  The next subsection gives the most relevant example of a smooth function.  An example in the meantime is defined by the rule

\begin{equation}
f \longrightarrow \int_{S^2} Gf^2
\tag{4.13}
\end{equation}
with $G$ being a smooth function on $\ctwon$.  

With regards to maps into $\rp$:  A map into a product chart $U_{\fp}\times\rp|_{\fp}$ is given by a pair $(\fq(\cdot),[f(\cdot)])$ of maps into $U_{\fp}$ and into $\rp|_{\fp}$ respectively.  Although both might be smooth, if $\fq(\cdot)$ maps to $U_{\fp}\cap U_{{\fp}^{\prime}}$,
 then the image of this map via the bundle transition function need not be differentiable as a map into $U_{\fp}\cap U_{{\fp}^{\prime}} \times \rp|_{{\fp}^{\prime}}$.  (This can happen even when the $[f(\cdot)]$ component is constant.)  However, such a map will be differentiable (at least of class $C^1$)
 if the $[f(\cdot)]$ component over $U_{\fp}$ lifts to a map from the domain space into a certain dense subvector space in $\hp$.  This vector subspace is the closure of the space of smooth sections of $\ip$ with compact support in $S^2-\fp$ using the norm whose square is defined by the rule
\begin{equation}
f \longrightarrow \int_{S^2 - \fp} (|\nabla df|^2 + |df|^2 + |f|^2).
\tag{4.14}
\end{equation}

In this regard:  The fiber-wise linear bundle transition map at each $\fq \in U_{\fp}\cap U_{{\fp}^{\prime}}$ sends this dense domain in $\hp$ in a 1-1 and onto fashion to the $\fp^\prime$ version of this domain in $\hpp$.  And in so doing, the $\fq$-derivatives of the bundle transition maps give a bounded linear operator from the dense domain to $\hpp.$  Said differently:  The bundle transition functions for $\rp$ are differentiable when restricted to the vector subspace in $\hp$ where the integral depicted in (4.14) is finite.

\subsection{The function $\mathcal{E}$}

Let $\mathcal{E}$ denote the function on $\rp$ whose value at any given $(\fp,[f])$ is this:

\begin{equation}
\mathcal{E}(\fp,[f]) \equiv \int_{S^2-\fp} |df|^2.
\tag{4.15}
\end{equation}
The upcoming Lemma 4.3 uses $\mathcal{E}$ to reinterpret part of what is said by Proposition 2.7.  It says in effect that $k$-critical points in $\ctwon$ in the sense of Definition 3.7 are the $\pi$-images of the critical points of $\mathcal{E}$ on $\rp$. (In this regard:  With $\mathcal{E}$ viewed as a function on $\rp$ as is done here and below, the critical points of $\mathcal{E}$ are critical points subject to the constraint that the integral of $f^2$ on $S^2-\fp$ is equal to 1.) \\

\noindent \bf Lemma 4.3: \it  The function $\mathcal{E}$ from (4.15) is a smooth function on $\rp$.  As such, its critical points are precisely the set of \{$\pm 1$\} equivalence classes of pairs $(\fp, f)$ such that $\fp \in \ctwon$ and $f \in \ssp$ is an $\ip$-eigensection with all numbers from Lemma 2.3's set $\{n_{p}(f)\}_{p\in \fp}$ being positive.\\

\noindent\bf Proof of Lemma 4.3:  \rm By definition, the function $\mathcal{E}$ is smooth on $\rp$ if and only if it is smooth when viewed on the restriction of $\rp$ to any given $\fp \in \ctwon$ version of $U_{\fp}$ with $\rp$ identified there as $U_{\fp}\times(\mathbb{S}_{\fp}/\{1\})$.  To see about that, suppose that $\fq \in U_{\fp}$ and $ [f] \in \mathbb{S}_{\fp}/\{\pm 1 \}$. The function $\mathcal{E}$ on $(\fq,[f])$ has the form depicted in (4.11) which is observedly jointly smooth with respect to variations in $\fq$ and $f$ (see the discussion in Section 2.4).

	With regards to the critical points:  The function is critical at $(\fp, [f])$ with respect to variations in the fiber direction iff and only if $f$ is an $\ip$-eigensection (this follows by definition).  The claim that it is also critical with respect to variations in $\fp$ if and only if all numbers from the set $\{n_p(f)\}_{p\in\fp}$ are positive follows from Proposition 2.5 and the formula in (2.21).  \\

	If $\mathcal{E}$ were a Palais-Smale function on $\rp$, which is to say proper, then its critical points would account for topological invariants of $\rp$ in a Morse theoretic sense (see e.g \cite{PS} and the review \cite{MW}).  Although the function $\mathcal{E}$ is proper along the fibers of $\rp$, it is not globally proper because $\ctwon$ is not compact.  As a consequence, it could be that $\mathcal{E}$ does not have any critical points notwithstanding the fact that $\rp$ has plenty of homology/cohomology.  (The cohomology of $\rp$ is described in the next subsection.)   What follows are some comments about this.

  	With regards to topology and critical points:  Let $C$ denote a non-zero class in $H_*(\rp)$.  Set $\mathcal{E}_C$ to denote the infimum over all singular cycles representing the class $C$ of the maximal value of $\mathcal{E}$ on the simplices that appear with non-zero coefficient in the cycle.  This number would be a critical value of $\mathcal{E}$ were $\mathcal{E}$ a proper function.  Since homology and cohomology in any given degree are dual to each other using $\integer/2\integer$ coefficients (all homology/cohomology here and in what follows uses $\integer/2\integer$ coefficients), there is also a cohomological point of view:  
	Letting $C$ denote a non-zero cohomology class, define $\mathcal{E}_C$ to be the infimum of values of $E\in(0,\infty)$ such that the class $C$ is non-zero on the $\mathcal{E} < C$ part of $\rp$.  This version of $\mathcal{E}_C$ would be a critical value of $\mathcal{E}$ 
	were $\mathcal{E}$ proper on $\rp$.   Homotopy invariants of $\rp$ also lead to putative critical values for $\mathcal{E}$:  Take the infimum over a homotopy class of maps from a compact space 
	to $\rp$ of the maximum value of $\mathcal{E}$ on the image of the map.  (This last version of a min-max value is used in Section 6.)  
	
Unfortunately $\mathcal{E}$ is not proper and therefore none of these hypothetical critical values of $\mathcal{E}$ must correspond to an actual critical points of $\mathcal{E}$.  Examples of this pathology are provided in Section 5.2.

The homology/cohomology of $\rp$ is the subject of the next subsection.

	\subsection{The homology/cohomology of $\rp$}

The upcoming Lemma 4.4 describes some non-zero classes in the $H^*(\rp)$.  Since the lemma refers to classes in $H^*(\ctwon)$, a digression comes first to describe the latter. 

The computations in \cite{Sc} (see Corollary 3.1) imply that $H^*(\ctwon)$ can be depicted as follows:  Let $\omega$ denote the obstruction class in $H^2(\ctwon)$ from Sections 4.1 and 4.2.  Then $H^*(\ctwon)$ can be written as	
\begin{equation}
H^*(\ctwon) \cong H^*(\mathbb{C}_{2n}) \oplus (\omega \wedge H^{*-2}(\mathbb{C}_{2n-1}))
\tag{4.16}
\end{equation}
where $\mathbb{C}_k$ denotes the space of $k$ unordered points in $\mathbb{C}$.  The identifications in (4.16) come from the splitting of the long exact sequence for the pair $(\ctwon,\ctwon-\mathcal{C}^x_{2n})$ with $x$ 
  being any chosen point in $S^2$ (and recall that $\mathcal{C}_{2n}^x$ denotes the set of configurations in $\ctwon$ that do not contain $x$). In terms of homology, (4.16) asserts in effect that each class in $H_*(\ctwon)$ can be written uniquely as a sum of a class from $\ctwon-\mathcal{C}^x_{2n}$ and a class with non-zero cap-product pairing with $\omega$.  With regards to $H^*(\mathbb{C}_k)$:  These modules are computed by \cite{Fu}, they are zero above a degree that is on the order of $n$ and, in general, non-zero in most degrees below the maximal one.

	The cohomology of $\rp$ can be computed using the fact that this space is a fiber bundle over $\ctwon$ with the infinite dimensional real projective space as fiber.  The next lemma summarizes.  The lemma refers to the class $\omega$ that appears in (4.16).  Keep in mind that $\omega\wedge\omega$ is zero in cohomology. \\

\noindent \bf Lemma 4.4: \it The cohomology of $\rp$ is (non-canonically) isomorphic to the vector space of polynomials in a degree 2-class $\tau$ with coefficient ring $\mathcal{K}$ as described in the subsequent bullets.  This is to say that 
$$H^*(\rp) = \mathcal{K}^* \oplus ( \mathcal{K}^{*-2} \wedge \tau )\oplus ( \mathcal{K}^{*-4} \wedge \tau \wedge \tau)\oplus \cdots
$$

\begin{itemize}
\item The definition of $\tau$:  Let $\hat\omega$  denote a given representative cocycle for $\omega$.  Then $\pi^*\hat\omega$  can be written as $d\nu$ with $\nu$ being a 1-cochain whose restriction to each fiber algebraically generates the cohomology of the fiber.  The cochain $\nu\wedge \nu$  is closed and non-zero in $H^*(\rp)$.  The class of $\nu\wedge \nu$ is $\tau$.

\item The definition of $\mathcal{K}$:  This module can be written as
$$\mathcal{K}^* = \pi^* H^* (\mathbb{C}_{2n}) \oplus (\Theta \wedge \pi^*H^{*-3}(\mathbb{C}_{2n})) $$
where $\pi^*$ and $\Theta$ are as follows:
\begin{itemize}
\item[$\mathrm{a)}$]	What is denoted by $\pi^*$ signifies the pull-back by the projection map $\pi$ to $\ctwon$.  This pull-back homomorphism is injective on the $H^*(\mathbb{C}_{2n})$ summand in (4.16) and zero on the other summand.

\item[$\mathrm{b)}$]	 To define $\Theta$, let $\hat o$ denote a 3-cocycle on $\ctwon$ obeying $d\hat{o} = \hat{\omega} \wedge \hat{\omega}$.  The 3-cocycle $\nu\wedge \pi^* \hat{\omega}  + \pi^*\hat{o}$ on $\rp$ is closed and non-zero in $H^*(\rp)$.  Its class is independent of the choice for the cocycle $\hat{o}$.  This cohomology class is $\Theta$.

\end{itemize}

\end{itemize}

\rm

Some implications of Lemma 4.4 concerning the critical point set of $\mathcal{E}$ are derived in Section 5.2.\\

\noindent \bf Proof of Lemma 4.4:  
\rm This lemma is an instance of Proposition A.1 in the appendix.

\section{The extension of $\rp$ and $\mathcal{E}$ to $\ctwonb$}

This section describes the sense in which the space $\rp$ and its function $\mathcal{E}$ extend to the whole of $\ctwon$'s compactification $\ctwonb$ from (3.2).  To say more about this extension, write the space $\ctwonb$ as done in (3.1).  Then each $k \in\{0, \ldots, n\}$ version of $\mathcal{C}_{2k}$ has its corresponding version of $\rp$ which will be denoted here as $\rp_{2k}$.  (The $k = 0$ version, $\rp_0$, is the quotient by the multiplicative action of \{$\pm 1$\} on the subspace of functions in the standard $L^2_1$ Sobolev space on $S^2$ that obey (4.4).)  The extension of $\rp$ over $\ctwonb$ is denoted by $\rpb$ and it is (as a point set)

\begin{equation}
\rpb = \rp_{2n} \cup \rp_{2n-2} \cup \cdots \cup \rp_0.
\tag{5.1}
\end{equation}

This extension maps to $\overline{\mathcal{C}}$ in the obvious way (the map is denoted by $\pi$).   The topology on $\rpb$ is defined from a basis of neighborhood open sets that are analogous to those used in Section 4.3 for the topology on $\rp$.  The details appear momentarily.  But note in advance the space $\rpb$ does not meet the technical definition of a fiber bundle because there aren't local product neighborhoods for the points in $\ctwonb-\ctwon$.  (Even so, the map $\pi$ is continuous and surjective.)  There is, however, a weak notion of a local product structure for $\ctwonb$ which is described in the upcoming Section 5.3.  

\subsection{A topology for $\rpb$}

A convention concerning elements in $\ctwonb$ is used for the specification of topology.  The statement of the convention refers to the elements in $\ctwonb$ as $\integer/2\integer$ divisors as done in Section 3.1:  Supposing that $[\fq]$ is such a divisor, then its minimal representative is implicitly viewed as a configuration in the relevant $k \in\{0, \ldots, n\} $ version of  $\mathcal{C}_{2k}$.  This number $k$ is said to be the \it degree \rm of the divisor.

To define the topology on $\rpb$, first recall from Part 2 of Section 3.1 how the topology on $\ctwonb$ was defined by a basis of open neighborhoods of any given $\integer/2\integer$ divisor:  Let $k$ denote an integer from $\{0, \ldots, n-1\}$ and let $[\fq]$ denote a given degree $k$ divisor from $\ctwonb$.  Fix a positive number to be denoted by $\epsilon$ which is chosen to be very small (much less than  ${1\over 1000}$ times the distance between the points in $[\fq]$'s minimal representative).  The corresponding open neighborhood was denoted there by $\mathcal{U}(\epsilon,k,[\fq])$ and it consists of the configurations that obey the conditions in (3.3).  (Remember that the minimal representative of a divisor as depicted in (3.1) is the configuration $\fp$ with all coefficient numbers $\{m_p\}_{p\in\fp}$ equal to 1.)

With the preceding in mind, let $\fq \in \mathcal{C}_{2k}$ denote the minimal representative of the divisor $[\fq]$, and let $[f]$ denote an equivalence class in $\rp_{2 k}|_{\fq}$.  A basis set for the open neighborhoods of the pair $([\fq],[f])$ in $\rpb$ are labeled by positive numbers $\epsilon, \delta$ and $\mathcal{E}$ with $\epsilon$ less than ${1\over 1000 n}$  and with $E$ greater than the $S^2$ integral of $|df|^2$.  This basis set consists of the pairs $([\fc],[h]) \in\rpb$  that obey the following constraints: \\
\it
\begin{itemize}
\item	$[\fc] \in \mathcal{U}(\epsilon,k,[\fq])$.
\item $\displaystyle \int_{S^2 - \fc} |dh|^2	 < E$.
\item Let $\fc$ denote the minimal divisor of $[\fc]$ and let $\Omega$ denote the $\dist( \cdot, \fc) > 200\epsilon$ part of $S^2$.  There is an isometry between $\mathcal{I}_{\fc}$ and $\mathcal{I}_{\fq}$ on $\Omega$ that identifies $h$ with a section of $\iq$ obeying $\displaystyle \int_\Omega (|d(f-h)|^2 + |f-h|^2)  < \delta$.
\end{itemize}\rm
\hfill (5.2)\\

It is left to the reader to verify that this topology restricts to each stratum in (5.1) to give the topology that is described in Section 4.3 for the relevant $k \in\{0, \ldots, 2n\}$ version of $\rp_{2k}$.  It is also left to the reader to verify that the map $\pi$ from $\rpb$ to $\ctwonb$ is continuous with this topology on $\rpb$.

\subsection{The extension of $\mathcal{E}$}
The lemma that follows makes a formal assertion to the effect that the function $\mathcal{E}$ extends as a continuous function to the whole of $\rpb$.\\

\noindent \bf Lemma 5.1:  \it The function on $\rpb$ that is defined by the rule 
$$([\fq],[f])  \to \mathcal{E}([\fq],[f]) \equiv  \int_{S^2-\fq} |df|^2$$
is continuous in general and smooth along any stratum in (5.1).  \\

\bf \noindent Proof of Lemma 5.1:  \rm The assertion that $\mathcal{E}$ is smooth along the strata of $\rpb$ are instances of Lemma 4.3.  The proof that $\mathcal{E}$ is continuous at a given pair $([\fq],[f])$ invokes the following observation:  \\

\begin{center} \it
Given a positive number $\mu$  at most ${1\over 1000n}$, there exists $\epsilon> 0$ such that the integral of $|df|^2$ on the 
union of any $n$ disks of radius $\epsilon$ is at most $\mu$.
\end{center}
\hfill (5.3)\\
\noindent By way of a proof:  There is a smooth section of $\iq$ with compact support in $S^2-\fq$ to be denoted by $f_\mu$ with the $S^2$ integral of $|d(f-f_\mu)|^2+|f-f_\mu|^2$ being less than  ${1\over 100}\mu$  .  (This follows from the definition of $\mathbb{H}_{\fq}$.)  The assertion in (5.3) holds for $f_\mu$ because $f_\mu$ is smooth and supported away from $\fq$.  Thus, it holds also for $f$.

	With (5.3) in hand, then the proof of Lemma 5.1 can be completed by copying with only notational changes what is said by (4.8)-(4.10) and what is said by the surrounding discussion in the last paragraph of Part 2 of Lemma 4.2's proof.\\
	
	To motivate the subsequent constructions and by way of an application of Lemmas 4.4 and 5.1, return to the context of the min-max values for $\mathcal{E}$ on classes in $H_*(\rp)$ as defined in the last paragraph in Section 4.5.  It follows from Lemma 5.1 that these min-max value for $\mathcal{E}$ are zero on the $\pi$-pull back of classes from the $H^*(\complex_{2n})$ summand in (4.16).  This is because there are representative cocycles for the duals of these classes in $H_*(\ctwon)$ that are entirely in $\ctwon$'s intersection with any given open neighborhood of the $\mathcal{C}_0$ stratum in $\ctwonb$.  
To explain:  These are classes from the base space $\ctwon$ that can be represented by continuously parametrized families of configurations that all lie in any given disk in the sphere.  In particular, the disk can have arbitrarily small radius.  If the radius is small enough, then the whole family can be lifted to $\rp$ via a section that assigns to each configuration in the family the (unique) eigenfunction with the smallest eigenvalue (which is nearly zero).  The eigenvalue is nearly zero and the eigenfunction is unique because all of the points in all of the configurations that comprise the family are very close together (near the $\mathcal{C}_0$ stratum in $\ctwonb$).  	
Said differently, the duals of these classes are sent to zero by the inclusion induced homomorphism from $H_*(\ctwon)$ into $H_*(\ctwonb)$.   By the same token, with $\tau$ and $\Theta$ as defined in Lemma 4.4, and supposing that $m$ is a non-negative integer, then the min-max values for $\mathcal{E}$ on the cup product of either $\wedge^m\tau$ or $\wedge^m\tau \wedge \Theta$ with the $\pi$-pull back of any class in the $H^*(\complex_{2n})$ summand in (4.16) is the same as the min-max value of $\mathcal{E}$ on just the class $\wedge^m\tau$ or $\wedge^m\tau \wedge \Theta$ the case may be.  These min-max values aren't known, but they are, in any event, no greater than their min-max values on $\rpb^2$.   In any event, a corollary of Proposition 3.3 is this:  Having fixed an integer $m$, there exists an integer $n_m$ with the following significance:  \\

\begin{center}\it
If $n\geq n_m$, then the min-max value of $\mathcal{E}$ for either $\wedge^m\tau$ or $\wedge^m\tau\wedge\Theta$ on the $\ctwon$ version 
of $\rp$ is the min-max value of $\mathcal{E}$ for these classes on some stratum of $\rpb-\rp$.
\end{center}
\hfill(5.4)\\
The point of the preceding is that $\rpb$'s cohomology is the relevant cohomology with regards to the relation between $\mathcal{E}$'s min-max values on $\rp$ and $\mathcal{E}$'s critical points on $\rp$.  

\subsection{$\rpb$ as a weak projective space bundle over $\ctwonb$}

The upcoming Lemma 5.3 describes a weak $\rp^\infty$ bundle structure for the space $\rpb$.   This notion of a weak $\rp^\infty$ bundle is defined in the context where there are topological spaces $\crr$ and $X$ with a surjective map $\pi: \crr \to X$ whose fibers are homeomorphic to $\rp^\infty$.  A weak $\rp^\infty$ bundle structure for $(\crr, X, \pi)$ consists of the following data.\\
\vskip -3mm
\begin{itemize}\it
\item For every non-negative integer $k$:  A locally finite open cover of $X$ (to be denoted by $\mathfrak{U}_k$) such that each $U \in \mathfrak{U}_k$ has a vector bundle $V_U \to U$ of rank greater than $k$ with a fiberwise embedding $\Psi_U$ from $(V_U-0)/\mathbb{R}^*$ to $\crr|_U$ that covers the projections to $U$
\begin{center}
\begin{tikzcd}[column sep=scriptsize](V_U - 0)/\real^* \arrow[dr,"\pi"] \arrow[rr,"\Psi_U"]& & \mathcal{R}|_U \arrow[dl,"\pi"'] \\& U \end{tikzcd}
\end{center}
and induces an isomorphism on the first homology of each fiber.

\item For each pair of positive integers $k$ and $k^\prime$:  If $U \in \uk$ and $U^{\prime} \in \mathfrak{U}_{k\prime}$ intersect, and supposing that $\dim(V_{U^{\prime}}) \geq \dim(V_U)$, then $\Psi_U$ factors through $\Psi_{U^{\prime}}$ over $U\cap U^{\prime}$ via an injective bundle map from $V_U|_{U\cap U^{\prime}}$ to $V_{U^{\prime}}|_{U\cap U^{\prime}}$.
\item Each $\integer/2\integer$ cohomology class of $\crr$ is non-zero when restricted to some sufficiently large $k$ version of $\crr$'s subspace  $\bigcup_{U\in\uk}\Psi_U((V_U-0)/\real^*)$.
\end{itemize}
\hfill (5.5)\\

A preliminary lemma is needed in order to define the weak $\rp^\infty$ structure for $\rpb$.  To set the stage:  The preliminary lemma (Lemma 5.2) refers to a $\integer/2\integer$ divisor from a stratum of $\ctwonb$ (it denotes the latter by $[\fq]$) and its minimal representative (denoted by $\fq$) which is a configuration from some version of $\mathcal{C}_{2k}$ (for $k \in\{0, \ldots., n\})$.  With $\fq$ in hand, the lemma refers to a positive number that is not an $\iq$-eigenvalue.  Letting $E$ denote that number, Lemma 5.2 uses $\mathbb{H}_{{\fq},E}$ to denote the span in the Hilbert space $\mathbb{H}_{\fq}$ of the $\iq$-eigensections with eigenvalue less than $E$.  The lemma also refers to $L^2$-orthogonal projection on $\mathbb{H}_{\fq}$.  This is the orthogonal projection defined by the $L^2$ inner product:
\begin{equation}
(f,f^\prime) \to \int_{S^2} ff^\prime.
\tag{5.6}
\end{equation}

The instances of $L^2$-orthogonal projection in the lemma are to the subspaces $\mathbb{H}_{(\cdot),E}$.   \\

\noindent \bf Lemma 5.2:  \it Fix a degree $2k$ divisor in $\ctwonb-\ctwon$ to be denoted by $[\fq]$ and let $\fq$ denote the minimal representative.   Then fix a positive number to be denoted by $E$ which is not an $\iq$-eigenvalue, and fix a positive number to be denoted by $\delta$.  Granted this data, there exists an open neighborhood of $[\fq]$ in $\ctwonb$ with the following properties:  Denote the neighborhood by $\overline{U}$.  Supposing that $\fp$ is the minimal representative of a $\integer/2\integer$ divisor in $\overline{U}$, then $E$ is not an $\ip$-eigenvalue; and there exists a bounded linear map from $\mathbb{H}_{\fq}$ to $\hp$ and a bounded linear map from $\hp$ to $\mathbb{H}_{\fq}$ with the properties listed next (both linear maps are denoted by $L$).
\begin{itemize}  
\item Both versions of L are injective on the domain $\mathbb{H}_{(\cdot),E}$ subspace, and their composition with the $L^2$-orthogonal projection to the range $\mathbb{H}_{(\cdot),E}$ subspace is an isomorphism.
\item If $f$ is from the domain $\mathbb{H}_{(\cdot),E}$, then
\begin{itemize} 
\item[$\mathrm{a)}$]	$\displaystyle (1-\delta)\int_{S^2} |f|^2 \leq  \int_{S^2} |Lf|^2 < (1+\delta)\int_{S^2} |f|^2$,
\item[$\mathrm{b)}$]	 $\displaystyle \int_{S^2} |df|^2 - \delta \int_{S^2} |f|^2 \leq   \int_{S^2} |d(Lf)|^2 <  \int_{S^2} |df|^2 + \delta \int_{S^2} |f|^2$.
\end{itemize}
\end{itemize}
\rm
A note about $\overline{U}$:  The subspace $\overline{U}\cap \mathcal{C}_{2m}$ is empty when $m$ is less than $k$ (which is half of the number of points that comprise $\fq$) and it is non-empty otherwise.   

	This lemma is proved momentarily.  A weak $\rp^\infty$ structure for $\rpb$ can be constructed using Lemma 5.2 according to the rules laid out below in Lemma 5.3.  \\

\noindent \bf Lemma 5.3:  \it What follows describes a weak $\rp^\infty$ structure for $\rpb$.   Fix a positive integer $k$.  For any given $\integer/2\integer$ divisor $[\fq]\in\ctwonb$ first fix a positive number $E$ which is not an $\iq$-eigenvalue and with $\dim(\mathbb{H}_{{\fq},E}) > k$.  
Then fix some very small, positive number $\delta$ which is much less than the distance from $E$ to the nearest $\iq$-eigenvalue.  With the preceding data as input, use Lemma 5.2 to define a version of the open set $\overline{U}$ (which is denoted here by $\overline{U}_q$).  
Having done this for each element in $\ctwonb$, take a locally finite cover of $\ctwonb$ from the collection $\{\overline{U}_{\fq}: \fq\in\ctwonb\}$ and denote this cover by $\mathfrak{U}_k$.  The required vector bundle $V_{\overline{U}_{\fq}}$  over a set $\overline{U}_{\fq}$ from the cover $\mathfrak{U}_k$ is the bundle whose fiber over any $[\fp] \in \overline{U}_q$ is $\mathbb{H}_{\fp,E}$.  The embedding  $\Psi_{\overline{U}_{\fq}}$ for the first bullet of (5.5) is the tautological inclusion map.\\

\rm

\noindent This lemma is proved using the next lemma.  \\

By way of terminology for the upcoming lemma:  A parametrized family of maps from one subset of $\rpb$ to another is said to be \it fiber preserving \rm when the projection via $\pi$ to $\ctwonb$ is not changed by any map from the family. \\

\noindent \bf Lemma 5.4: \it Fix a positive number to be denoted by $\mathfrak{C}$ and let $\rpb^{\mathfrak{C}}$ denote the $\mathcal{E}< \mathfrak{C}$ part of $\rpb$.  If $k > 2\mathfrak{C}$, then there exists a fiber preserving deformation retract of $\rpb^{\mathfrak{C}}$ onto an $\rpb$ version of the subspace  $\bigcup_{U\in\mathfrak{U}_k}\Psi_U((V_U - 0)/\real^*)$.\\

\noindent\rm Lemma 5.4 is proved after Lemma 5.2.  \\

\noindent \bf Proof of Lemma 5.3:  \rm The conditions in the first two bullets of (5.5) follow directly from the construction.  The third bullet of (5.5) holds for $\rpb$ if any given (singular) homology class is represented by a closed chain whose image lies entirely in some $\rpb$ version of  $\bigcup_{U\in \mathfrak{U}_k}\Psi_U((V_U-0)/\real^*)$ for $k$ sufficiently large.   To establish this, note first that any given closed singular chain has compact support in $\rpb$ which implies that $\mathcal{E}$ has an upper bound on that chain.  Let $\mathfrak{C}$ denote this upper bound.  Granted this upper bound, use the deformation retract from Lemma 5.4 to obtain a homologous chain that lies in a  $\bigcup_{U\in \mathfrak{U}_k}\Psi_U((V_U-0)/\real^*)$ subspace of $\rpb$. \\

\bf \noindent Proof of Lemma 5.2:  \rm The construction of $L$ borrows much by way of strategy and notation from Parts 1-4 of the proof of Proposition 3.2.   There are three parts to the construction of $L$ and the verification of its properties.  \\

\it Part 1: \rm  Fix a positive number to be denoted by $\epsilon$ which is chosen to be very small (much less than ${1\over 1000n}$), and then reintroduce from Part 2 of Section 3.1 the set $\mathcal{U}(\epsilon,k,[\fq])$.  Let $\fp$ denote a configuration of $2j$ points in this set for some $j \in\{k, \ldots, n\}$.  Define the function $\chi_\epsilon$ from (2.2) using $\fp$, but denote it here by $\chi_{\fp,\epsilon}$.  Note in particular that the function is zero at all points from $\fq$ and that the line bundles $\ip$ and $\iq$ are isomorphic on the support of $\chi_{p,\epsilon}$. 

	By virtue of $\epsilon$ being small, there is a point to be denoted by $x$ where $\chi_{\fp,\epsilon}$ is equal to 1.  Use this point to define an isomorphism between the line bundles $\ip$ and $\iq$ on the support of $\chi_{\fp,\epsilon}$.   With this isomorphism in hand, the map $L: \mathbb{H}_{\fq} \to \hp$ is given by the rule
\begin{equation}f\longrightarrow \chi_{\fp,\epsilon}f.
\tag{5.7}
\end{equation}
This same formula likewise defines the corresponding version of $L$ sending $\hp$ to $\mathbb{H}_{\fq}$.  In either direction, the arguments for (3.5) and (3.6) can be repeated with only notational changes to see that
\begin{itemize}
\item $\displaystyle \int_{S^2} |d(\chi_{\fp,\epsilon}f)|^2 = \int_{S^2} |df|^2   + \fe_1 \int_{S^2} |f|^2$,
\item $\displaystyle \int_{S^2} |\chi_{\fp,\epsilon}f|^2 =(1+\fe_2) \int_{S^2} |f|^2$,
\end{itemize}
\hfill (5.8)\\
where $|\fe_1| \leq c_0 { 1\over |\ln \epsilon| } n (1+\lambda)$ and $|\fe_2| \leq c_0\epsilon n(1+\lambda)$ when $f$ is in the span of the respective $\iq$-eigensections or $\ip$-eigensections with eigenvalue less than $\lambda$.\\

\it	Part 2:  \rm Choose a positive number to be denoted by $\mu$  so that there are no $\iq$-eigenvalues in the interval between $E-\mu$ and $E+\mu$.  Supposing that $\epsilon$ is small, assume that there exists a configuration $\fp$ from $\mathcal{U}(\epsilon,k,[\fq]$) with an $\ip$-eigenvalue in the interval between $E- {1\over2}\mu$  and $E+ {1\over2}\mu$.  The plan is to derive nonsense from this assumption.  To this end, let $f$ denote such an eigensection with the $S^2$ integral of $f^2$ equal to 1.  Write $\chi_{\fp,\epsilon}f$ as a sum of $\iq$-eigensections.  This depiction has the form $f_<+f_>$ where $f_<$  denotes the contribution to this sum from the eigensections with eigenvalue less than $E-\mu$ (and thus $f_>$ contains the contribution from the eigensections with eigenvalue greater than $E+\mu$).  Apply (5.8) to $\chi_{\fp,\epsilon}f_<$ (viewed as a section of $\ip$) to see that it is, for the most part, a sum of $\ip$-eigensections with eigenvalue less than $ E- {1\over 2}\mu$  when $\epsilon$ is very small.  Therefore, the $L^2$-inner product of $\chi_{p,\epsilon}f_<$ with $f$ is very small when $\epsilon$ is very small. 
  Thus, $f_<$ itself must be very small when $\epsilon$ is very small because the $L^2$ inner product between $\chi_{\fp,\epsilon}f_<$ and $f$ is the same as that between $f_<$ and $\chi_{\fp,\epsilon}f$ and the latter inner product is the $S^2$ integral of $f_<^2$ (by definition). As explained directly, the detailed analysis leads to this: 
\begin{equation}
\int_{S^2} |f_<|^2 \leq c_0  {1\over\mu|\ln\epsilon|}n (1+E)
\tag{5.9}
\end{equation}
%
%

To derive the preceding inequality, start with the observation that
\begin{equation}
\int_{S^2} |df_<|^2 \leq (\mathrm{E} - \mu) \int_{S^2} |f_<|^2.
\tag{5.10}
\end{equation}
With this in hand, write its left hand side as the $S^2$ integral of $\langle df_<, d(\chi_{\fp,\epsilon},f)\rangle$ (the integral of the inner product between $df_<$ and $d(\chi_{\fp,\epsilon}f)$).  Then move the $\chi_{\fp,\epsilon}$ from multiplying $f$ to multiplying $f_<$ so as to write that integral as  
\begin{equation}
\int_{S^2} \langle d(\chi_{\fp,\epsilon}f_<), df\rangle + \mathfrak{\tau}_1.
\tag{5.11}
\end{equation}
where $\mathfrak{\tau}_1$ is a term whose norm is bounded by $c_0$ times the $S^2$-integral of the function $|d\chi_{\fp,\epsilon}|(|f_<||df|+|f||df_<|)$; and thus by $c_0{1\over |\ln \epsilon|}n$ times the $S^2$ integral of $|df|^2$.  (Use Lemma 2.1 and the fact that $|d\chi_{\fp,\epsilon}|$ is bounded by $c_0 {1\over |\ln \epsilon|} {1\over \mathrm{dist}(\fp, \cdot)}$.  The factor of $n$ appears because $\fp$ has $2n$ points. Let $\lambda$ denote $f$’s eigenvalue and write (5.11) as
\begin{equation}
\lambda \int_{S^2} |f_<|^2+ \mathfrak{\tau}_1.
\tag{5.12}
\end{equation}
Because $\lambda > E - {\mu\over 2}$, the inequality in (5.10) with (5.12) imply that
\begin{equation}
{\mu\over 2} \int_{S^2} |f_<|^2 \leq c_0 {1\over |\ln \epsilon|} n (1 + \mathrm{E})
\tag{5.13}
\end{equation}
which is what is claimed by (5.9).

	The inequality in (5.9) implies (because $\chi_{\fp,\epsilon}f$ is nearly $f_>$) that
\begin{equation}
\int_{S^2} |d (\chi_{\fp,\epsilon}f)|^2 \geq  (\mathrm{E} + \mu) \left(1 - c_0 \mu^{-1} {1\over |\ln \epsilon|} n (1 + \mathrm{E})\right)
\tag{5.14}
\end{equation}
The latter inequality implies in turn (via 5.8) that 
\begin{equation}
\int_{S^2} |d f|^2 \geq c_0 (\mathrm{E} + \mu) \left(1 - c_0 \mu^{-1} {1\over |\ln \epsilon|} n (1 + \mathrm{E})\right)
\tag{5.15}
\end{equation}which exhibits the desired nonsense if $\epsilon$ is sufficiently small (given $\mu$ and E) because the left hand side of this inequality is, by assumption, no greater than $\mathrm{E} +{1\over 2\mu}$.\\

	\it Part 3:  \rm Let $\Pi_{\fp,E}$ denote the projection operator on $\hp$ onto $\hpe$ for the $L^2$-inner product on $\hp$.   It follows from the $\mathbb{H}_{\fq}$ to $\hp$ version of (5.8) that the composition of first $L$ and then $\Pi_{\fp,E}$ is injective on $\mathbb{H}_{\fq,E}$ when $\epsilon$ is small (given $E$ and $\mu$); and nearly isometric with respect to the $\mathbb{H}_{(\cdot)}$ and the $L^2$ inner products as indicated in the statement of the proposition.  
	
Much the same arguments can be used with the $\hp$ to $\mathbb{H}_{\fq}$
 version of $L$ to see that the composition of first this version of $L$ and then $\Pi_{\fq,E}$ maps $\hpe$ injectively into $\mathbb{H}_{\fq,E}$ when $\epsilon$ is small (given $E$ and $\mu$); and that this composition is nearly an isometry when $\epsilon$ is small (given $E$ and $\mu$).\\

\noindent\bf Proof of Lemma 5.4:   \rm The three parts of the proof that follow explain how to construct the desired deformation retract.\\

\it Part 1:  \rm Here is the key observation:  Suppose that $[\fq]$ is a $\integer/2\integer$-divisor and let $\fq$ denote its minimal representative.  Suppose also that $E_{\fq} > 2\mathfrak{C}$ and that $E_{\fq}$ is not an $\iq$-eigenvalue.  Let $f$ denote an element in $\mathbb{H}_{\fq}$ with $\mathfrak{C}$ bounding the $S^2$-integral of $|df|^2$.  Then the $S^2$-integral of $|\Pi_{\fq,E_{\fq}} f|^2$ is greater than ${1\over 2}$ if the $S^2$-integral of $|f|^2$ is equal to 1.  

Given a positive integer $k$, then $E_{\fq}$ can be chosen larger if necessary so that the dimension of $\mathbb{H}_{\fq, E_{\fq}}$ is greater than $k$.  This lower bound is assumed in what follows.\\

\it Part 2: \rm Keeping Part 1 in mind, there is an open neighborhood of $[\fq]$ in $\ctwonb$ to be denoted by $\overline{U}_{\fq}$ such that if $[\fq^\prime]$ is from this neighborhood, then $E_{\fq}$ is not an $\mathcal{I}_{\fq^{\prime\prime}}$-eigenvalue.  Now each $[\fq] \in\ctwonb$ has its version of $\overline{U}_{\fq}$ and since the union of all of these sets covers $\ctwonb$, there is a finite subcover.  Choose such a finite cover and denote this finite collection of sets by $\uk$.  Let $\{\chi_{\fq}: \overline{U}_{\fq} \in \uk\}$ denote a chosen partition of unity subbordinate to the cover.\\

\it Part 3: \rm Supposing that $([\fp],[f])$ is in $\rpb^{\mathfrak{C}}$, define the new $\ip$ section $f_1$ by the rule
\begin{equation}
f_1 \equiv \sum_{U_{\fq} \in\mathfrak{U}_k} \chi_{\fq} \Pi_{\fp, E_{\fq}}f.
\tag{5.16}
\end{equation}
What follows are two important facts about $f_1$.  To state the first, let $E$ denote the minimum of those $E_{\fq}$ that appear with non-zero $\chi_{\fq}$ in (5.16).  Then the $S^2$-integral of  $|\Pi_{\fp,E}f_1|^2$ is greater than ${1\over 2}$  because $\Pi_{\fp,E}f_1$ is the same as $\Pi_{\fp,E}f$.   This implies in particular that $f_1$ is not identically zero.  With this last point understood, introduce $Z_1$ to denote the square root of the $S^2$-integral of $|f_1|^2$.

To state the second fact, let $E^\prime$ denote the maximum of those $E_{\fq}$ that appear with non-zero $\chi_{\fq}$ in (5.16).  Then $\Pi_{\fp,E^\prime} f_1 = f_1$ which is to say that $f_1$ is in $\mathbb{H}_{\fp, E^\prime}$.  This second fact has implications with regards to Lemma 5.3's weak $\rp^\infty$ structure for $\rpb$:  The $\rpb$ element $([\fp], {1\over Z_1} [f_1])$ is in the image of at least one $\overline{U}_{\fq} \in \uk$ version of $\Psi_{\overline{U}_{\fq}}$.  (Remember that  $\Psi_{\overline{U}_{\fq}}$ is a fiber preserving embedding of $(V_{\overline{U}_{\fq}} -0)/\real^*$ to $\rpb$), and remember that the fiber of $V_{\overline{U}_{\fq}}$ over $\fp$ is $\Pi_{\fp, E_{\fq}}$.)

It follows from the preceding fact that the map from $\rpb^\mathfrak{C}$ to $\rpb$ that sends any given element $([\fp],f)$ to $([\fp],f_1)$ maps $\rpb^{\mathfrak{C}}$ into the subspace $\bigcup_{\overline{U}_{\fq} \in \mathfrak{U} } \Psi_{\overline{U}_{\fq}  } (V_{\overline{U}_{\fq}  } - 0)/ \real$. \\

\it	Part 4:  \rm The desired deformation retract is a map from $[0, 1] \times\rpcb$ to $\rpcb$ whose $t=1$ end member is the map from Part 3 that send $([\fp], f)$ to $([\fp],[f_1])$.  The time $t \in [0,1]$ member of the family sends $([\fp],f)$ to $([\fp], [f_1+(1-t)(f-f_1)])$.  This is a well defined map from $\rpcb$ because $f_1$ is non-zero and it is fiber preserving because the $[\fp]$ component of an given element doesn't change.  Note also that any $ t \in[0, 1]$ version maps $\rpcb$ to itself.  Indeed, supposing that $([\fp],[f])$ is from $\rpcb$, then the projection of $f$ onto the span of the $\ip$-eigensections with eigenvalue less than or equal to $\mathfrak{C}$ does not change as $t$ increases whereas the $L^2$ norm of the projection of $f$ to the span of the $\ip$-eigensections with eigenvalues greater than $\mathfrak{C}$ is a non-increasing function of $t$.  And, by construction, the $t=1$ end member of this family maps into $\bigcup_{\overline{U}_{\fq} \in \mathfrak{U} } \Psi_{\overline{U}_{\fq}  } (V_{\overline{U}_{\fq}  } - 0)/ \real$.

\subsection{The obstruction class}

	Returning to the context of (5.5), suppose that $(\mathcal{R}, X, \pi)$ has a weak $\rp^\infty$ bundle structure.  This weak $\rp^\infty$ bundle structure comes from a \it weak vector bundle \rm structure when the collection $\{V_U: U \in \uk\}_{k\in \mathbb{N}}$ from (5.5) has the additional property described below in (5.17).  To set notation, suppose that $k, k^\prime$ are non-negative integers, and then suppose that $U$ is from $\uk$ and $U^{\prime}$ is from $\mathfrak{U}_{k^\prime}$, that these sets with non-empty intersection, and that their corresponding vector bundles obey $\dim(V_{U^{\prime}}) \geq \dim(V_U)$.  Let $T_{U^{\prime},U}$ denote the injective bundle map from $V_U|_{U\cap U^{\prime}}$ to $V_{U^{\prime}}|_{U\cap U^{\prime}}$ that is required by the second bullet of (5.5).   \\

\begin{center}\it
For any triple $k, k^\prime, k^{\prime\prime} \in \mathbb{N}$ and sets $U, U^{\prime}$ and $U^{\prime\prime}$ from $\uk, \ukp$ and $\ukpp$ that share points and are such that $\dim(V_{U^{\prime\prime}}) \geq \dim(V_{U^{\prime}})\geq \dim(V_U)$, the corresponding vector bundle maps obey the identity $T_{U^{\prime\prime}U} = T_{U^{\prime\prime}U^{\prime}}T_{U^{\prime}U}$.
\end{center}
\hfill (5.17)\\

\rm

Now a weak $\rp^\infty$ bundle structure need not come from weak vector bundle structure (which turns out to be our situation with our bundle over $\ctwonb$).  To see the obstruction, note first that whether or not the weak $\rp^\infty$ bundle structure comes from a weak vector bundle structure, the bundle maps for the weak $\rp^\infty$ structure do obey 
\begin{equation}T_{U^{\prime\prime}U }= x_{_{U^{\prime\prime}U^{\prime}U}}T_{U^{\prime\prime}U^{\prime}}T_{U^{\prime}U}
\tag{5.18}\end{equation}
with $x_{_{U^{\prime\prime}U^{\prime}U}}$ being a nowhere zero, continuous map from $U^{\prime\prime}\cap U^\prime \cap U$ to $\mathbb{R}$. 
Let $z_{_{U^{\prime\prime}U^{\prime}U}}$ denote the sign (either 1 or $-1$) of $x_{_{U^{\prime\prime}U^{\prime}U}}$.  This is a locally constant map from $U^{\prime\prime}\cap U^{\prime}\cap U$ to $\integer/2\integer$.  These locally constant maps are such that the collection  
\begin{equation}
\left\{z_{_{U^{\prime\prime}U^{\prime}U}}: U^{\prime\prime}, U^{\prime}, U \in \bigcup_{k\in \mathbb{N}} \uk\right\} 
\tag{5.19}
\end{equation}
when restricted to any locally finite subcover defines a (multiplicative) degree two  \v{C}ech cocycle for the subcover whose image in $H^2(X; \integer/2\integer$) is independent of the chosen subcover.  This cohomology class is said to be the \it obstruction class.  \rm As explained in the next paragraphs, if the first $\integer/2\integer$ cohomology of all sets from $\knuk$ is zero, then this weak $\rp^\infty$ structure comes from a weak vector bundle structure if and only if this class is zero.  (In general, if the class is zero, then the collection of $T_{U^{\prime},U}$ maps can be modified 
by multiplying each by a suitable nowhere zero function
so that they obey the identity in (5.17).)  

To explain why the data in (5.19) always defines a cochain (a closed cocycle):  Consider four sets $U_1$, $U_2$, $U_3$ and $U_4$ that share points and with corresponding vector bundles whose dimensions don't decrease with increasing label.  Write $T_{ji}$ for the $U=U_i$ and $U^{\prime}=U_j$ version of $T_{U^{\prime}U}$ and consider $T_{43}T_{32}T_{21}$.  This can be written on the one hand as $x_{_{321}}x_{_{431}} T_{41}$ and on the other as $x_{_{432}}x_{_{421}} T_{41}$.  Thus, $z_{_{321}}z_{_{431}}z_{_{432}}z_{_{421}} = 1$ which is the statement that the collection in (5.19) defines a closed cocycle.  

With regards to changing (5.19) by a (multiplicative) coboundary:  Changing (5.19) by the differential of a 1-cocycle $\left\{z_{_{U^{\prime}U}}: U, U^{\prime} \in \knuk\right\}$, accounts for the change of each $T_{U^{\prime}U}$ to $z_{_{U^{\prime}U}} T_{U^{\prime}U}$.   If this change can be made so as to make all $x_{_{U^{\prime\prime}U^{\prime}U}}$ positive, then a partition of unity can be employed to modify each $T_{U^{\prime}U}$  so that the condition in (5.17) holds.  (To do this, write $x_{_{U^{\prime\prime}U^{\prime}U}}$ as the exponential of a real valued function.)

If the first cohomology of the sets from the collection $\knuk$ are not all trivial, then there is an additional equivalence that can be used to obtain the given weak $\rp^\infty$ structure from a weak vector bundle structure:  Each $U \in \knuk$ version of $V_U$ can be modified by tensoring it with a real line bundle over $U$ with the proviso that any two respective $U$ and $U^{\prime}$ line bundles are isomorphic on $U\cap U^{\prime}$.  Assuming this, let  $\tilde V_U$ denote the modified $V_U$.  The modified version of $T_{U^{\prime}U}$ is obtained from the original by tensoring with a chosen isomorphism from the respective $U$ line bundle to the respective $U^{\prime}$ line bundle.  Let $\eta_{_{U^{\prime}U}}$ denote a chosen isomorphism between these line bundles.   
When $U, U^{\prime}$ and $U^{\prime\prime}$ are three sets that share points, then $\eta_{_{ U^{\prime\prime}U}}$ can be written as  $\hat x_{_{U^{\prime\prime}U^{\prime}U}}\eta_{_{ U^{\prime\prime}U^{\prime}}}\eta_{_{U^{\prime}U}}$ with  $\hat x_{_{U^{\prime\prime}U^{\prime}U}}$ mapping to $\real^*$.  Let  $\hat z_{_{U^{\prime\prime}U^{\prime}U}}$ denote the sign of  $\hat x_{_{U^{\prime\prime}U^{\prime}U}}$.  The collection of  $\hat z$'s also defines a multiplicative, degree 2  \v{C}ech cocycle.  Modifying each $U \in \knuk$ versions of $V_U$ in this way changes each $T_{U^{\prime}U}$ to $\eta_{_{ U^{\prime}U}} T_{U^{\prime}U}$ and this corresponds to changing each $z_{_{U^{\prime\prime}U^{\prime}U}}$ to the product  $\hat z_{_{U^{\prime\prime}U^{\prime}U }}z_{_{U^{\prime\prime}U^{\prime}U}}$.

	The relevant example is $\rpb$ where the obstruction class in $H^2( \ctwonb)$ is non-zero (its restriction to $\rp$ is already non-zero.)  As was the case with $\rp$, the $\rpb$ obstruction class is denoted by $\omega$.  
	
The next lemma says more about this class (it is needed for the next subsection).\\

\noindent \bf Lemma 5.5: \it  The self-cup product $\omega^{2n}$ is non-zero in $H^{4n}(\ctwonb)$.\\ \rm

\noindent \bf  Proof of Lemma 5.5:  \rm A preliminary step is a proof that $H^{4n}(\ctwonb)$ is non-zero.  To prove that this is so, consider that the exact sequence for the pair $(\ctwonb, \ctwonmtb)$ appears in part as
\begin{equation}
H^{4n-1}(\ctwonmtb) \to H^{4n}(\ctwonb,\ctwonmtb) \to H^{4n}(\ctwonb) \to H^{4n}(\ctwonmtb)
\tag{5.20}
\end{equation}
whose left-most and right-most terms are zero.  Thus, $H^{4n}(\ctwonb)$ is isomorphic to $H^{4n}(\ctwonb,\ctwonmtb)$.  Meanwhile, $H^{4n}(\ctwonb,\ctwonmtb)$ is isomorphic to the compactly supported cohomology of $\ctwon$ which is isomorphic to $\integer/2\integer$ since $\ctwon$ is a manifold. 

	The calculation for the $2n$-fold self-cup product of $\omega$ has three parts.\\

\it Part 1:  \rm The standard compactification of $\ctwon$ is to view it as a subset of $\mathbb{CP}^{2n}$ as follows:  First view $S^2$ as $\cp^1$ which is to say $(\complex^2-0)/\complex^*$.  Letting $(z, w)$ denote the coordinates for $\complex^2$ and thus homogeneous coordinates for $\cp^1$,
 a point in $\cp^1$ is determined by homogeneous coordinates $[\alpha, \beta]$ as the zeros of the linear function $\alpha w - \beta z$.  Then, $2n$ points in $\cp^1$ (counting multiplicities) are the zeros of the $2n$'th order homogeneous polynomial
\begin{equation}(\alpha_1w - \beta_1z)(\alpha_2w - \beta_2z) \cdots (\alpha_{2n}w- \beta_{2n}z)
\tag{5.21}
\end{equation}
Multiplying this out, then the coefficients of the various powers of $w$ are the homogeous coordinates for $\cp^{2n}$:
\begin{equation}
a_{2n}w^{2n}+ a_{2n-1}w^{2n-1}z + \cdots + a_0z^{2n} .
\tag{5.22}
\end{equation}
Doing this identifies $\ctwon$ as the open subspace in $\cp^{2n}$ where the polynomial in (5.22) has distinct zeros (the zeros are the homogeneous coordinates $([\alpha_1,\beta_1], \cdots [\alpha_{2n},\beta_{2n}]$).  Viewed in this light, there is a continuous map (to be denoted by $\mathbb{F}$) from $\cp^{2n}$ to $\ctwonb$ that is the identity map on the top strata $\ctwon$ and collapses various lower dimensional varieties in $\cp^{2n}$ corresponding to points where (5.22) has roots with multiplicity onto corresponding strata in the $\ctwonb$ compactification depicted in (3.2).   \\

\it	Part 2:  \rm Let $w$ denote the generator of $H^2(\cp^{2n})$ with it understood that the coefficient ring is $\integer/2\integer$.  The powers of $w$ generate $H^*(\cp^{2n})$, and in particular, $w^{2n}$ is non-zero.  If $\mathbb{F}^*\omega = w$, then $\omega^{2n}$ can't be zero (because pull-back and cup-product commute).  Now the class $\mathbb{F}^*\omega$ is equal to $w$ if and only if its restriction to a complex line in $\cp^{2n}$ is non-zero.   To see about that, consider first the degree $2n-1$ homogeneous polynomial $w^{2n-1}-z^{2n-1}$ which has $2n-1$ distinct roots which are the points where $z = \eta w$ with $\eta$  being a ($2n-1$)'st root of unity.  Supposing that $a \in \complex$, then the polynomial
\begin{equation}
(aw - z)(w^{2n-1}- z^{2n-1}) = a w^{2n} - w^{2n-1}z + awz^{2n-1}-z^{2n} 
\tag{5.23}
\end{equation}
has distinct roots as long as $a^{2n -1} \neq 1$.  In any event, as $a$ varies in $\complex\cup\infty$, the corresponding locus in $\cp^{2n}$ traces out a degree 1 rational curve.  Let $\Sigma$ denote this curve.\\

\it	Part 3:  \rm The $\mathbb{F}$-image of this curve in $\ctwon$ is a set of configurations where $2n$-1 points are fixed to be the points in $\cp^1$ where $z = \eta w$ with $\eta^{2n-1} = 1$, and with the last point (parametrized in homogeneous coordinates by $[a,1]$) varying at will except that it must avoid the $2n-1$ fixed points.  The $\mathbb{F}$-image of $\Sigma$ in $\ctwonb$ adds $2n-1$ configurations in the stratum $\ctwonmt$ from (3.2) which correspond to the collisions between the $a$-parametrized point and each of the $2n-1$ fixed points.  

Now note that if $x$ denotes the $z = 0$ point in $\cp^1$, then $\mathbb{F}(\Sigma) \cap \ctwonmt$ is disjoint from $\mathcal{C}^x_{2n}$.  With that understood, recall that the class $\omega$ is dual in $\ctwon$ to $\mathcal{C}^x_{2n}$.  Since the image $\mathbb{F}(\Sigma)$ intersects $\mathcal{C}^x_{2n}$ just once (where the parameter $a = 0$), the image curve $\mathbb{F}(\Sigma)$ has pairing 1 with $\omega$ in $\ctwon$.  

The preceding calculation is not the end of the story because the closure of the divisor $\mathcal{C}^x_{2n}$ in $\ctwonb$ consists of the whole of $\ctwonmtb$ and in particular, the stratum $\ctwonmt$.  Since $\mathbb{F}(\Sigma)$ intersects $\ctwonmt$ an odd number of times, it is necessary to show that each such intersection counts as zero when evaluating the pairing of $\mathbb{F}(\Sigma)$ with $\omega$.   

This is done by deforming the stratum $\ctwonmtb$ so that a neighborhood of the points in $\mathbb{F}(\Sigma)\cap \ctwonmt$ are pushed into $\ctwon$.  Here is how to do that:  Let $\beta$ denote a chosen function with compact support on $\ctwonmt-\mathcal{C}_{2n-2}^x$ that is equal to 1 at the $2n-1$ configurations that comprise the intersection of $\mathbb{F}(\Sigma)$ with $\ctwonmt$, and equal to zero on any configuration that contains a point with distance less than ${1\over100n}$  from $x$.  Define $\mathcal{F}: [0,1] \times\ctwonmtb \to \ctwonb$ by the rule whereby 
\begin{equation}
\mathcal{F}(t,\fq) = \fq \cup \left(\left[{1\over 1000n} t\beta( \fq),1\right], \left[-{1\over 1000n} t\beta(\fq), 1\right]\right).
\tag{5.24}
\end{equation}
Of particular note is that $\mathbb{F}(\Sigma)$ is disjoint from $\mathcal{F}(t=1,\fq)$.  Thus, the $\mathbb{F}$-image of $\Sigma$ in $\ctwonb$ also has pairing 1 with the class $\omega$.  This is to say that $\mathbb{F}^*\omega$ has pairing 1 with $\Sigma$.  As noted previously, this implies that $\mathbb{F}^*\omega = w$.

\subsection{Some non-zero classes in the cohomology the space $\rpb$}

	A part of the cohomology of $\rpb$ is described in the next lemma.  The lemma uses $\omega$ to denote the obstruction class in $H^2(\ctwonb)$ as was done in the previous subsection. \\

\noindent \bf Lemma 5.6:  \it The cohomology of $\rpb$ contains the vector space of polynomials in a degree two class $\tau$ with coefficient ring $\mathcal{Q}$ as described in the subsequent bullets.  Thus,
$$\mathcal{Q} \oplus \mathcal{Q}\tau \oplus \mathcal{Q} \tau^2 \oplus \mathcal{Q}\tau^3 \oplus \cdots$$
\begin{itemize}
\item The definition of $\tau$:  Let  $\hat \omega$ denote a given representative cocycle for $\omega$.  Then $\pi^*\hat\omega$  can be written as $d\nu$ with $\nu$ being a 1-cochain whose restriction to each fiber algebraically generates the cohomology of the fiber.  The cochain $\nu\wedge\nu$ is closed and non-zero in $H^*(\rpb)$.  The class of  $\nu \wedge \nu$  is the class $\tau$.
\item The definition of $\mathcal{Q}$:  What is denoted by $\mathcal{Q}$ is the 2-dimensional vector space spanned by the degree zero generator 1 and a degree $4n+1$ class to be denoted by  $\Xi$.  To define the latter, let $\hat\iota$ denote a degree $(4n+1)$-cocycle on $\ctwonb$ obeying $d\hat{\iota} = \hat\omega^{ 2n+1}$.  The degree $(4n+1)$-cocycle $\nu\wedge\pi^*\hat\omega^{ 2n} + \pi^*\hat{\iota}$ on $\rpb$ is closed and non-zero in $H^{2n+1}(\rpb)$.  Its class is independent of the choice for the cocycle $\hat{\iota}$ and it is independent of the choice for the cochain $\nu$ and the choice for the representative cocycle  $\hat\omega$ for $\omega$.  This cohomology class is  $\Xi$.\\
\end{itemize}

\noindent \bf
Proof of Lemma 5.6: \rm Except for the assertion to the effect that  $\Xi$  is independent of the various choices, the assertion of the lemma is an instance of Proposition A.1 in the appendix.  To see about invariance:  Changing any of $\hat\omega$,  $\nu$ and $\hat{\iota}$  maintaining the relations between them changes the cocycle representing  $\Xi$  by something that is exact (by inspection) plus the $\pi$-pull back of a closed cocycle on $X$ of degree $4n+1$.  Meanwhile, all of those are exact because the top dimensional strata of $X$ has dimension $4n$.

Note that if the cohomology of $\rpb$ is computed using a CW decomposition for $\ctwonb$, then $\hat\omega^{2n+1}$ is zero and then $\hat{\iota} \equiv 0$ too.  If simplicial cohomology is used, then there is no guarantee that  $\hat\omega^{2n+1}$ is the zero cochain; but it is cohomologous to zero.

\section{Min-max for $\mathcal{E}$ on $\rpb$}
	There are some hints that the classes from the set $\{\Xi \wedge \tau^m\}_{m\in\{0,1,\ldots\}}$ have min-max values that correspond to critical points of $\mathcal{E}$ on large $n$ version of $\ctwon$.  There are two indications, the first being this:  If it were the case that the $(2m+1)$'st $\ip$-eigenvalue for $\mathcal{E}$ had multiplicity 1 at each configuration (which can't actually happen), then its supremum on $\ctwon$ would be the min-max value for the cohomology class $ \Xi\wedge \tau^m$ .  Thus, min-max using the classes $\{ \Xi\wedge \tau^m\}$ constitute a cohomological work-around for the fact that the functions on $\ctwon$ given by ordering of the $\mathcal{I}_{(\cdot)}$-eigenvalues are not everywhere differentiable.   The second reason for hope concerns the appearance of the $\hat\omega^{2n}$ factor in 
	the definition of
	 $\Xi$ :  This factor sees the whole of $\ctwon$, not just some low codimension subset that can be pushed into $\ctwonb-\ctwon$.
	
	\subsection{The modified definition of min-max value}
	What follow directly is the definition of `min-max value' for the context at hand.  To set the stage, let $\yb$ denote a fiber bundle over $\ctwonb$ with compact, finite dimensional fiber (in practice, the fiber will be $\rp^{2m+1}$) which is smooth along each stratum of $\ctwonb$. A continuous map from $\overline{\mathcal{Y}}$ to $\rpb$ is said to be \it strata preserving \rm if, for each $k \in\{0, \ldots, n\}$, the map sends $\yb|_{\mathcal{C}_{2k}}$  to $\rpb|_{\mathcal{C}_{2k}}$.   Let $\mathfrak{F}$ denote a chosen homotopy class of strata preserving maps from $\yb$ to $\rpb$.  The associated min-max value of $\mathcal{E}$ for the set $\fcc$ is this:

\begin{equation}
\inf_{\Psi \in \fcc} \max_{y\in \yb} \mathcal{E}(\Psi (y)).
\tag{6.1}
\end{equation}

The upcoming Lemma 6.1 is used to obtain the relevant version of $\yb$ and set $\fcc$.  With regards to terminology:  When the lemma speaks of an $\rp^{2m+1}$ fiber bundle, this refers to a special sort of bundle:  Each point in the base has a neighborhood that identifies the fiber bundle over the neighborhood with the $\real^*$ quotient of the complement of the zero section in a $(2m+2)$-dimensional vector bundle.  In addition, the transition functions between two such intersecting open sets are induced by fiberwise linear maps that identify the corresponding vector bundles over the intersection of the two base sets.\\  

\bf \noindent Lemma 6.1:  \it Fix a positive integer for $n$ and a non-negative integer for $m$.  There is an $\rp^{2m+1}$-fiber bundle $\pi_{\mathcal{Y}}: \yb_{n,m} \to \ctwonb$ with a fiber preserving map to $\rpb$ that pushes forward the fundamental class of $\yb_{n,m}$ to a class that has non-zero pairing with $\nu^{m} \wedge \Xi$. (The space $\yb_{n,m}$ is the $X = \ctwonb$ instance of what is denoted by $\mathcal{Y}_m$ in the appendix.)\\

\noindent \bf Proof of Lemma 6.1:  \rm This is an instance of Lemma A.6 in the appendix (Section A.6). \\

Granted Lemma, 6.1, take $\fcc$ to be the homotopy class of maps $\Psi$  as described by Lemma 6.1.  The corresponding version of the min-max value in (6.1) is denoted by $\mathcal{E}_{n,m}$.  

\subsection{Min-max values $\{\mathcal{E}_{n,m}\}_{n>0,m\geq 0}$}

	Fix a positive integer for $n$.  The lemma that follows concerns the size of the min-max values from the set $\{\mathcal{E}_{n,m}\}_{m=0,1,\ldots}$.\\

\noindent \bf Lemma 6.2:  \it Fix a positive integer for $n$. 
\begin{itemize}
\item For any non-negative integer $m$:  The min-max value $\mathcal{E}_{n,m}$ is not less than the maximum over $\ctwonb$ of the function $\mathrm{E}_{(\cdot)}$ from (1.2).  Therefore, it is greater than $\kappa^{-1}n$ with $\kappa$ positive and independent of $n$.  
\item The set $\{\mathcal{E}_{n,m}\}_{m\in\{0, 1, \ldots\}}$ is unbounded.\\
\end{itemize}

\noindent\bf Proof of Lemma 6.2: \rm  With regards to the lemma's first bullet:  Supposing that the min-max value $\mathcal{E}_{n,m}$ is as least as large as the supremum of $\mathrm{E}_{(\cdot)}$ on $\ctwon$, then it is a priori greater than $c_0^{-1}n$ because the supremum of $\mathrm{E}_{(\cdot)}$ on $\ctwon$ is greater than $c_0^{-1}n$.  (See Proposition 3.3.) 

To prove that the min-max value is not less than $\mathrm{E}_{(\cdot)}$, suppose that this isn't the case so as to derive some nonsense.  If $\Psi$  is from $\fcc$, then the push-forward by $\Psi$  of the fundamental class of $\yb$ is a closed, $2m+4n+1$ dimensional chain in $\rpb$ that has non-zero pairing with $\nu^{m} \wedge \Xi$ .  Suppose that the maximum of the values of $\mathcal{E}$ on the image of $\Psi$  is less than the maximum of value of $\mathrm{E}_{(\cdot)}$ on $\ctwon$.   If this is the case, then the $\pi$ image of $\Psi (\yb)$ misses a configuration in $\ctwon$.   (Remember that the maximum of $\mathrm{E}_{(\cdot)}$ on $\ctwonb$ is at configurations in $\ctwon$.)  Now comes an important point:  The cocycle  $\hat\omega^{2n}$ is cohomologous to zero in the complement of any configuration in $\ctwonb$ because the fundamental class of the complement of any configuration is zero.  This is to say that  $\hat\omega^{2n}$ on the complement of a configuration can be written as $d\beta$ with $\beta$ being a degree degree $4n-1$ cochain. Also,  $\hat\omega^{2n+1}$ can be written there as $d(\hat \omega\wedge \beta)$.  This implies that the pairing between 
a cocycle representative of
$\nu^{m}\wedge  \Xi$  and the $\Psi$ -push-forward of the fundamental class of $\yb$  is the same as the pairing between $\nu^{2m+1}\wedge d\pi^*\beta - \nu^{2m}\wedge d\pi^*(\hat\omega \wedge \beta)$ and that class, which is zero because the class is closed and
\begin{equation}
\nu^{2m+1}\wedge d\pi^*\beta - \nu^{2m}\wedge d\pi^*(\hat \omega \wedge \beta) = -d(\nu^{2m+1}\wedge \beta).
\tag{6.2}
\end{equation}
This conclusion is the desired nonsense because it runs afoul of the assumptions about $\Psi$.

	With regards to the second bullet:  Let $\mathfrak{c}_{2m}$ denote a homology class on $\rpb$ that is dual to $\nu^{2m}$.  This class is represented by the fundamental class of any embedding of $\rp^{2m}$ in any fiber of $\rpb$ that induces an isomorphism on the first homology of the fiber.  An important point is that the min-max value of $\mathcal{E}$ for $\mathfrak{c}_{2m}$ is an unbounded function of $m$.  This follows from instances of Proposition 1.1 and Lemmas 5.2 and 5.3 which together imply the following:  Given any positive number $E$, there exists an integer (call it $k_E$) such that the $k_E$'th lowest $\ip$-eigenvalue is greater than $E$ at each configuration in $\ctwonb$.  This implies that the min-max value for $\mathfrak{c}_{2m}$ is greater than $E$ if $2m > k_E$.  As explained next, the number $E$ (assuming $2m > k_E$) is a lower bound for $\mathcal{E}_{n,m}$.
	
To prove that $\mathcal{E}_{n,m} \geq E$, suppose to the contrary that $\mathcal{E}_{n,m}$ is less than $E$ to generate nonsense.  In this event, there is a map $\Psi: \yb\to \rpb$ from $\fcc$ with $\Psi^*\nu^{2m}$ evaluating to zero on the generator of $2m$'th homology of the fibers of the projection from $\yb$ to $\overline{\mathcal{C}}_{2m}$.  Meanwhile, the cohomology of $\yb$ is described by Lemma A.2 in the appendix, and it follows from that description, and from what is said in Sections A.2 and A.3, that Lemma A.2's cohomology class $\tau$ can be taken to be that of $\Psi^*\nu^2$.  And, according to Lemma A.2, the class $\tau^m$ does have non-zero pairing with the $2m$'th homology of the fibers of the projection to $\ctwonb$.  

\subsection{Convergence to a critical point}
	The proposition that follows directly makes an assertion to the effect that the min-max value of $\mathcal{E}$ for any $m \geq 0$ version of $\nu^{m}\wedge \Xi$  is a critical value of $\mathcal{E}$ on some strata of $\ctwonb$.  This is to say that there exists $k \in\{0, 1, \ldots, n\}$ and a critical point of $\mathcal{E}$ on $\mathcal{C}_{2k}$ where $\mathcal{E}$'s value is equal to the min-max value $\mathcal{E}_{n,m}$.\\

\noindent \bf Proposition 6.3:  \it For each positive integer $n$ and non-negative integer m, there exists $k\in\{0,\ldots, n\}$ and a critical point of $\mathcal{E}$ on $\rp|_{\mathcal{C}_{2k}}$ where $\mathcal{E}$ is equal to $\mathcal{E}_{n,m}$. \\

\noindent \bf Proof of Proposition 6.3:  \rm The proof of this proposition has eight parts.  Part 1 has a preliminary observation. Parts 2-4 focus on just the differential of $\mathcal{E}$ along the fibers of the projection to $\ctwonb$.  The final parts consider the behavior of the differential of $\mathcal{E}$ on lifts to $\rpb$ of vector fields along the strata of $\ctwonb$. \\

\it	Part 1: \rm  Here is an important point to keep in mind:  Fix a positive number to be denoted by $E$ and let $\rpb^E$ denote the subspace of elements in $\rpb$ that are characterized as follows:  An element $(\fp,[f])$ is in $\rpb^E$ if and only if $f$ is in the span of the $\ip$-eigenvectors with eigenvalue at most $E$.  The point to keep in mind is that the subspace $\rpb^E$ is compact.\\

\it Part 2: \rm Suppose for the moment that $(\fp,[f])$ is a given element in $\rpb$.  Let $E$ denote a positive number that is not an $\ip$-eigenvalue.  Write $f$ as $f_<$ + $f_>$ where $f_<$ is the $L^2$ orthogonal projection of $f$ onto the subspace $\mathbb{H}_{\fp,E}$ (the subspace in $\hp$ that is spanned by the $\ip$-eigensections with eigenvalue less than $E$).  An observation follows:  If $\mu   > 0$ and if the $S^2$-integral of |$f_>|^2$ is greater than $\mu$, then $\mathcal{E}(\fp,[f])$ will be greater than $\mu E$.   

The preceding observation has the following consequence:   Suppose that $\delta > 0$ and that $\Psi$  is a map from $\fcc$ that maps the whole of $\yb$ to the $\mathcal{E}(\cdot)< \mathcal{E}_{n,m}+\delta$ part of $\rpb$.  If $(\fp,[f])$ is from the image of $\Psi$ , and if $E$ is greater than $\mathcal{E}_{n,m}+\delta$, then the $S^2$ integral of the corresponding |$f_>|^2$ is strictly less than 1 $\left(\text{it is at most }{\mathcal{E}_{n,m}+ \delta \over E}\right)$.  This implies, in particular, that $f_<$ is non-zero.

	As explained momentarily, the observation in the previous paragraph can be used to continuously deform the map $\Psi$  to a new strata preserving map that has the properties listed in the next lemma.\\

\noindent \bf Lemma 6.4: \it  Fix $\delta > 0$ and suppose that $\Psi$  is a map from $\fcc$ when $\mathcal{E}(\Psi (\cdot)) < \mathcal{E}_{n,m}+\delta$.  Then there exists a homotopy of $\Psi$  in $\fcc$ whose end member (denoted by $\Psi^\prime$) obeys the following:
\begin{itemize}
\item The map $\Psi^\prime$ is from $\fcc$ .
\item Supposing that $y \in\yb$, then $\mathcal{E}(\Psi^\prime(y)) \leq \mathcal{E}(\Psi (y))$
\item Supposing that $y \in\yb$ write $\Psi(y)$ as $(\fp,[f])$.  Then $f$ is in the span of the $\ip$-eigensections with eigenvalue less than $\mathcal{E}_{n,m}+2\delta.$ \\
\end{itemize}

\noindent \bf Proof of Lemma 6.4:  \rm To set up the construction of $\Psi^\prime$, suppose for the moment that $\fp$ is a given element in $\ctwonb$.  Fix a number to be denoted by $\ep$ that is greater than $\mathcal{E}_{n,m} +\delta$ but less than $\mathcal{E}_{n,m}+2\delta$ and which is not an $\ip$-eigenvalue.  According to Lemma 5.2, there is an open neighborhood of $\fp$ in $\ctwonb$ with the property that if $\fq$ is in this neighborhood, then $\ep$ is not an $\iq$-eigenvalue.  Fix such a neighborhood and denote it by $\mathrm{U}_{\fp}$.  The various $\fp\in\ctwonb$ versions of $\ub _{\fp}$ form an open cover of $\ctwonb$ and thus there is a finite subcover.  This is to say that there is a finite set of points whose corresponding $\ub _{(\cdot)}$ sets cover $\ctwonb$.  Choose such a finite set and denote it by $\Lambda$. Then choose a partition of unity subbordinate to the associated cover.  When $\fp$ is from $\Lambda$, then the corresponding function from the partition of unity is denoted below by $\chi_{\fp}$.

Supposing now that $y \in \yb$, write $\Psi(y)$ as $(\fq,[f])$.  Fix for the moment an element $\fp\in\Lambda$ whose corresponding set $\ub _{\fp}$ contains $\fq$.  Then write $f$ as $f_{\fp,<} + f_{\fp,>}$  with the former being the projection of $f$ to the span of the $\iq$-eigenvalues with eigenvalue less than $\ep $. For any given $t \in[0,1]$, define $f_t$  by the rule whereby
\begin{equation}
f_t \equiv {1\over Z_t} \left( 
\sum_{\fp \in\Lambda} \chi_{\fp} f_{\fp, <} + (1-t) \sum_{\fp \in \Lambda} \chi_{\fp} f_{\fp, >} 
\right)
\tag{6.3}
\end{equation}
where $Z_t$ is positive and chosen so that the $S^2$-integral of $|f_t|^2$ is equal to 1.  Because $f_{\fp,<}$ is non-zero, this deformation is continuous with respect to $t$.
  
Lemma 5.2 guarantees that the assignment of $t \in[0,1]$ and $y \in \yb$  to $\Psi_t(y)$ defines a continuous map from $[0,1]\times \yb$ into $\rpb$.  Since it is strata preserving, this is a homotopy of $\Psi$  in $\fcc$.  Set $\Psi^\prime$ to be the $t = 1$ member.  The construction is such that if $y\in \yb$, then the map $t\to\mathcal{E}(\Psi_t(y))$ (which is $\mathcal{E}(\fq,f_t))$ is non-increasing and it is strictly decreasing if $f_{\fp,>}$ is non-zero for at least one $\fp\in\Lambda$.  By construction, $\mathcal{E}(\Psi_1(y))\leq \mathcal{E}_{n,m}+2\delta$ because $f_1$ has zero projection to the span of the $\iq$-eigenvectors with eigenvalue greater than $\mathcal{E}_{n,m}+2\delta$.\\

\it Part 3:  \rm The preceding lemma has two important implication with the first being this:  Having fixed $\delta$ for Lemma 6.4, suppose that $(\fp,[f])$ is in the image of Lemma 6.4's map $\Psi^\prime$.  If there are \it no  \rm $\ip$-eigenvalues between $\mathcal{E}_{n,m}-\delta^\prime$ and $\mathcal{E}_{n,m}+2\delta$ with $\delta^\prime$ any given positive number, then $\mathcal{E}(\fp,[f])$ is at most $\mathcal{E}_{n,m}-\delta^\prime$.
  
Here is the second implication: \it \\

\begin{center}
\it
With $\Psi^\prime$ as described by Lemma 6.4, there exists $(\fp, [f])$ from the image of $\Psi^\prime$ with $\mathcal{E}(\fp,[f]) \geq \mathcal{E}_{n,m}$ such that the projection of $f$ to the span of the $\ip$-eigensections with eigenvalues between
$\mathcal{E}_{n,m} - \sqrt{\delta}$ and $\mathcal{E}_{n,m} + 2\delta$ has $L^2$ norm greater than $(1-3\sqrt{\delta})^{1\over 2}.$
\end{center}
\hfill (6.4)

\noindent \rm To prove this:  Supposing that (6.4) fails for a given element $(\fp, [f])$ from the image of $\Psi^\prime$, write $f = f_< + f_>$ with $f_<$ denoting the $L^2$ orthogonal projection of $f$ to the span of the $\ip$-eigensections with eigenvalue less than $\mathcal{E}_{n,m} - \sqrt{\delta}$.  Using this decomposition, one has that
\begin{equation}
\mathcal{E} ( \fp, [f]) \leq (\mathcal{E}_{n,m} - \sqrt{\delta}) \|f_<\|^2_{\mathbb{L}} + (\mathcal{E}_{n,m} + 2\delta) (1 - \|f_<\|^2_{\mathbb{L}});
\tag{6.5}
\end{equation}
and this is strictly less than $\mathcal{E}_{n,m} - \delta$ if $\|f_<\|^2_{\mathbb{L}}$ is greater than $3\sqrt{\delta}$ which is when $\|f_>\|^2$ is less than $1-3\sqrt{\delta}$.  Thus, (6.4) can’t fail for every $f$ with $\mathcal{E}(\fp,[f]) \geq \mathcal{E}_{n,m}$ because the maximum value of $\mathcal{E}$ on the image of $\Psi^\prime$ would then be less than $\mathcal{E}_{n,m}$ which violates the assumptions.\\

\it Part 4:  \rm What is said by Part 1 about compactness, and what is said in Lemma 5.6 and (6.4) have the following implications:  Fix any decreasing sequence of positive numbers with limit zero (denoted by $\{\delta_j\}_{j=1,2,\ldots}$) and then a corresponding sequence of maps from $\fcc$ (denoted by $\{\Psi_j\}_{j=1,2,\ldots}$ ) such that for each index $j$, the conditions in Lemma 5.4 are obeyed with $\delta = \delta_j$ and with $\Psi^\prime = \Psi_j$.  For each $j$, set $\mathcal{A}_j$ to be the $\mathcal{E} \geq \mathcal{E}_{n,m}$ part of $\Psi_j(\yb)$.  The set of limit points of the sequence $\{A_j\}_{j=1,2,\ldots}$ contains at least one $(\fp,[f])\in\rpb$ with $f$ being an $\ip$-eigenvector with eigenvalue is $\mathcal{E}_{n,m}$, which is to say that $-\Delta f= \mathcal{E}_{n,m} f$.   (Remember that the sequence $\{A_j\}_{j=1,2,\ldots}$ is a sequence in a compact subset of $\rpb$.  As a consequence, any corresponding sequence of elements (the first from $\mathcal{A}_1$, the second from $\mathcal{A}_2$, and so on) has a convergent subsequence.  By way of terminology:  The set of limits of these subsequences is the set of limit points of $\{\mathcal{A}_j\}_{j=1,2,\ldots}.$) \\

\it	Part 5: \rm This part of the proof is a digression to present some background and then introduce notation for the remaining parts.  To start the digression, suppose for the moment that $\{\phi_t\}_{t\in\real}$ is a 1-parameter family of area preserving diffeomorphisms of $S^2$ with the $t = 0$ member being the identity diffeomorphism.  This family acts in the obvious way on $\ctwonb$ so as to preserve the stratification in (3.2).  To elaborate:  Supposing that $k \in\{0, \ldots, n\}$ and $\fp$ is in the stratum $\mathcal{C}_{2k}$, then $\phi_t(\fp)$ is the configuration whose points are the $\phi_t$-images of the points in $\fp$.  This 1-parameter family can be lifted to act on $\rpb$ as follows:  The parameter $t$ member $\phi_t$ sends any given element $(\fp,[f])$ to the element $(\phi_t(\fp),[f \circ \phi_t^{-1}])$.  This lifted map is denoted by $\bar\phi_t$.

What follows are some remarks about this definition.  The first remark is that an identification between respective line bundles $(\phi_t^{-1})^*\ip$ and $\mathcal{I}_{\phi_t(\fp)}$  is not required to define $\bar\phi_t$ because the two possible isometric identifications differ by a sign.  The second remark is that the lifted map $\bar\phi_t$ does not map $\ip$-eigensections to  $\mathcal{I}_{\phi_t(\fp)}$-eigensections.  The third remark is that the lifted map is not differentiable with respect to $t$; it is only continuous.  This is because the $t$-derivative of the $\phi_t^{-1}$-pull back of a generic element in $\mathbb{H}_{\phi_t(\fp)}$  will not have finite $\mathbb{H}_{\fp}$-norm.  (The $\hp$ norm of the $t$-derivative of the pull-back will be finite however if the integral on the right hand side of (4.14) is finite for the element in question.  See Section 4.4 for more about this issue.)  

	Now suppose that $\Psi : \yb \to \rpb$ is a map from the set $\fcc$ and that $\tau: \yb \to [0,1]$ is a continuous map.  The preceding constructions can be used to obtain a 1-parameter family of deformations of $\Psi$  with any given $t\in\real$ version being the map (denoted by $\Psi_{\tau,t}$) that is defined by the rule
\begin{equation}
y \to \Psi_{\tau,t} (y) \equiv \bar\phi_{\tau (y)t} (\Psi (y)).
\tag{6.6}
\end{equation}

Each member of this family is from $\fcc$.   (Each member is strata preserving because both $\bar\phi_{(\cdot)}$ and $\Psi$  are strata preserving. \\

\it	Part 6:  \rm This part of the proof starts with a digression to serve as a reminder with regards to observations from Section 3.  In particular, remember that when $(\fp,[f])$ is a given element in $\rpb$, then the value of $\mathcal{E}$ on $\bar\phi_t(\fp,[f])$ (the latter is $(\phi_t(\fp), [f\circ \phi_t^{-1}])$ can be written as 
\begin{equation}
\mathcal{E}(\bar\phi_t(\fp,[f])) = \int_{S^2}\langle df, df\rangle_{\fm_t}
\tag{6.7}
\end{equation}
where $\langle , \rangle_{\fm_t}$  denotes the $\phi_t$-pull back of the round metric's inner product on $T^*S^2$.  

An important conclusion from (6.7) is that the map $t \to \mathcal{E}(\bar\phi_t(\cdot))$ is differentiable with respect to $t$, this is not-withstanding the fact that $\bar\phi_t$ is only continuous.  Moreover, the $t$-derivative of $\mathcal{E}(\bar\phi_t(\fp,[f]))$ can be written using (6.7) as
\begin{equation}
\int_{S^2} \langle \mathfrak{k}_t, df \otimes df \rangle_{\fm_t}
\tag{6.8}
\end{equation}
where $\mathfrak{k}_t$ denotes the Lie derivative of $\fm_t$ along the vector field on $S^2$ that is obtained by first viewing $\phi_{(\cdot)}$ as a map from $\real\times S^2$ to $S^2$ and then using this map to push forward the tangent vector to the $\real$ factor.  

The next lemma supplies a key observation about (6.7) and (6.8).  \\

\noindent \bf Lemma 6.5: \it  Suppose that $\delta$ is a positive number and that the $t = 0$ and $(\fp,[f]) \in\rpb$ version of (6.8) is more negative than $-\delta$.  Then there is an open neighborhood of $(\fp,[f])$ in $\rpb$ (denoted by $\ub $) and a positive number (denoted by $\epsilon$) with the following significance:  If $\xi \in \ub $ and $t\in(-\epsilon,\epsilon)$, then $\mathcal{E}(\bar\phi_t(\xi)) \leq \mathcal{E}(\xi)- {1\over 4}\delta t$.\\

\noindent \bf Proof of Lemma 6.5:  \rm The first observation for the proof concerns (6.8):  If this is more negative than $-\delta$ at $(\fp,[f])$, then there is an open neighborhood of $(\fp,[f])$ in $\rpb$ (this will be the neighborhood $\ub $) such that the $t = 0$ and $(\fp^{\prime},[f^\prime])$ version of (6.8) is more negative than $- {1\over 2}\delta$ when $(\fp^{\prime}, [f^\prime])$ is in $\ub $.  Except for notation, the proof of that this is so is identical to the proof of Lemma 5.1.  Then, looking at the $t$-dependence, it follows from the continuity of $\mathfrak{k}_{(\cdot)}$ that there exists $\epsilon>0$ such that any $t \in[0, \epsilon)$ and $(\fp^{\prime},[f^\prime]) \in \ub $ version of (6.8) is more negative than $-{1\over 4} \delta$.  The assertion in the lemma follows directly from this last observation by integrating (6.8) with respect to $t$.\\ 

\it Part 7:  \rm Take $(\fp,[f])$ and $\delta$ and then $\ub $ and $\epsilon$ from Lemma 6.5.  Suppose that $\Psi$  is from the set $\fcc$.  Let $\tau_*(\cdot)$ denote a continuous function on $\yb$ mapping to the interval [0,1] with compact support on $\Psi^{-1}(\ub $).  A path of deformations of $\Psi$  parametrized by $[0,\epsilon]$ is defined by the rule whereby $(t,y)\in[0,\epsilon]\times\overline{\mathcal{Y}}$ is sent to $\bar\phi_{\tau_*(y)t}(\Psi (y))$.  The $t \in[0, \epsilon]$ member of this deformation family is denoted by $\Psi_{\tau_*,t}$.  The behavior of $\mathcal{E}$ with respect to these deformations is the topic of the next paragraphs.

	With regards to $\mathcal{E}$:  If $y$ is a point from $\yb$ that is not mapped by $\Psi$  to $\ub $, then $\Psi_{\tau_*,t}(y)$ is $\Psi(y)$ because $\tau_*(y)$ is zero; and thus $\mathcal{E}(\Psi_{\tau_*,t(y)}) =\mathcal{E}(\Psi (y))$.  With this understood, suppose now that $y$ is mapped by $\Psi$  to $\ub $.  In this case, Lemma 6.5 finds 
\begin{equation}
\mathcal{E}(\Psi_{\tau_*,t}(y)) \leq \mathcal{E}(y)-  {1\over 4}\tau_*(y) t \delta .
\tag{6.9}
\end{equation}

This implies (among other things) that the maximum of $\mathcal{E}$ on the image of $\Psi_{\tau_*,t}$ is no greater than its maximum on the image of $\Psi$.\\

\it	Part 8: \rm  Let  $\mathfrak{C}_{n,m}$ denote the subset in $\rpb$ consisting of elements $(\fp,[f])$ with $f$ being an $\ip$-eigensection with eigenvalue $\mathcal{E}_{n,m}$.  A key point to keep in mind for what follows that  $\mathfrak{C}_{n,m}$ is a closed, compact set.  \\
	
	Suppose for the sake of argument that there is not one $k \in\{0, 1, \ldots, n\}$ such that the intersection of  $\mathfrak{C}_{n,m}$ with the stratum $\mathcal{C}_{2k}$ contains a critical point of $\mathcal{E}$ on  $\mathcal{C}_{2k}$ .  The four steps that follow use this assumption to generate nonsense.\\
 
\underline{Step 1}:  Suppose now that $k \in\{0, 1, \ldots, n\}$ and that $\fp$ is in the stratum $\mathcal{C}_{2k}$ of the compactification  $\overline{\mathcal{C}}_{2k}$ .  Let $v$ denote a tangent vector to $\fp$ along this stratum (this being a vector in $\otimes_{p\in\fp} TS^2|_{p}$).  Section 2.6 describes how to extend $v$ from the points in $\fp$ as a divergence free vector field on $S^2$.  Denote this vector field by $\mathfrak{v}$.  Integrating this vector field gives a 1-parameter family of area preserving diffeomorphisms of $S^2$.  The family is parameterized by the coordinate $t$ on the real line $\real$ with the $t = 0$ member being the identity diffeomorphism and with the $t$-derivative at $t = 0$ being $\mathfrak{v}$.  The diffeomorphism labeled by any given $t \in\real$ is denoted by $\phi_t$. \\

\underline{Step 2}:  Because of the assumptions, for each $k \in\{0, 1, \ldots, n\}$ the following is true:  If $(\fp,[f])$ is from  $\mathfrak{C}_{n,m} \cap \rp_{2k}$, then there is a corresponding $\phi_t$ from Step 1 for which (6.8) is negative.  Therefore, each $(\fp,[f])$ in  $\mathfrak{C}_{n,m}$ has a corresponding $\epsilon, \delta$ and $\ub$ as described in Lemma 6.5.  Since  $\mathfrak{C}_{n,m}$ is compact, there exists positive $\epsilon_*$ and $\delta_*$ such that each $(\fp,[f])\in \mathfrak{C}_{n,m}$ version of $\epsilon$ and $\delta$ is greater than the respective $\epsilon_*$ and $\delta_*$.  In addition, there exists a finite set of elements in  $\mathfrak{C}_{n,m}$whose corresponding versions of $\ub$ supply an open cover of  $\mathfrak{C}_{n,m}$.  Let $\Lambda$ denote this finite set of elements and let $\overline{\mathcal{U}}$ denote the union of these $\Lambda$-versions of $\ub$.  \\

 \underline{Step 3}:  By virtue of what is said in Part 4, there exists some positive number to be denoted by $\epsilon_{n,m}$ with the following significance:  If $\Psi$  is from $\fcc$ and if the maximum of $\mathcal{E}$ on $\Psi$  is less than $\mathcal{E}_{n,m}+\epsilon_{n,m}$, then the part of the image of $\Psi$  where $\mathcal{E} > \mathcal{E}_{n,m}-\epsilon_{n,m}$ is contained in $\overline{\mathcal{U}}$.  With this understood, suppose henceforth that the maximum of $\mathcal{E}$ on the given version of $\Psi$  is indeed less than $\mathcal{E}_{n,m}+\epsilon_{n,m}$.
 
Let $\tau$ denote the function on $\yb$ given by the rule
\begin{equation}
\tau_0(\cdot) \to \chi\left( {2\over \epsilon_{n,m}} \left( \mathcal{E}_{n,m}-  {1\over 2} \epsilon_{n,m}- E(\Psi (\cdot))\right)\right).
\tag{6.10}
\end{equation}
This function is equal to 1 on the part of  $\overline{\mathcal{Y}}$ where $\mathcal{E}(\Psi(\cdot)) \geq \mathcal{E}_{n,m} -  {1\over 2}\epsilon_{n,m}$ and it is equal to 0 on the part where $\mathcal{E}(\Psi(\cdot)) < \mathcal{E}_{n,m} - \epsilon_{n,m}$.  In particular, this implies that the support of $\tau$ is in $\bar U$.\\

\underline{Step 4}:  Label the elements in $\Lambda$ consecutively starting from 1 and use the corresponding labels for the $\Lambda$-versions of  $\bar U$.  Fix for each label $i$ from $\Lambda$, fix an open set $V_i$ in $U_i$ having compact closure in $U_i$ and such that the union of these sets also covers $\mathfrak{C}_{n,m}$.  Then, for each index $i$ from $\Lambda$, fix a continuous function with values between 0 and 1 to be denoted by $\sigma_i$ which is compactly supported on $U_i$ and is equal to 1 on $V_i$.
  
Apply the homotopy from Part 7 to $\Psi$ using  $\bar U_1$ and with $\tau_*$ being $\sigma_1 \tau_0$ with $\tau_0$ defined via (6.10).  Let $\Psi_1$ denote the $t = 1$ end member of this homotopy.  Now apply the homotopy from Part 7 again using $\bar U_2$ starting from $\Psi_1$ and with $\tau_*$ being $\sigma_2 \tau_0$.  Continue in this vein until all of the $\Lambda$-versions of  $\bar U$ are accounted for.  Doing that results in a map from $\mathfrak{F}$ (denoted by $\Psi_*$) that has non-zero pairing with $\tau^{2m} \wedge \Xi$.  

By virtue or (6.9), this map $\Psi_*$ obeys
%
\begin{equation}
\max_{\yb} \mathcal{E}(\Psi_*(\cdot)) \leq \max_{\yb} (\Psi (\cdot)) - c_0^{-1} \epsilon_* \delta_*.
\tag{6.11}
\end{equation}
This result is the desired nonsense because if maximum of $\mathcal{E}(\Psi (\cdot))$ is less than $ \mathcal{E}_{n,m} + {1\over 2}c_0^{-1}\epsilon_*\delta_*$, then the maximum of  $\mathcal{E}$  on  $\Psi_*$ will be less than  $\mathcal{E}_{n,m}$  which can't happen by virtue of  $\mathcal{E}_{n,m}$  being the infimum of the set of maxima of  $\mathcal{E}$  on the maps from $\fcc$.

\subsection{On the existence of critical points of  $\mathcal{E}$  on $\ctwon$}

	The analysis so far allows for the following pathology:  There exists a positive integer (call it $j$) such that all of the critical point of  $\mathcal{E}$  that are obtained via Proposition 6.3 using every possible pair of integers  $n$  and $m$  (with  $n  \geq  1$ and  $m  \geq  0$) lie in $\overline{\mathcal{C}}_{2j}$ (which is a subspace of all  $n  \geq  j$ versions of $\ctwonb$).   We can't rule this out at present.  In any event, the question directs us to derive a more refined picture of the behavior of  $\mathcal{E}$  on  $\ctwon$  near the lower dimensional strata in its closure  $\ctwonb$.  The upcoming Lemma 6.6 is an example of what might be needed.

	To set the stage for Lemma 6.6, suppose that $j \in\{0, \ldots, n-1\}$ and that  $\fq$  is the minimal representative of a  $\integer/2\integer $ divisor in the $\mathcal{C}†_{2j}$ stratum of  $\ctwonb$.  Let $\{p_1, \ldots, p_{2n -2j}\}$ denote distinct points in  $S^2-\fq$  and set  $\fp  \equiv \{p_1, \ldots, p_{2n -2j}\}\cup \fq$.  Fix  $n$  disjoint, embedded arcs in  $S^2-\fq$  whose endpoints pair the points in the set $\{p_1, \ldots, p_{2n -2j}\}$.  Let  $\Sigma$  denote their union.  Then choose a point to be denoted by  $x$  in  $S^2 -( \Sigma\cup \fq).$  Since the bundle $\iq$ and  $\ip$  are isomorphic on  $S^2 -( \Sigma \cup \fq)$, the choice of an isomorphism of their fibers at  $x$  extends uniquely to an isomorphism between the two line bundle over their common domain  $S^2 -( \Sigma\cup \fq)$.  A choice of such a point  $x$  and a corresponding isomorphism will be made momentarily.
	
	To continue with the stage setting:  Let $f_{\fq}$ denote an $\iq$-eigensection and let $f_{\fp}$ denote an  $\ip$-eigensection; and let $\mathrm{E}_{\fq}$ and $\ep$  denote their respective eigenvalues.  Take  $x$  to be a point where neither $f_{\fq}$ nor $f_{\fp}$ is zero and then choose the isomorphism between $\iq$ and  $\ip$  at  $x$  so that $f_{\fq}f_{\fp}$ is positive at $x$.  (The product $f_{\fq}f_{\fp}$ is a section of the line bundle $\ipp$ with  $\fp^\prime  \equiv \{p_1, \ldots, p_{2n -2j}\}$.  The isomorphism between  $\ip$  and $\iq$ on  $S^2 -( \Sigma\cup \fq)$ gives an isomorphism between $\ipp$ and the product real line bundle on  $S^2 - \Sigma$.)
	
Now comes a crucial observation:  The preceding isomorphism between $\iq$ and $\ip$  on  $S^2 -( \Sigma \cup \fq)$ defines two isomorphisms along  $\Sigma$  which differ via multiplication by $-1$.  These are obtained by following the isomorphism defined in the preceding paragraph along a path from  $x$  to  $\Sigma$  that starts at  $x$  and ends on  $\Sigma$  (transversely to  $\Sigma$) whose interior is disjoint from  $\Sigma$.  In particular, a choice of normal vector to  $\Sigma$'s arcs defines one of these two isomorphisms between $\ip$  and $\iq$ along  $\Sigma$:  Take a path from  $x$  that hits  $\Sigma$  at its endpoint from the side of the chosen normal vector.  Choose this normal vector to define the isomorphism along  $\Sigma$.   With regards to notation, Lemma 6.6 has $\partial_n$ denoting the directional derivative on  $\Sigma$  in the chosen normal direction.   \\

\noindent \bf Lemma 6.6: \it Suppose that $\fq,  \fp$  and $f_{\fq}, f_{\fp}$, and $\mathrm{E}_{\fq}, \mathrm{E}_{\fp}$, and also the isomorphism between  $\ip$  and $\iq$ are as described above.  Then there is the following identity:
\begin{equation*}
(\ep - \eq) \int_{S^2 - (\fq \cup \Sigma)} f_{\fq} f_{\fp} = -2 \int_\Sigma (f_{\fq} \partial_n f_{\fp} - f_{\fp}\partial_n f_{\fq}).
\end{equation*}

\noindent \bf Proof of Lemma 6.6:  \rm Multiply the equation $\ep f_{\fp} = -\Delta f_{\fp}$ by $f_{\fq}$ and likewise multiply the equation $\eq f_{\fq} = -\Delta f_{\fq}$ by $f_{\fp}$ and then subtract the second from the first.  Having done that, integrate both sides of the result over  $S^2 -( \Sigma \cup \fq)$ and then integrate by parts.\\

	This lemma is most most useful if  $\Sigma$  is part of the zero locus of $f_{\fp}$ in which case there is only the $f_{\fq}\partial_nf_{\fp}$ term in the lemma's right hand side integral whose sign on any component of  $\Sigma$  can be determined by the sign of $f_{\fq}f_{\fp}$ on the chosen normal's side of  $\Sigma$.  (If this sign is $+$, then the integrand is $-$ and vice versa.)  For example, if $j = n-1$ so that there are just two points in  $\fp^\prime$  which are very close to each other (and one arc in  $\Sigma$), and if $f_{\fq}f_{\fp}$ is positive very near to  $\Sigma$  on the normal direction side of  $\Sigma$  then the lemma's right hand side integral is negative and so the right hand side in the lemma's identity is positive.
	  
The preceding observation implies the claim made in the first sentence of the third bullet of (1.5) to the effect that the supremum of the function  $\mathrm{E}_{(\cdot)}$  on  $\ctwon$  is always achieved at some configuration on $ \ctwon$.  Use an induction argument with the integer  $n$  starting from  $n  = 1$ as follows:  Supposing that the supremum of  $\mathrm{E}_{(\cdot)}$  is taken at a configuration in  $\ctwon$.  Take  $\fq$  in Lemma 6.6 to be such a configuration and take $f_{\fq}$ to be a $\iq$-eigensection with  $\iq$-eigenvalue equal to $\eq$.  Now take  $\fp \in \mathcal{C}_{2n+2}$ as follows:  Let  $x$  denote a point in  $S^2$  where $f_{\fq} \neq 0$ and let $p_1$ and $p_2$ denote two points very close to $x$, equidistant from  $x$  on opposite sides of a great circle arc through $x$.  Let $r$ denote the distance between them.  Set  $\fp  = \{p_1, p_2\}\cup \fq$ and let $f_{\fp}$ denote the  $\ip$-eigensection with eigenvalue $\ep$.  If $r$ is sufficiently small, then there is a component of the zero locus of $f_{\fp}$ in a small radius ball centered on  $x$  that connects $p_1$ to $p_2$.
 (This can be proved using Lemma 5.2.)  Take  $\Sigma$  to be such a component.  As explained in the preceding paragraph, the formula in Lemma 6.6 implies in this case that $\ep$  is greater than $\eq $.   

To elaborate on what is going on here:  The set up in the preceding paragraph is designed so that $|f_{\fp}|$ is very close to $|f_{\mathfrak{q}}|$ except in a small radius disk containing the points $p_1$ and $p_2$.  Thus, over most of $S^2 - (\mathfrak{q} \cup \{p_1, p_2\})$, the product $f_{\fp}f_{\fq}$ will be positive except very near where either $f_{\fq}$ is zero or $f_{\fp}$ is zero with it understood that the identification between $\ip\otimes \iq$ and the product line bundle on $S^2 \cup (\fq \cup \{p_1, p_2\})$ is obtained from the identification at $x$ by following paths from $x$ that don’t intersect $\fq \cup \Sigma$.  In particular, the proportion of $S^2- (\fq - \Sigma)$ where $f_{\fp} f_{\fq} > 0$ approaches 1 as the distance from $p_1$ to $p_2$ shrinks to zero.  This is why the integral of $f_{\fp}f_{\fq}$ is positive when the distance from $p_1$ to $p_2$ is sufficiently small.

 \section{The case of $\mathcal{C}_2$}

	The case where  $n  = 2$ already illustrates much of what is described in the preceding sections.  The subsections that follow have some observations on this case.

\subsection{Identifying  $\ctwo$  and $\ctwob$.}

The space  $\ctwo$  can be identified with the unit disk bundle in the tangent bundle to $\rp^2$.  To make this explicit, first view  $S^2$  as the radius 1 sphere about the origin in  $\mathbb{R}^3$, and then identify the unit disk bundle in its tangent bundle (which is denoted by $T_{<1} S^2$) as
\begin{equation}
T_{<1} S^2 = \{(n,v) \in \mathbb{R}^3 \times \mathbb{R}^3: |n| = 1, |v| < 1 \text{ and } \langle n,v\rangle  = 0\}
\tag{7.1}
\end{equation}
The projection to  $S^2$  sends any given pair $(n,v)$ to $n$.  There is an action of $\{ \pm 1 \}$ on $T_{<1} S^2$  that covers the $\{ \pm 1 \}$ action on  $S^2$  whose quotient is $\mathcal{C}_2$, the quotient identification being
\begin{equation}
(n,v) \sim (-n,v).
\tag{7.2}
\end{equation}
The quotient, $\mathcal{C}_2$, is the unit disk bundle in an $\real^2$ bundle over $\rp^2$ that is isometric to the tangent bundle of $\rp^2$.  Explicitly, the inverse identification sends the equivalence class $([n],v)$ to the pair of points
\begin{equation}
v \pm (1-|v|^2)^{1\over 2} n
\tag{7.3}
\end{equation}
on the unit sphere about the origin in  $\mathbb{R}^3$.  (Here and below, when  $n$  is a unit length vector in $\mathbb{R}^3$, then $[n]$ denotes the corresponding point in $\rp^2$.)  Note in particular that the $v  = 0$ locus is the copy of $\rp^2$ given by those configurations where the constituent points are antipodal.  With  $\ctwo$  viewed as the disk bundle in an $\real^2$ bundle over $\rp^2$, this  $v  = 0$ locus is the zero section of that vector bundle.  

As for $\ctwob$, this is the Thom space for this $\real^2$-bundle over $\rp^2$:  It is the space of pairs $(n, v) \in  \mathbb{R}^3 \times  \mathbb{R}^3$  with $|n| = 1, |v| \leq 1$ and  $\langle n,v\rangle  = 0$ modulo the equivalence relation in (7.2) and the equivalence
\begin{equation}
(n,v) \sim (n^\prime, v^\prime) \text{  when }  |v| = |v^\prime| = 1.
\tag{7.4}
\end{equation}
Thus, a path in  $\ctwo$  where $|v| \to 1$ converges in  $\ctwob$  to the stratum  $\mathcal{C}_0$ in $ \ctwob  = \mathcal{C}_2\cup \mathcal{C}_0$.

	In this case, the cohomology of  $\ctwo$  is that of $\rp^2$ which is $\integer/2\integer$ in degree 0, 1 and 2.  Meanwhile, an instance of the Thom isomorphism theorem says that the cohomology of  $\ctwob$  is  $\integer/2\integer$  in degree 0, it is zero in degree 1, then  $\integer/2\integer$  in degrees 2, 3 and 4.  The degree 2 class is the obstruction class $ \omega$, and (according to the Thom isomorphism theorem), it's self cup product is non-zero.
	
\subsection{The fiber bundle $\mathcal{F}$}
	The fiber bundle  $\mathcal{F}  \to  \ctwo$  consists of the set of pairs of the form $(\fp,x)$ with $\fp\in\ctwo$ and with  $x  \in S^2-\fp$.  The fiber of  $\mathcal{F}$  over any given configuration is diffeomorphic to an open annulus.  This annulus deformation retracts onto a circle bundle whose fiber over a given configuration  $\fp  = ([n], v)$ is the set of unit length vectors in the plane orthogonal to the vector  $n$  (this is the great circle whose points are equidistant from the two points in $\fp$).  Viewed in this light, it is the pull-back from $\rp^2$ of the unit disk bundle in the tangent bundle to $\rp^2$.   
	
The assertion in Section 4.1 to the effect that there is no real line bundle over  $\mathcal{F}$  that restricts to the fiber of each configuration  $\fp$  from $\ctwo$ as  $\ip$  is equivalent to the assertion that there is no degree 1 class in the $ \integer/2\integer$  cohomology of the unit circle bundle in $\mathrm{T}\rp^2$ that restricts non-trivially to each fiber.  In turn, the latter is equivalent to the assertion that any given fiber circle represents the zero class in the  $\integer/2\integer$  first homology of the unit circle bundle in $\mathrm{T}\rp^2$.  That this is so can be seen by constructing a section of this bundle over the complement in $\rp^2$  of any given point in $\rp^2$  and then looking at how the section behaves near the omitted point (see below).  

For example, let $n_*$ denote a given unit length vector in   $\mathbb{R}^3$   and let $[n_*]$ denote the corresponding point in $\rp^2$ .  A section over $\rp^2 -[n_*]$ is given by the rule
\begin{equation}
[n] \to \left([n],  v  = {1\over \sqrt{1-\langle n, n_*\rangle^2}} (n_* - n \langle n, n_*\rangle )\right).
\tag{7.5}
\end{equation}
If  $n$  is close to $n_*$, then  $n$  can be written as  ${1\over\sqrt{1+\epsilon^2}}(n_*+\epsilon u )$ where $u$ has norm 1 and is orthogonal to $n_*$ and where  $\epsilon  \in(0,1)$.  The section in (7.5) near  $\epsilon  = 0$ when written in terms of  $\epsilon$  and $u$ has the form

\begin{equation}
\left[{1\over \sqrt{1+\epsilon^2}} (n_*+\epsilon u )\right] \to  ([n_*],  v  =u) + \mathcal{O}(\epsilon)
\tag{7.6}
\end{equation}
which directly illustrates that the fiber over $n_*$ is a boundary in the unit tangent bundle.  Indeed, the preceding construction depicts the fiber over $[n_*]$ as the boundary of an embedded M\"obius band.  Meanwhile, the boundary of a M\"obius band is zero in the  $\integer/2\integer$  homology of the M\"obius band since it is twice the class of the central circle in the band's $\integer$-homology.

\subsection{The space  $\rpb$  and the function $\mathcal{E}$}
	As explained in the appendix of \cite{TW}, the only critical points of  $\mathcal{E}$  on the  $\ctwo$  version of  $\rp$  are over the configurations of antipodal points in $\ctwo$.   In this case, the set of  $\ip$ -eigenvalues is the set $\left\{m^2- {1\over 4}\right\}_{m=1,2 \ldots}$; any given integer $m$  eigenvalue has multiplicity $2m$.   The eigenspace with eigenvalue $m^2-{1\over 4}$ for a given antipodal configuration is spanned by eigensections that can be described as follows:  Take spherical coordinates for  $S^2$  with the longitude being rotation about the axis through the pair of points from the configuration.  With these coordinates denoted by $(\theta,\varphi)$, then the eigenspace consists of the linear combinations of the real and imaginary parts of the eigensection that is depicted below in (7.6) and then the real and imaginary parts of what is obtained from that eigensection by acting by the generators of the $SO(3)$ action no more than m-1 times.
\begin{equation}
f_m = (\sin\theta)^{m-{1\over 2}}e^{\left(m-{1\over 2}\right)\varphi}  
\tag{7.7}
\end{equation}
Note that the minimal eigenvalue ($m = 1$) has multiplicity two and it does not obey the criteria to be a critical point of $\mathcal{E}$.  The next highest eigenvalue is a critical value (the $m  = 2$ case in (7.7)).  But note also that there is a dimension 4 eigenspace with this value of $\mathcal{E}$.  All of the  $m  \geq  2 $ eigenvalues are also critical values of $\mathcal{E}$.

	The compactification $\ctwob$  adds the one point  $\mathcal{C}_0$  to $\ctwo$.   The spectrum of the Laplacian over the $\mathcal{C}_0$ stratum is the set $\{m(m+1)\}_{m=0,1,\ldots}$ with any given $m(m+1)$ eigenvalue having multiplicity $2m+1$.  This being the case,  $\mathcal{I}_{(\cdot)}$ -eigenvalues for configurations in $\ctwo$ with the two points very close will be very nearly of the $m(m+1)$ form with $m$  being a non-negative integer with its values ranging from zero to some positive integer that gets ever larger as the points come together.
		
	To elaborate:  If the two points are moved from being almost antipodal to being almost coincident along a great circle arc, then the eigenvalues from the antipodal set $\left\{m^2- {1\over 4}\right\}_{m=1,2 \ldots}$ must move in a continuous fashion to generate the coincident set $\{m(m+1)\}_{m=0,1,\ldots}$ and likewise for the eigensections.  Keeping in mind that the antipodal eigenvalues are alternating with the coincident eigenvalues on the non-negative real axis, here is how this works:  When the points are moved together along a great circle arc,  $m$  of $2m$-linearly independent eigensections for a given antipodal eigenvalue $m^2-{1\over 4}$ have decreasing eigenvalue which limit to $m(m-1)$, and  $m$  of them have increasing eigenvalue which limit to $m(m+1)$.
	  
With regards to the preceding:  The precise eigenspaces that give decreasing or increasing eigenvalues depends on the choice of the great circle arc that is used to bring the two points together.   This can be seen most clearly for the $m  = 1$ case of the antipodal configuration eigenspace (eigenvalue  $3\over 4$) whose normalized eigensections have the form 
\begin{equation}
z\sin^{1\over2}\theta \sin\left( {1\over 2}(\varphi-\alpha )\right) 
\tag{7.8}
\end{equation}
with $\alpha  \in \real/2\pi \integer$ and with  $z$  being independent of $\alpha$.  Consider a path in  $\ctwo$  that starts with an antipodal configuration and then changes the configuration so as to bring the two points together along a constant $\varphi$ half great circle of longitude (say $\varphi = \varphi_*$). It follows from Proposition 2.6 that the $\alpha  = \varphi_*$ version of (7.6) deforms to a corresponding family of $ \mathcal{I}_{(\cdot)}$-eigensections along the configuration path with the corresponding path of eigenvalues decreasing to 0 as the points converge.  Meanwhile, the $\alpha =\varphi_*+\pi$ version of (7.6) deforms to a family of  $\mathcal{I}_{(\cdot)}$-eigensections along the configuration path with the corresponding path of eigenvalues increasing to 2 along the configuration path (keep in mind that 2 is the $m=1$ version of $m(m+1)$). 

	One can imagine that there is a topological explanation for the interesting spectral flow over  $\ctwo$  via an analysis of the effect on  $\rpb$  of the rotation group's action on  $S^2$ .

\appendix
\section{The cohomology of a weak   $\rp^\infty$   bundle}

	This appendix states and then proves a proposition that describes the  $\integer/2\integer$  cohomology of a weak   $\rp^\infty$   bundle over a space $X$.   An additional lemma at the end of the appendix (Section A.6) points out a useful corollary from the proof.  (A proposition much like the one below might well be in the literature in which case we apologize to its authors for failing to find it and hence reference it.)  
	
	To set the stage for the proposition, suppose that $X$ is a compact CW complex and that $\pi: \mathcal{R} \to X $ is a weak   $\rp^\infty$   bundle.  (The assumption holds for the case $X =  \ctwonb$.)  

There is an obstruction class in $H^2(X)$ which is defined just as in Section 5.4.  Let $ \omega$  denote this obstruction class.  Now introduce $\mathcal{A}^*$ to be the kernel in $H^*(X)$ of the homomorphism  $\omega \wedge (\cdot) $ to $H^{*+2}(X)$, and let $\mathcal{B}^*$ denote the cokernel.  By way of an example, let $a$ denote the highest non-zero self-cup product of  $\omega$.  This is to say that $  \omega^a  \neq 0$ but $  \omega^{ a+1}= 0$.  Then $  \omega^a$  is in $\mathcal{A}^{2a}$.    \\

\noindent \bf Proposition A.1: \it The  $\integer/2\integer$  cohomology of $\mathcal{R}$ is isomorphic (non-canonically) to the vector space of polynomials in a degree two class  $\tau$  with coefficient ring $\mathcal{C}$ as described in the subsequent bullets.  Thus, 
$$H^*(R) \approx  C^*\oplus (C^{*-2}\wedge  \tau ) \oplus (C^{*-4}\wedge  \tau^2) \oplus \cdots$$
\begin{itemize}
\item	The definition of $ \tau$:  Let $\hat\omega$  denote a given representative cocycle for $  \omega$.  Then $\pi^*\hat\omega$  can be written as $d\nu$ with  $\nu$  being a 1-cochain whose restriction to each fiber algebraically generates the cohomology of the fiber.  The cochain  $ \nu \wedge \nu $   is closed and non-zero in $H^*(\mathcal{R}$).  The class of   $\nu \wedge \nu$   is the class  $\tau$.  This class can depend on the choice for $\hat\omega$  and  $\nu$  but any such change has the form $\pi^*\lambda$ with $\lambda$ a  non-zero degree 2-class in $H^*(X)$ that can be represented as the square of a 1-cochain (as  $\mu \wedge  \mu )$.
\item	The definition of $\mathcal{C}$:  The degree $*$ summand of $\mathcal{C}$ is zero if $*$ is negative, it is isomorphic to  $\integer/2\integer$  if $* = 0$, and for positive degrees, it is isomorphic (non-canonically) to $\mathcal{B}^*\oplus  \mathcal{A}^{*-1}$ with an isomorphism defined in the three parts that follow:  
\begin{itemize}
\item 	Fix a basis for $\mathcal{B}^*$ and then a lift of each basis element as a closed cocycle that represents a class in $H^*(X)$ which projects back to the basis element (the lift for a given basis element $\rb$ is denoted $\hat{\rb}$).  
\item Fix a basis for $\mathcal{A}^{*-1}$, and for each basis element, fix a cocycle that represents that element (the cocycle for the basis element $\mathrm{A}$ is denoted by $\hat{\mathrm{A}}$.  Having fixed $\hat{\mathrm{ A}}$, then fix degree $*$ cochain $\hat{i}_\mathrm{A}$ with the property that $d\hat{i}_\mathrm{A} =  \hat \omega\wedge  \hat {\mathrm{A}}$.  
\item The isomorphism $\mathcal{C}^* \approx \mathcal{B}^* \oplus  \mathcal{A}^{*-1}$ sends any given pair of basis elements $\rb$ and $\ra$ to the class in $H^*(\mathcal{R})$ of the cocycle  $\pi^*\hat{\rb}  +  \hat\nu\wedge \pi^*(\hat{\ra}) + \pi^*\hat{i}_A$.
\end{itemize}
\end{itemize}
\rm
The proof of this proposition is given in Section A.5 of this appendix.  The proof employs the auxiliary $\rp^{2m+1}$ bundles over $X$ that are constructed in the next subsection.  As explained in Section A.6 of this appendix, these finite dimensional fiber bundles can also be used to represent the homology of $\mathcal{R}$.

\subsection{The $\rp^{2m+1}$ bundles}
	Fix a non-negative integer $m$.  The two parts of this subsection use the obstruction cochain to construct a useful $\rp^{2m+1}$ bundle over $X$.  This $\rp^{2m+1}$ bundle is denoted by $\mathcal{Y}_m$.  \\

	\it Part 1: \rm  The class $  \omega$  on  $\ctwonb$  can be represented by a  \v{C}ech 2-cocyle as defined with respect to a finite cover.  Let $\mathcal{U}$ denote such a cover.  The corresponding (multiplicative) cocyle assigns either $-1$ or 1 to each component of each mutual intersection of three sets from $\mathfrak{U}$.  If $U^{\prime\prime}$, $U^{\prime}$ and $U$ are three sets from $\mathfrak{U}$ with points in common, then the associated number ($-1$ or 1) is denoted by $z_{_{U^{\prime\prime}U^{\prime}U}}$.   

	Now let $\{\varphi_{_U}\}_{_{U\in\mathfrak{U}}}$ denote a partition of unity subordinate to the cover $\mathfrak{U}$.  This can be used to produce maps from pairwise intersections of sets from $\mathfrak{U}$ to the Lie group $SO(2m+2)$ as follows:  Supposing that $U$ and $U^{\prime}$ are sets from $\mathfrak{U}$ that share points, set 
\begin{equation}
w_{_{U^{\prime}U}} =  {1\over 2}\sum_{U^{\prime\prime}\in \mathfrak{U}} \varphi_{_{U^{\prime\prime}}}(1-z_{_{U^{\prime\prime}U^{\prime}U}}),
\tag{A.1}
\end{equation}
which is a map from $U\cap  U^{\prime}$ to [0, 1].  The latter map is used to construct a map from $U\cap U^{\prime}$ to the group $SO(2)$:

\begin{equation}
A_{U,U^\prime} \equiv \begin{pmatrix}
\cos (\pi w_{_{U^\prime U}}) & -\sin (\pi w_{_{U^\prime U}})\\
\sin (\pi w_{_{U^\prime U}}) & \cos (\pi w_{_{U^\prime U}})
\end{pmatrix}.
\tag{A.2}
\end{equation}

This is the desired map in the case when $m  = 0$.  If  $m  > 0$, the desired map is the tridiagonal, $(2m+2)\times(2m+2)$ matrix that is depicted below in (A.3) as a diagonal $(m+1)\times(m+1)$ matrix whose entries are $2 \times 2$ matrices with zero off of the diagonal and $A_{U^{\prime}U}$ on the diagonal:
\begin{equation}
\begin{pmatrix}
A_{U^\prime U} & 0 & 0\\
0 & \ddots & 0\\
0 & 0 & A_{U^\prime U}\end{pmatrix}
\tag{A.3}
\end{equation}
The latter matrix is denoted also by $A_{U^{\prime}U}.$  The important point now is that if $U^{\prime\prime}$, $U^{\prime}$ and $U$ have points in common, then 
\begin{equation}
A_{U^{\prime\prime}U^{\prime}}A_{U^{\prime}U}A_{UU^{\prime}} = Z_{U^{\prime\prime}U^{\prime}U}\mathbb{I},
\tag{A.4}
\end{equation}
where $\mathbb{I}$ denotes the $(2m+2)\times(2m+2)$ identity matrix.\\

\it	Part 2:  \rm The plan for this part of the proof is to use the constructions from Part 1 to build an $\rp^{2m+1}$ fiber bundle over $ \ctwonb$  that has two key properties:  It is not the $\{1, -1\}$ quotient of a sphere bundle in a $2m+2$ dimensional vector bundle; and the corresponding 2-dimensional obstruction class is the class $  \omega$ .  (The statement that this $\rp^{2m+1}$ bundle is not the quotient of a sphere bundle is saying in effect that there is no class in the first cohomology of the total space that restricts to the fiber over any given point in $X$ as the generator of first cohomology of that fiber.  Equivalently:  There is no real line bundle over the total space that restricts to any given fiber as the one non-trivial line bundle over the fiber.)

The desired fiber bundle is obtained from the disjoint union of the sets $\{U \times \rp^{2m+1}\}_{U\in\mathfrak{ U}}$  by the equivalence relation given momentarily.  To set the notation: Depict $\rp^{2m+1}$ as the quotient of the radius 1 sphere centered at the origin in $\real^{2m+2}$ by the action of multiplication by $\{1, -1\}$.  When a point in this sphere is denoted by $x$, then its image in $\rp^{2m+1}$ is denoted by $[x]$.   Now for the equivalence relation:  The relation identifies a point $(\fp, [x])$ in the space $U\times \rp^{2m+1}$ with a point $(\fp^{\prime},[x^\prime])$ in $U^{\prime} \times \rp^{2m+1}$ if and only if  $\fp  =  \fp^\prime$  and $[x] = [A_{_{U^{\prime}U}}x^\prime]$.  This is a valid equivalence relation (and hence defines a topological space) because each instance of $z_{_{U^{\prime\prime}U^{\prime}U}}$ in (A.4) is either 1 or $-1$.  

It also follows from (A.4) that the bundle $\mathcal{Y}_m$ cannot be lifted to a $(2m+1)$ dimensional sphere bundle because (A.4) identifies its obstruction class in $H^2(X)$ as the class $  \omega$.

\subsection{The pull-back of the obstruction class $  \omega$  to $\mathcal{Y}_m$}

Let $\pi_m$ denote the projection map from $\mathcal{Y}_m$ to $X$.   Fix a 2-cocycle to be denoted by  $\hat\omega$ on $X$ that gives the class $\omega$ .  Part 1 of this section proves that $\pi_m^*\hat\omega$  is exact.  Part 2 proves that if  $\hat\nu$ is a 1-cochain with $d\hat\nu  = \pi^*$, then  $\hn$ must restrict to each fiber of  $\pi$  so as to algebraically generate the cohomology of the fiber (which is $\rp^{2m+1}$).  \\

\it Part 1: \rm To prove that $\pi_m^*\hat\omega$  is exact:  View $\pi_m^*\mathcal{Y}_0$ as an $\rp^1$ bundle over $\mathcal{Y}_m$.  It's obstruction class is $\pi_m^* \omega$.  On the other hand, $\pi^*_m \mathcal{Y}_0$ can be viewed as the fiber product $\mathcal{Y}_m \times_{_X}\mathcal{Y}_0$ with the convention being that the the left hand $\mathcal{Y}_m$ is the `base' manifold for the fiber bundle $\pi_m^*\mathcal{Y}_0$.  Let $\pi_L$ denote the projection to this left hand $\mathcal{Y}_m$ and let $\pi_{_R}$ denote the projection to the right hand $\mathcal{Y}_0$.  Viewed as an $\rp^{2m+1}\times \rp^1$ bundle over $X$, there are two line bundles over any given fiber, $\mathcal{I}_L$ and $\mathcal{I}_R$ with the former being the tautological line bundle over $\rp^{2m+1}$ and the latter being the tautological line bundle over $\rp^1$.  The obstruction cocycle for the tensor product bundle $\mathcal{I}_L\otimes\mathcal{I}_R$ to extend over $\mathcal{Y}_m\times_{_X}\mathcal{Y}_0$ is (multiplicatively) given on triple intersections of sets $U^{\prime\prime}, U^{\prime}$ and $U$ from the cover $\mathfrak{U}$ by the  \v{C}ech 2-cocycle $z^2_{_{U^{\prime\prime}U^{\prime}U}}$, which is 1.  As a consequence of this, there is a line bundle (denoted by $\mathcal{I}_{LR}$) on $\mathcal{Y}_m\times_{_X}\mathcal{Y}_0$ that restricts to each fiber over $X$ as $\mathcal{I}_L\otimes\mathcal{I}_R$.  With $\mathcal{Y}_m\times_{_X}\mathcal{Y}_0$ now viewed as $\pi_m^*\mathcal{Y}_0$, this line bundle $\mathcal{I}_{L,R}$ restricts to each fiber $\rp^1$ of the projection $\pi_{_L}$ to $\mathcal{Y}_m$ as the tautological line bundle over the fiber.  Since the obstruction for this is $\pi_m^* \omega$, the latter class is zero in $H^2(\mathcal{Y}_m)$.   This implies that $\pi_m^*\hat \omega$  can be written as $d\hat\nu$.\\

Part 2:  This part of the proof explains why  $\hat\nu$ restricts to each fiber $\rp^{2m+1}$ as the generator of the first cohomology of the fiber.  To do this, consider first the  $m  = 0$ case.   This case is simpler because $\rp^1$ is diffeomorphic to $S^1$ and thus $\mathcal{Y}_0$ is a circle bundle over $X$.  As such, there is the Gysin sequence (see e.g. Hatcher's book \cite{Ha}):
\begin{equation}
0 \to H^1(X) 
\stackrel{\pi_0^*}{\longrightarrow}  H^1(\mathcal{Y}_0) 
\stackrel{\iota_0}{\longrightarrow}   H^0(X)
\stackrel{\hat e\wedge(\cdot)}{\longrightarrow}   H^2(X)  
\stackrel{\pi_0^*}{\longrightarrow} H^2(\mathcal{Y}_0)\to \cdots
\tag{A.5}
\end{equation}
where $\hat e$ denotes the  $\integer/2\integer$  Euler class and $ \hat e \wedge (\cdot)$ denotes the cup product with this class.  According to (A.5), the homomorphism $\iota_0: H^1(\mathcal{Y}_0)\to H^0(X)$ is surjective if and only if there is a class in $H^1(\mathcal{Y}_0)$ whose pairing with the generator of the first homology of each fiber is non-zero.  In this case, there is no such generator because the obstruction class $  \omega$  is zero.  Therefore,  $\hat e \wedge  (\cdot)$ maps $H^0(X)$ isomorphically to $H^2(X)$ and the class  $\hat e$  is the unique class in $H^2(X)$ whose pull-back to $\mathcal{Y}_0$ is zero.  Since the pull-back of the class $  \omega$  to $\mathcal{Y}_0$ is zero, it follows that $  \omega$  must be the Euler class  $\hat e$.  

	With  $\hat e$  understood to be $\omega$ , fix  $x  \in X$ and consider the following diagram:

\begin{center}
\begin{tikzpicture}[-stealth,
  label/.style = { font=\footnotesize }]
  \matrix (m)
    [
      matrix of math nodes,
      row sep    = 2em,
      column sep = 2em
    ]
    {
      H^1(\mathcal{Y}_0|_x) &H^2(\mathcal{Y}_0, \mathcal{Y}_0|_x) & H^2(\mathcal{Y}_0) & \textcolor{white}{b} \\
      0 &H^2(X,x) & H^2(X) & 0\\
      0&H^0(X,x) & H^0(X) & H^0(x)  \\
    };
  \foreach \i in {1,...,3} {
    \path
      let \n1 = { int(\i+1) } in
        (m-1-\i) edge node [above] {} (m-1-\n1)
        (m-2-\i) edge node [below, label] {} (m-2-\n1)
        (m-3-\i) edge node [below, label] {} (m-3-\n1)
        (m-2-\i) edge node [left,  label] {} (m-1-\i)
        (m-3-\i) edge node [left,  label] {} (m-2-\i);
  }
 \end{tikzpicture}\end{center}
\hfill{(A.6)}

The top row is from the exact sequence of the pair $(\mathcal{Y}_0, \mathcal{Y}_0|_x)$ and the second and third rows are from the exact sequence of the pair $(X, x)$.  The vertical arrows are the Gysin sequence homomorphisms; the top row of vertical arrows are $\pi_0^*$; and the lower vertical arrows are cup product with the Euler class (which is $  \omega$  in our case) or the relative Euler class as the case may be.

Now observe that the middle row says that $  \omega$  which is in $ H^2(X) $ comes from a unique class in $H^2(X,x)$ which will be denoted by $  \omega_x$.  Since $\pi_0^* \omega  = 0$, it follows that $\pi_0^* \omega_x$ is in the image of the connecting homomorphism from $H^1(\mathcal{Y}_0|_x)$.  In particular, this implies the claim about  $\hat\nu$'s restriction to the fiber if it is the case that $\pi_0^* \omega_x$ is non-zero.  But the fact that $\pi_0^* \omega_x$ is not zero follows from the fact that $H^0(X,x)$ is zero since the vertical arrows make for an exact sequence. 

The assertion about $\hat \nu$  for the case when  $m  > 0$ follows directly from the  $m  = 0$ case because of the block diagonal form of the matrix in (A.3).  Indeed, by virtue of (A.3) being block diagonal, the fiber bundle $\mathcal{Y}_0$ sits as a subbundle in $\mathcal{Y}_m$ whose fiber over any given point is a standard $\rp^1$ in the $\rp^{2m+1}$ fiber of $\mathcal{Y}_m$.

\subsection{The cohomology of $\mathcal{Y}_m$}


This section first states and then proves a lemma that describes the $\mathbb{Z}/2\mathbb{Z}$ cohomology of $\mathcal{Y}_m$.  The notation uses $\mathcal{A}^*$ to denote the kernel in $H^*(X)$ of the homorphism $\omega\wedge$ and $\mathcal{B}^*$ denotes the cokernel of this homomorphism.  The lemma that follows uses $\mathcal{C}^*$ to denote $\mathcal{B}^* \oplus \mathcal{A}^{*-1}$. A convention in what follows is that all modules with negative degrees are equal to zero.\\

\noindent\bf
Lemma A.2:  \it If  $m  = 0$:  There is a non-canonical isomorphism between $H^*(\mathcal{Y}_0)$ and $\mathcal{C}^*$.  If $ m  \geq  1$:  There is a class  $\tau  \in H^2(\mathcal{Y}_m)$ with  $\tau^{ 2k}$ non-zero when $k\in\{1, \ldots,m\}$ (its pairing is 1 with the generator of the $2k$'th homology of any given fiber); and there is a (non-canonical) direct sum decomposition
$$H^*(\mathcal{Y}_m) \approx \mathcal{C}^* \oplus (C^{*-2}\wedge  \tau ) \oplus (C^{*-4}\wedge  \tau^2) \oplus \cdots.$$
with it understood that $\mathcal{C}^* \equiv 0$ when $*$ is negative, that $\mathcal{C}^0 = \mathcal{B}^0$ and that the highest power of  $\tau$  that can appear is $ \tau^m$.\\

\noindent \bf Proof of Lemma A.2: \rm The proof has four parts.\\

\it	Part 1: \rm  The assertion of the lemma for the case of $\mathcal{Y}_0$ follows from the Gysin sequence in (A.5) which leads to the exact sequence
\begin{equation}
0 \to \mathcal{B}^* \to H^*(\mathcal{Y}_0) \to \mathcal{A}^{*-1}  \to 0 .
\tag{A.7}
\end{equation}
This sequence splits because the coefficients are  $\integer/2\integer$  and $H^*(\mathcal{Y}_0)
\approx \mathcal{B}^* \oplus \mathcal{A}^{*-1}$, but the splitting is not canonical.  Here is one isomorphism:  Fix a cocycle $\hat\omega$  to represent $  \omega$  
and then fix  $\hat\nu$ with $d \hat \nu = \pi_0^* \hat\omega$.  Now fix a basis for $\mathcal{B}^*$ and then a lift of each basis element as a closed cocycle that represents a class in $H^*(X)$ which projects back to the basis element (the lift for a given basis element B is denoted $\hat{\mathrm{B}}$ ).  
Meanwhile, fix a basis for $\mathcal{A}^{*-1}$, and for each basis element, fix a cocycle that represents that element (the cocycle for the basis element A is denoted by $\hat{\ra}$).  Having fixed $\hat{\ra}$, then fix degree $*$ cochains  $\{\hat{i}_A\}$ with the property that $d\hat{i}_A = \hat\omega \wedge \hat{\ra}$.  An isomorphism to split (A.7) sends any given pair of basis elements B and A to the classes in $H^*(\mathcal{Y}_0)$ of the respective cocycles
\begin{equation}
(\mathrm{B}, \mathrm{A}) \to  \pi_0^* \hat{\mathrm{B}} +  \hat \nu \wedge \pi_0^*(\hat{\mathrm{A}}) + \pi_0^*\hat{i}_A .
\tag{A.8}
\end{equation}
By way of a parenthetical remark:  The splitting in (A.7) is non-canonical because $\hat{i}_A$ can be changed to $\hat{i}_A+ \hat{\rc}$ with $\hat{\rc}$  any closed cochain.  If  $\hat{\rc}$ is not exact, and if it projects to a non-zero class in $\mathcal{B}^*$, then this change in $\hat{i}_A$ is accounted for by suitable changes to the basis for $\mathcal{B}^*$.  In addition, the cochain $\hn$  can be changed by adding $\pi_0^*\hat{\rc}$  with $\hat{\rc}$  being a closed degree 1 cochain.  If $\hat{\rc}$  is not exact and if  $\hat{\rc}\wedge \hat{\ra}$ represents a non-zero element in $\mathcal{B}^*$, then this change is also accommodated by changing the basis $\mathcal{B}^*$.

By way of a relevant example:  Let $a$ denote the positive integer such that $  \omega^a$  is non-zero and $  \omega^{a+1}$ is zero.  Fix a cocycle representative for $  \omega$  to be denoted by $\hat\omega$  and then fix a degree $2a+1$ cochain $\hat{i}$ such that $d\hat{i} =\hat\omega^{  a+1}$.  The cocycle  $\hat\nu \wedge  \pi^*\hat\omega^a + \pi^*\hat{i}$ represents a non-zero class in $H^{2a+1}(\mathcal{Y}_0)$.  Note in this regard that when $a ={1\over 2}  \dim(X)$, then the resulting class in $H^{\dim(X)+1}(\mathcal{Y}_0)$ is independent of the choice for the representative cocycles, $\hat \omega$, $\hat\nu$  and $\hat{i}$ because $H^*(X) = 0$ when $*$ is greater than $\dim(X)$.

By way of a second example:  The degree 2 cochain $\hat\nu \wedge \hat\nu$  is closed; so what is this class?  In general, the cohomology class of   $\hat\nu \wedge \hat\nu$   has the form 
\begin{equation}
\pi_0^*\hat{\rb}_0 +  \hat \nu \wedge \pi_0^*\hat{\ra}_0 + \pi_0^*\hat i_{\ra_0}  
\tag{A.9}
\end{equation}
where $\hat{\rb}_0$ is a closed form that represents a non-zero class in $\mathcal{B}^2$, where $\hat{\ra}_0$ represents a degree 1 class on $X$ that is annihilated by wedge product with $  \omega$; and where $\hat i_{\ra_0}$  is a degree 2-cochain with $d \hat i_{\ra_0} =  \hat\omega\wedge \hat{\ra}_0$.\\

\it 	Part 2:  \rm To prove the lemma when  $m  > 0$, some preliminary observations are needed.  The first is this:  If A denotes a non-zero class from $\mathcal{A}^{*-1}$ (which is the kernel in $H^{*-1}(X)$ of  cup product with $  \omega$) and $\hat{\ra}$ is a representative cocycle and $\hat{i}_A$ is a cochain that obeys $d\hat{i}_A= \hat\omega\wedge \hat{\ra}$, then  $\hat\nu_R\wedge \pi_m^*\hat{\ra} + \pi_m^*\hat{i}_A$ is a closed cocycle on $\mathcal{Y}_m$ that represents a non-zero class; it is non-zero because it is already non-zero when restricted to the sub-fiber bundle $\mathcal{Y}_0$.  For the same reason:  If $\rb$ is from $\mathcal{B}^*$ and $\hat{\rb}$ is a representative cocycle, then $\pi_m^*\hat{\rb}$  represents a non-zero class in $H^*(\mathcal{Y}_m)$.  The second observation is that the preceding constructions define a summand in $H^*(\mathcal{Y}_m)$ that is isomorphic to $\mathcal{B}^*\oplus  \mathcal{A}^{*-1}$ (although not canonically).   This summand is denoted by $\mathcal{C}^*$.
	
	Here is a third observation:  If $k$ is any integer from 1 through $m$, then $\hat\nu_R^{2k}$ is a closed cocycle which is non-zero in $H^{2k}(\mathcal{Y}_m)$.  It is closed because 2 equals 0 in  $\integer/2\integer$; and it is a non-zero class because it has non-zero pairing with the non-zero class in the degree $2k$ homology of any given fiber of the projection to $X$.   This class is not in $\mathcal{C}^*$.  Let  $\tau$  henceforth denote the class that is represented by $\hat\nu_R^2$ in $H^2(\mathcal{Y}_0)$. \\
	
	Part 3:  To prove the lemma when  $m  > 0$, return to the $S^1$ bundle $\pi_m^*: \mathcal{Y}_0 \to \mathcal{Y}_m$.  Write $\pi_m^*\hat\omega$  as $d\hat\nu_m$ and write $\pi_0^*\hat\omega$ as $d \hat\nu_0$.   Then write $\pi_L^* \hat\nu_m$  as  $\hat \nu_L$ and write $\pi_R^*\hat\nu_0$ as  $\hn_R$.  Writing cocycles additively now, note that $d( \hn_L- \hn_R) = 0$ because $\pi_L^*\pi_m^*\ho$  is the same cocycle as $\pi_R^*\pi_0^*\ho$.  This implies that the 1-cochain  $\hn_L - \hn_R$ is a closed cochain on $\mathcal{Y}_m\times_{_X}\mathcal{Y}_0$.   In this regard:  The cohomology class of  $\hn_L -  \hn_R$ represents the first Stieffel-Whitney class of the line bundle $\mathcal{I}_{L,R}$ from Part 1 of Section A.2.  (It follows from what is said in Part 2 of Section A.2 about the $\hat\nu_m$  and $\hat\nu_0$ versions of what is denoted there by $\hat\nu$ that $\hat u $ restricts to $\rp^{2m+1}\times  \rp^1$ fibers of the projection map from $\mathcal{Y}_m\times_{_X}\mathcal{Y}_0$ to $X$ as the first Stieffel-Whitney class of the tensor product of the respective $\rp^{2m+1}$ and $\rp^1$ tautological line bundles.) 
	
The cochain  $\hn_L - \hn_R$ is denoted by $\hat u $ and $u$ is used to denote its cohomology class.\\

\it	Part 4: \rm The key observation is that $\mathcal{Y}_m\times_{\real}\mathcal{Y}_0$ can be viewed on the one hand as an $\rp^{2m+1}$ bundle over $\mathcal{Y}_0$ and on the other, as a $\rp^1$ bundle of $\mathcal{Y}_m$.  In either case, the pull-backs of the powers of the class $\hat u$ generate the cohomology of the fiber.  As a consequence, the Leray-Hirsch theorem (see \cite{Ha}) can be used to depict the cohomology of $\mathcal{Y}_m\times_{_X}\mathcal{Y}_0$ in two ways which is what is done in (A.10) below. \\

\begin{itemize}
\item $H^*(\mathcal{Y}_m\times_{_X}\mathcal{Y}_0)  \approx H^*(\mathcal{Y}_0)\oplus H^{*-1}(\mathcal{Y}_0) \oplus H^{*-2}(\mathcal{Y}_0) \oplus  \cdot\cdot\cdot  \oplus H^{*-2m-1}$. 
\item $H^*(\mathcal{Y}_m\times_{_X}\mathcal{Y}_0) \approx   H^*(\mathcal{Y}_m)\oplus H^{*-1}(\mathcal{Y}_m)$.
\end{itemize}
\hfill (A.10)\\
(The notation is such that modules with negative degree are equal to 0.)  These isomorphism work as follows:  In the case of the top bullet, a class $\rc_k  \in H^{*-k}(\mathcal{Y}_0)$ is sent to the class $\pi_R^* \rc_k \wedge u^k$; and in the case of the lower bullet, a class $\rc_0\in H^{*}(\mathcal{Y}_0)$ or $\rc_1\in H^{*-1}(\mathcal{Y}_m) $ is sent to $\pi_L^*\rc_0$ or $\pi_L^*\rc_1\wedge u$ as the case may be.

	To exploit the two direct sum representations in (A.10), note first that
\begin{equation}
\hat u ^2 = \pi^*\hat{\rb}_0 +  \hn_L\wedge \pi^*\hat{\ra}_0 +  \hn_R^2  
\tag{A.11}
\end{equation}
where  $\pi$  denotes the projection from $\mathcal{Y}_m\times_{_X}\mathcal{Y}_0$ to $X$, and where $\hat{\rb}_0$ and $\hat{\ra}_0$ are as depicted in (A.9).   Note that the identity in (A.11) can also be written as  
\begin{equation}
\hat u ^2 = \pi^*\hat{\rb}_0 +  \hn_R\wedge \pi^*\hat{\ra}_0 +  \hn_R^2 + \hat u \wedge \pi^*\hat{\ra}_0
\tag{A.12}
\end{equation}

	The preceding identity can be used to write the summands that appear in the top bullet of (A.10) using the direct sum depiction in the lower bullet of (A.10) and vice-versa.  Doing that (and remembering what is said by Part 2) leads in a direct line to Lemma A.2's assertion.

\subsection{Fiber preserving maps from $\mathcal{Y}_m$ to $\mathcal{R}$}
	The purpose of this subsection is to construct a fiber preserving map from $\mathcal{Y}_m$ to $\mathcal{R}$ with certain desirable properties.  The lemma that follows summarizes.  (The lemma refers to notions from (5.5).)\\

\noindent\bf
Lemma A.3:  \it Fix a non-negative integer for $m$.  There exists a fiber preserving, fiber-wise embedding of $\mathcal{Y}_m$ into $\mathcal{R}$ with the following property:  If $k$ is a sufficiently large integer, then this embedding over any given open set $U$ from the cover $\uk$ maps $\mathcal{Y}_m|_U$ into the  $\Psi_U$-image of $(V_U-0)/\real^*$.  Moreover, the latter map is the  $\Psi_U$ image of a fiber-wise fiber bundle embedding from $\mathcal{Y}_m|_U$ to $(V_U-0)/\real^*$ that identifies $\mathcal{Y}_m|_U$ with the $\real^*$-quotient of the complement of the zero section in an $m$-dimensional sub-vector bundle in $V_U$.\\

\rm The rest of this subsection has the proof of this lemma.  It is used in the next subsection to prove Proposition A.1.\\

\bf \noindent Proof of Lemma A.3: \rm	There are two parts to the proof\\

\it Part 1: \rm  Start by fixing a large integer (to be denoted by $k$) that is in any event greater than $m$.  Let $\mathfrak{U}_k$ denote the corresponding open cover of $X$ that is described in (5.5).  Supposing that $U$ is a set from $\mathfrak{U}_k$, the first observation is that the obstruction cocycle $  \omega$  is zero upon restriction to $H^2(U)$.  This is because (5.5)'s embedding  $\Psi_U$ induces an isomorphism on the first homology of the fibers of $\mathcal{R}$ over $U$.  With this understood, it follows that $\mathcal{Y}_m|_U$ can also be written as the $\real^*$ quotient of the complement of the zero section in a $2m+2$ dimensional vector bundle over $U$.  Denote that bundle by $V_{mU}$.  An instance of the upcoming Lemma A.4 implies this:  If $k$ is sufficiently large (a lower bound is determined by  $m$  and $X$), then there is an injective vector bundle homomorphism $\varphi_{_{mU}}:V_{mU } \to V_U$.  Such a map induces then a fiberwise embedding of $\mathcal{Y}_m|_U$ into $(V_{mU}-0)/\real^*$ which composes with $ \Psi_U$ to give a fiberwise injective embedding of $\mathcal{Y}_m|_U$ into $\mathcal{R}|_U$.   \\

\noindent\bf Lemma A.4: \it  Fix a positive integer to be denoted by $m$.  There exists a positive integer (denoted by $k$) with the following significance:  If $V\to \mathcal{X}$ is a vector bundle with fiber dimension  $m$  and if $E \to \mathcal{X}$ is a vector bundle with fiber dimension $k$, then there is an injective vector bundle homomorphism from  $V$  to $E$.  Moreover, any two injective vector bundle homomorphisms are homotopic through injective bundle homomorphism.\\

\noindent\bf Proof of Lemma A.4:  \rm Let \sc{Hom}$_{m,k}$ \rm denote the vector space of linear maps from $\real^m$ to $\real^k$ (the space of matrices with  $m$  rows and $k$ columns, thus $\real^{mk}$.  Inside \sc{Hom}$_{m,k}$ \rm sits the subspace of injective maps which is denoted by \sc{Inj}$_{m,k}$. \rm The complement is a closed subvariety whose codimension is $k-m-1$.   It follows as a consequence that the homotopy groups of \sc{Inj}$_{m,k}$ \rm vanish in dimensions less than $k-m$.

With the preceding as background, the vector bundle \sc{Hom}$(V_{mU}, V_U)$ \rm is a vector bundle over $U$ whose fiber is \sc{Hom}$_{m,k}$ \rm; and inside this bundle sits the fiber bundle \sc{Inj}$(V_{mU}, V_U)$ \rm whose fiber over any given point  $x  \in U$ is the space of linear, injective maps from $V_{mU}|_x$ to $V_U|_x$ which is a copy of \sc{Inj}$_{m,k}$. \rm Remembering that $X$ is assumed to be a CW complex, there is an upper bound to the dimension of its cells.  Denote that by $z$.  If $k-m$ is greater than $z$, then an inductive construction starting from the 0-cells, then the 1-cells and so on will construct a section of \sc{Hom}$(V_{mU}, V_U)$ \rm that sits entirely in \sc{Inj}$(V_{mU}, V_U)$, \rm which is to say that it is injective on each fiber.  If $k-m > z+1$, the same sort of inductive construction can be used to prove that any two injective vector bundle homomorphisms are homotopic through a path of fiber-wise injective homomorphisms.\\

Part 2:  To continue with $\mathcal{Y}_m$ and $\mathcal{R}$:  For any choice of $k$, the obstruction cocycle restricts as zero to any given set from $\uk$.  As a consequence, the fiber bundle $\mathcal{Y}_m$ over any given such set (call it $U$) can be written as $(V_{mU}-0)/\real^*$ with $V_{mU} \to U$ denoting a vector bundle with fiber dimension $m$. 

To choose $k$:  Take $k^\prime$ sufficiently large (given  $m$  and the dimension of $X$ as a CW complex) so that if $k \geq  k^\prime$, then Lemma A.4 can be used to find an injective vector bundle homomorphism over each $U \in\uk$ from $V_{mU}$ to $V_{U}$.  With $k \gg k^\prime$, chose such an injective vector bundle homomorphism for each such $U$ and denote it by $\eta_{_U}$.
  
To patch these together:  Label the sets in $\uk$ by consecutive integers starting from 1 so that when $U_j$ and $U_{j^\prime}$ are from $\uk$ with $j^\prime > j$, then $\dim( V_{U_{j^\prime}}) \geq  \dim(V_{U_j} )$.  Supposing that $U_1$ and $U_2$  share points, then over $U_1\cap U_2$  sits the homomorphism $\eta_{_{U_2}}$ that sends $V_{mU_2}$ to $V_{U_2}$; and also sitting there is the homomorphism $T_{U_2U_1}\eta_{_{U_1}}$ that sends $V_{mU_1}$ to $V_{U_2}$.  Meanwhile, $V_{mU_2}$ and $V_{mU_1}$ are canonically isomorphic over $U_1\cap U_2$  up to the action of $\real^*$.  Having fixed such an isomorphism, then $T_{U_2U_1}\eta_{_{U_2}}$  and $\eta_{_{U_2}}$ are two injective bundle homomorphisms over $U_1\cap U_2$  from $V_{mU_2}$ to $V_{U_2}$.   With this understood, Lemma A.4 can be invoked over $U_1\cap U_2$  to construct (with the help of cut-off functions) a modification of $\eta_2$ so that the new version and $T_{U_2U_1}\eta_{_{U_1}}$ agree on $U_1\cap U_2$  up to the action of $\real^*$.  Having done this, much the same construction can be used with $\eta_3$ and Lemma A.4 to extend $\eta_1$ and $\eta_2$ from $U_1\cup U_2$  to $U_1\cup  U_2  \cup  U_3$.   Continuing in this same vein with $U_4$ and then $U_5$ and so on will construct the fiber preserving map for Lemma A.3.

\subsection{ Proof of Proposition A.1}
	The proof of Proposition A.1 has four parts.  By way of a look ahead, the basic strategy is this:  The third bullet of (5.5) says in effect that any given cohomology class on $\mathcal{R}$ is visible in sufficiently large $k$ versions of   $\displaystyle\bigcup_{U\in\uk}\Psi_U((V_{U}-0)/\real^*)$.  Meanwhile, the upcoming Lemma A.5 says in effect that if  $m$  is sufficiently large, then the pull-back by one of Lemma A.3's embeddings (with $k$ sufficiently large) will make this class visible in $\mathcal{Y}_m$.  With that understood, it then follows from Lemma A.2's description of the cohomology of $\mathcal{Y}_m$ that this class has a cocycle representative that can be written as a sum with summands of the form $\pi^*\mathbbm{q} \wedge  \tau^b$ with $\mathbbm{q}$  a cocycle that represents an element from $\mathcal{C}^*$, with $ \nu$  as described in the second bullet of Proposition A.1 and with $b$ a non-negative integer.  Conversely, any cochain on $\mathcal{R}$ of the form just described is closed and, it follows from Lemmas A.2, A.3 and A.5 that such a cochain is not exact. \\

\it	Part 1: \rm  This part of the proof introduces and then proves the afore-mentioned Lemma A.5.\\

\bf\noindent Lemma A.5:  \it Fix a non-negative integer for  m  and then a sufficiently large, positive integer k so that there is a fiberwise embedding $\varphi: \mathcal{Y}_m\to \mathcal{R}$ that factors through   $\displaystyle \bigcup_{U\in\uk} \Psi_ U((V_{U}-0)/\mathbb{R}^*)$ as described by Lemma A.3.  The homomorphism of pull-back by $\varphi$ is an isomorphism from $\displaystyle H^*\left(\bigcup_{U\in\uk}  \Psi_ U((V_{U}-0)/\mathbb{R}^*)\right)$ to $H^*(\mathcal{Y}_m)$ if the degree $*$ is $2m+1$ or less.\\

\noindent \bf Proof of Lemma A.5:  \rm The proof goes by way of an inductive Mayer-Vietoris argument that has three steps to it.  \\

\underline{Step 1}:  \rm Supposing that $U$ is a set from $\uk$, let  $\varphi_{_U}$  denote the fiberwise embedding of $\mathcal{Y}_m|_U$ into $(V_{U}-0)/\mathbb{R}^*$.   Instances of the Leray-Hirsch theorem (see \cite{Ha}) compute the cohomology of both $(V_{U}-0)/\mathbb{R}^*$ and $\mathcal{Y}_m|_U$ to have the (non-canonical) form
\begin{equation}
H^*(U)\oplus (H^{*-1} (U)\wedge \nu_{_U}) \oplus (H^{*-2}(U)\wedge \nu_{_U}^2) \oplus \cdots \oplus \nu_{_U}^b
\tag{A.13}
\end{equation}
where the notation is as follows:  First, $H^*(U)$ is identified with its pull-back via the projection to $U$, which is an injective homomorphism.  Second, $\nu_{_U}$ is a degree 1 class that restricts to each fiber as the generator of the degree 1 cohomology of that fiber, and where $b$ is the smaller of $*$ and the dimension of the fiber of $(V_{U}-0)/\mathbb{R}^*$ or $\mathcal{Y}_m$ as the case may be.   

The respective isomorphisms in (A.13) imply in particular that the homomorphism of pull-back via  $\varphi_{_U}$  identifies the respective $\mathcal{Y}_m|_U$ and $(V_{U}-0)/\mathbb{R}^*$ summands of $H^*(U)\wedge \nu^h_{_U}$ in (A.13) for values of  $h$  starting at 0 and ending at the minimum of $*$ and $(2m+1)$.  This is because the pull-back of the $(V_{U}-0)/\mathbb{R}^*$ version of $\nu_{_U}$ can be used for the $\mathcal{Y}_m|_U$ version due to how  $\varphi_{_U}$  embeds the fibers of $\mathcal{Y}_m|_U$.   In particular, the pull-back by $\varphi$ is an isomorphism on cohomology of degree $2m+1$ or less.\\

\underline{Step 2}:  Label the sets in $\uk$  consecutively from 1 in the manner of Part 2 from the proof of Lemma A.4.  For any labeling integer $j$, let $X_j$ denote the union of the sets from $\uk$  with label $j$ or greater, thus $U_j\cup  U_{j+1}\cup \cdots  \cup U_N$ where  $N$  is the largest label of the sets from $\uk$.   Supposing that $U_{j-1}$ shares points with $X_j$, let $\hat U$ denote their intersection.  The intersection of  $\Psi_{U_{j-1}} ((V_{U_{j-1}} - 0) / \real^*)$ with $\displaystyle\bigcup_{j\leq i\leq N} \Psi_{U_i} ((V_{U_{j-1}} - 0) / \real^*)$ is the restriction to $\hat U$ of the fiber bundle   $\Psi_{U_{j-1}} ((V_{U_{j-1}} - 0) / \real^*)$ because the dimension of the fiber of $V_{U_{j-1}}$  is not greater than that of any $V_{U}$ when $U$'s label is greater than $j$.  Thus, this intersection is homeomorphic to  $((V_{U_{j-1}} - 0) / \real^*)|_{\hat U}$.  As a consequence, another instance of the Leray-Hirsch theorem describes the cohomology of the intersection between the sets  $\Psi_{U_{j-1}}  ((V_{U_{j-1}} - 0) / \real^*) $ and $\displaystyle\bigcup_{j\leq i\leq N} \Psi_{U_j} ((V_{U_{j-1}} - 0) / \real^*)$  via a version of (A.13) that has $U$ replaced by $\hat U$ and has $b$ being the minimum of $*$ and the dimension of the fiber of  $(V_{U_{j-1}} - 0) / \real^*$.  Meanwhile the cohomology of $\mathcal{Y}_m|_{\hat U}$ also has the form of (A.13) with $U$ replaced with $\hat U$ and with $b$ being the minimum of $*$ and $2m+1$.  

As was the case in Step 1, this implies that pull-back via $\varphi$ identifies the respective $\mathcal{Y}_m|_{\hat U}$ and $(( V_{U_{j-1}}-0)/\mathbb{R}^*)|_{\hat U}$ summands of $H^*(\hat U)\wedge \nu_{\hat U}^b$ in the respective $\hat U$ versions of (A.13) for values of $b$ starting at 0 and ending at the minimum of $*$ and $(2m+1)$.  This last conclusion implies in particular that pull-back by $\varphi$ is an isomorphism on cohomology of degree $2m+1$ or less.\\

	\underline{Step 3}:  Consider now an induction on the integer $j$ starting from $j =  N$  and working down to $j = 1$.  The induction assumption is this:  \\

\begin{center}\it For a given integer $j \in\{2, \ldots, N\}$, the embedding $\varphi$ as a map from $\mathcal{Y}_m|_{X_j}$ to $\displaystyle \bigcup_{j\leq i\leq N} \Psi_{U_i} ((V_{U_{i}} - 0) / \real^*)$  induces via pull-back an isomorphism
 between the respective cohomology up to degree $2m+1$.
\end{center} 
\hfill (A.14) 

\rm  Note that it follows as an instance of what is said in Step 1 that this assumption holds when $j = N$.  As explained directly, it follows from Steps 1 and 2 that if (A.14) holds for a given integer $j$, then it also holds with $j$ replaced by $j-1$.   The assertion of the lemma follows if this claim is true for all $j$.

The preceding claim is proved by comparing the respective Meyer-Vietoris exact sequences for $\mathcal{Y}_m|_{X_{j-1}}$ and $\displaystyle\bigcup_{j-1\leq i\leq N} \Psi_{U_{i}} ((V_{U_{i}} - 0) / \real^*)$  which come from the decomposition 
of $X_{j-1}$ as $U_{j-1}\cup X_j$.   This comparison has the form depicted below in (A.15).  The notation has $\mathcal{Z}_j$ denoting $\displaystyle\bigcup_{j\leq i\leq N}   \Psi_{U_{i}} ((V_{U_{i}} - 0) / \real^*)$.  The upcoming (A.15) also writes the intersection between $\mathcal{Z}_{j-1}$ and $(V_{U_{j-1}} - 0) / \real^*$ as
$((V_{U_{j-1}} - 0) / \real^*)|_{\hat U}$ with $\hat U \equiv U_{j-1}\cap  X_j$.

\begin{center}
\begin{tikzpicture}[-stealth,
  label/.style = { font=\footnotesize }]
  \matrix (m)
    [
      matrix of math nodes,
      row sep    = 2em,
      column sep = 0.1em
    ]
    {
   \cdots &   H^{*-1}(\mathcal{Y}_m|_{\hat U}) &H^{*}(\mathcal{Y}_m|_{X_{j-1}}) & H^{*}(\mathcal{Y}_m|_{X_j})\oplus H^{*}(\mathcal{Y}_m|_{U_{j-1}}) &\cdots \\
      \cdots &H^*((V_{U_{j-1}} - 0) / \real^*)|_{\hat U} & H^*(\mathcal{Z}_{j-1}) & H^*(\mathcal{Z}_{j})\oplus H^*((V_{U_{j-1}} - 0) / \real^*))&\cdots\\
      };
  \foreach \i in {1,...,4} {
    \path
      let \n1 = { int(\i+1) } in
        (m-1-\i) edge node [above] {} (m-1-\n1)
        (m-2-\i) edge node [below, label] {} (m-2-\n1);
  }
    \foreach \i in {2,...,4} {
      \path
      let \n1 = { int(\i+1) } in
        (m-2-\i) edge node [left,  label] {} (m-1-\i);
  }
 \end{tikzpicture}\end{center}
\hfill{(A.15)}
\\
To explain the arrows:  The horizontal sequences are the two Mayer-Vietoris sequences, whereas the vertical arrows represent the homomorphism of pull-back by $\varphi$.
  
It follows from what is said in Step 2 that the left most vertical arrow is an isomorphism on degrees $2m+1$ or less.  Meanwhile, the induction hypothesis and what is said in Step 1 imply the same for the right most pair of arrows.  This implies that the middle arrow is also an isomorphism for cohomology degrees $ 2m+1 $ or less.\\

\it	Part 2:  \rm This part of the proof establishes the first bullet of Proposition A.1.  The first claim to be proved in this regard is that $\pi^* \omega$  is zero in the cohomology of $\mathcal{R}$.  To prove this, assume to the contrary that $\pi^* \omega$  is non-zero so as to generate nonsense.  To start:  If $\pi^* \omega  \neq 0$, then the class $\pi^* \omega$  would restrict to all sufficiently large $k$ versions of   $\displaystyle\bigcup_{U\in\uk}\Psi_U((V_{U}-0)/\mathbb{R}^*)$ as a non-zero class (this is the third bullet assumption in (5.5)).  Supposing that $ m  > 1$, and supposing that $k$ is sufficiently large, then the pull-back of this non-zero class to $\mathcal{Y}_m$ by one of Lemma A.3's fiber preserving embeddings would be non-zero (see Lemma A.5).  But that event is nonsense because the latter pull-back is the same as $  \omega$'s pull-back by the projection $\pi_m $ from $\mathcal{Y}_m$ to $X$ which is zero (see Lemma A.2).

	With the preceding understood, now fix a degree 2-cocycle on $X$ to be denoted by $\hat\omega$  to represent the class $  \omega$ .  Having chosen $\hat\omega$, then fix a degree 1 cochain on $\mathcal{R}$ (to be denoted by $\hat \nu$) with $d\hat\nu  = \hat\omega$.  The 2-cochain  $\hat\nu\wedge\hat \nu$   is closed and it is not exact because it evaluates to 1 on the generator of the second homology of each fiber of $\pi$.  
	
The class of $\hat\nu \wedge \hat \nu$  is denoted by  $\tau$.  If a coboundary, say $d \mu$, is added to $\hat\omega$, then this is accommodated by changing $\hat \nu$  to $\hat\nu^\prime = \hat \nu +\pi^* \mu$.  This in turn changes $ \tau$  when $ \mu \wedge  \mu$  is a non-zero class in $ H^2(X) $ because $\hat\nu^\prime\wedge\hat\nu^\prime = \hat\nu \wedge \hat\nu  +\pi^*( \mu \wedge  \mu )$.   In general,  $\hat\nu$ can be changed by adding a closed 1-cocycle (call it $\sigma$) which changes $\hat \nu \wedge \hat\nu$  to   $\hat \nu \wedge \hat\nu + \sigma\wedge \sigma$.  But, as explained in the next paragraph, the cohomology class of $\sigma$  must be the pull-back via $ \pi$  of a class from $H^1(X)$ because these are the only 1 dimensional classes in $H^1(\mathcal{R})$; and thus the class of $\sigma \wedge \sigma$  is also pulled back by $\pi$.
   
To prove that the class of $\sigma$  is pulled back by $\pi$:  Any class in $H^1(\mathcal{R})$ is visible in $\displaystyle H^1\left( \bigcup_{U\in\uk} \Psi_U((V_{U}-0)/\mathbb{R}^*)\right)$ when $k$ is sufficiently large by virtue of the assumption in the third bullet of (5.5).  Therefore, such a class is visible in $H^1(\mathcal{Y}_m)$ for  $m  \geq  1$ when $k$ is sufficiently large (given $m$) by virtue of Lemma A.5.  And, according to Lemma A.2, all such classes are the pull-backs from $H^1(X)$ via the projection map from $\mathcal{Y}_m$ to $X$.\\

\it	Part 3: \rm  This part proves that any class in $H^*(\mathcal{R})$ can be written as a sum of those that can be represented by cocycles that have the form
\begin{equation}
(\pi^* \hat{\rb} + \hat\nu  \wedge \pi^*(\hat A) + \pi^*\hat{i}_A )\wedge   \hat\nu^b
\tag{A.16}
\end{equation}
where the notation is as follows:  What is denoted here by $\hat\nu$  is the 1-cochain from Part 1 that obeys $d\hat \nu = \hat\omega$.  Meanwhile, $\hat{\rb}$, $\hat{\ra}$ and $\hat{i}_A$ are cochains on $X$.  In particular, what is denoted by $\hat{\rb}$  is closed and represents a class in $H^*(X)$ that projects to a non-zero class in $\mathcal{B}^*$; what is denoted by $\hat{\ra}$ is a cocycle on $X$ that represents a cass in $H^{*-1} (X)$ that is annihilated by cup product with $  \omega$ ; and what is denoted by $\hat{i}_A$ is a degree $*$ cochain obeying $d\hat{i}_A = \hat{\ra}\wedge\hat\omega$.

To prove the preceding claim, let H denote a given cohomology class on $\mathcal{R}$.  If $k$ is sufficiently large, then by virtue of the third bullet assumption in (5.5), this class is non-zero in all sufficiently large $k$ versions of $\displaystyle\bigcup_{U\in\uk}  \Psi_U((V_{U}-0)/\mathbb{R}^*)$.  Thus, its pull-back by one of Lemma A.3's maps is non-zero on $\mathcal{Y}_m$ for  $m$  sufficiently (by virtue of Lemma A.5).  It then follows from Lemma A.2 that this class can be represented as a sum of classes that have the form depicted in (A.16) with   $\hn$ denoting the pull-back by the Lemma A.3 embedding of its namesake on $\mathcal{R}$.  (Note in this regard that the pull-back of a cochain on $X$ by the projection from $\mathcal{Y}_m$ to $X$ is identical to the pull-back of that cochain by the composition of first the projection from $\mathcal{R}$ to $X$ and then Lemma A.3's embedding.)\\

	Part 4:  This part proves that any class on $\mathcal{R}$ that is represented by a cocycle that can be written in the form of (A.16) with either $\hat{\rb}$'s class in $\mathcal{B}^*$ non-zero or $\hat{\ra}$'s class in $\mathcal{A}^{*-1}$ non-zero (or both) represents a non-zero class in $H^*(\mathcal{R})$.  To do this, fix a pair $(\rb, \ra)$ in $\mathcal{C}^* = \mathcal{B}^*\oplus \mathcal{A}^{*-1}$ and construct the cocycle depicted in (A.16).  The corresponding cohomology class is non-zero in $H^*(\mathcal{R})$ because its pull-back to $H^*(\mathcal{Y}_m)$ for  $m$  sufficiently large via one of Lemma A.3's fiber preserving embeddings is non-zero (take $k$ very large in Lemma A.3 and then see Lemma A.5).

\subsection{  Representing homology in  $\mathcal{R}$ }
	The lemma that follows states some implications from the preceding proof of Proposition A.1.  \\
	
	\noindent\bf Lemma A.6:  \it Fix a positive integer to be denoted by $m$.  There is a fiber preserving embedding from $\mathcal{Y}_m$ to $\mathcal{R}$ that induces (via pull-back) an isomorphism between $H^*(\mathcal{Y}_m)$  and the (finite) summand in $H^*(\mathcal{R})$ consisting of the polynomials in  $\tau$  of degree at most  $m$  with $\mathcal{C}^*$ coefficients.\\

\rm Turning this around, this lemma asserts in effect that the homology of $\mathcal{R}$ that is dual to the summand in question is obtained via push-forward from the homology of $\mathcal{Y}_m$ via a fiber preserving embedding of $\mathcal{Y}_m$ into $\mathcal{R}$.\\

\noindent\bf
Proof of Lemma A.6:  \rm Let $\mathcal{K}_m \subset H^*(\mathcal{R})$ denote the summand consisting of the polynomials in  $\tau$ of degree at most  $m$  with $\mathcal{C}^*$ coefficients.  Note in particular (from Lemma A.2) that this summand \it is \rm the cohomology of $\mathcal{Y}_m$. 
 By virtue of the third bullet assumption in (5.5), this summand $\mathcal{K}_m$ is visible in all sufficiently large $k$ versions of the subspace   $\displaystyle \bigcup_{U\in\uk}\Psi_U((V_{U}-0)/\mathbb{R}^*)$. 
  It then follows from Lemmas A.5 that if $m^\prime$ is sufficiently large, then pull-back to $\mathcal{Y}_{m^\prime}$ of $\mathcal{K}_m$ via sufficiently large $k$ versions of Lemma A.3's fiber preserving embeddings of $\mathcal{Y}_{m^\prime}$ into $\mathcal{R}$ defines an injective homomorphism from $\mathcal{K}_m$ into the cohomology of $\mathcal{Y}_{m^\prime}$.  
  Because these maps are fiber preserving, the image of $\mathcal{K}_m$ is the summand in the cohomology of $\mathcal{Y}_{m^\prime}$ given by the polynomials in  $\tau$  of degree at most  $m$  with $\mathcal{C}^*$  coefficients.  With this understood, note that if $m^\prime > m$, then there is a canonical fiber preserving embedding of $\mathcal{Y}_m$ into $\mathcal{Y}_{m^\prime}$ which is fiberwise the standard embedding of $\rp^{2m+1}$ in $\rp^{2m^\prime+1}$.  This follows from the fact that the matrix in (A.3) is block diagonal.  As a consequence, this embedding pulls back the $\mathcal{K}_m$ summand of the cohomology of $\mathcal{Y}_{m^\prime}$ isomorphically onto the cohomology of $\mathcal{Y}_m$ to give Lemma A.2's identification of this cohomology with $\mathcal{K}_m$.

\begin{appendices}

\end{appendices}


\end{document}